\documentclass[twoside]{irmaems}
\usepackage{amsmath, amssymb, latexsym, eucal, xypic, makeidx}

\ifx\pdftexversion\undefined
  \usepackage[dvips]{graphics}
\else
  \usepackage[pdftex]{graphics}
\fi

\renewcommand\theequation{\arabic{section}.\arabic{equation}}

\theoremstyle{definition} 

 \newtheorem{definition}{Definition}[section]

 \newtheorem{remark}[definition]{Remark}

 \newtheorem{example}[definition]{Example}

\theoremstyle{plain}      

 \newtheorem{proposition}[definition]{Proposition}

 \newtheorem{theorem}[definition]{Theorem}

 \newtheorem{corollary}[definition]{Corollary}

 \newtheorem{lemma}[definition]{Lemma}

\newcommand{\C}{\mathbb C}
\newcommand{\R}{\mathbb R}
\newcommand{\Z}{\mathbb Z}
\newcommand{\HH}{\mathbb H}
\newcommand{\PP}{\mathbb P}
\newcommand{\model}{\mathbb M}
\newcommand{\barmodel}{\overline{\mathbb M}}
\newcommand{\Teich}{\mathcal{T}(S)}
\newcommand{\metTeich}{\mathcal{T}_{hyp.}(S)}
\newcommand{\Fricke}{{\mathfrak F}(S)}

\newcommand{\QC}{{\mathcal Q}{\mathcal C}}
\newcommand{\Mod}{\mathcal{M}(S)}
\newcommand{\mcg}{\mathsf{Mod}(S)}

\newcommand{\met}{\mathsf{Met}(S)}
\newcommand{\hypmet}{\mathsf{Met}_{hyp.}(S)}

\newcommand{\barTeich}{\overline{\mathcal{T}}(S)}
\newcommand{\barMod}{\overline{\mathcal{M}}(S)}
\newcommand{\thurston}{\overline{\mathcal{T}}^{th.}(S)}
\newcommand{\bthurston}{\partial{\mathcal{T}}^{th.}(S)}
\newcommand{\qf}{\mathcal{QF}(S)}
\newcommand{\dwp}{\mathop{{\rm d}_{wp}}\nolimits}
\newcommand{\dt}{\mathop{{\rm d}_{\mathcal T}}\nolimits}

\newcommand{\tl}{\mathop{{\mathfrak L}_\rho}\nolimits}
\newcommand{\twp}{\mathop{{\mathfrak L}_{wp}}\nolimits}
\newcommand{\Diff}{\mathsf{Diff}(S)}
\newcommand{\Diffo}{\mathsf{Diff}_0(S)}
\newcommand{\Diffp}{\mathsf{Diff}^+(S)}

\newcommand{\MF}{\mathsf{MF}(S)}
\newcommand{\PMF}{\mathsf{PMF}(S)}

\newcommand{\id}{\mathsf{id}}

\newcommand{\CM}{\mathcal{CM}}

\newcommand{\lra}{\longrightarrow}
\newcommand{\simrightarrow}{\buildrel \sim\over\longrightarrow}

\newcommand{\Tr}{\operatorname{\rm Tr}}
\newcommand{\rk}{\operatorname{\rm rk}}
\newcommand{\ad}{\operatorname{\rm ad}}

\newcommand{\Hom}{\operatorname{\rm Hom}}
\newcommand{\iso}{\operatorname{\mathsf Iso}}
\newcommand{\ext}{\operatorname{\mathsf Ext}} 
\newcommand{\Aut}{\operatorname{\mathsf Aut}}

\newcommand{\Out}{\operatorname{\mathsf Out}(S)}

\newcommand{\PSL}{\operatorname{\mathsf PSL}}
\newcommand{\SL}{\operatorname{\mathsf SL}}
\newcommand{\SU}{\operatorname{\mathsf SU}}

\newcommand{\QD}{\operatorname{\mathsf QD}}

\newcommand{\ord}{\operatorname{\rm ord}}
\newcommand{\slope}{\operatorname{\rm slope}}
\newcommand{\length}{\operatorname{\rm length}}

\newcommand{\End}{\operatorname{\rm End}}

\newcommand{\Div}{\operatorname{\rm div}}
\newcommand{\Sym}{\operatorname{\rm Sym}}
\newcommand{\Ric}{\operatorname{\mathsf Ric}}
\newcommand{\Riem}{\operatorname{\mathsf Riem}}
\newcommand{\Hopf}{\operatorname{\mathsf Hopf}}

\newcommand{\real}{\operatorname{\mathsf{Re}}}
\newcommand{\imag}{\operatorname{\mathsf{Im}}}

\begin{document}

\title{Harmonic Maps and Teichm\"uller Theory}

\author{Georgios D. Daskalopoulos\thanks{Work partially
supported by NSF Grant DMS-0204191}\, and Richard A. Wentworth\thanks{Work partially
supported by NSF Grants DMS-0204496 and DMS-0505512}}

\address{Department of Mathematics \\
		Brown University \\
		Providence,  RI  02912 \\
email:\,\tt{daskal@math.brown.edu}
\\ \\
Department of Mathematics \\
   Johns Hopkins University \\
   Baltimore, MD 21218\\
email:\,\tt{wentworth@jhu.edu}
}

\maketitle




\tableofcontents

\section{Introduction}  \label{intro}

Teichm\"uller theory is rich in applications to topology and physics.  
By way of the mapping class group the subject is closely related to knot theory and three-manifolds.  From the uniformization theorem, Teichm\"uller theory is part of the more general study of Kleinian groups and character varieties.  Conformal field theory and quantum cohomology make use of the algebraic and geometric properties of the  Riemann moduli space.

At the same time, analytic techniques have been important in Teichm\"uller theory almost from the very beginning of the subject.  Extremal maps and special metrics give alternative perspectives to moduli problems and clarify our understanding of a wide range of results.  In some cases they can be used to obtain new properties.

The goal of this paper is to present some of the more recent activity using analysis, and in particular harmonic maps, in the context of Teichm\"uller theory, representations of surface groups, mapping class groups, and Weil-Petersson geometry.  Topics have been selected in order to illustrate the theme that the analytic and topological points of view complement each other in a useful way.  For example, we will present four different proofs of the fact that Teichm\"uller space is a cell, and we will discuss the recent completion of a harmonic maps approach to Teichm\"uller's existence and uniqueness theorems on extremal quasiconformal maps.
Instead of a systematic survey of the subject, we have chosen  to present the ideas behind the results through examples and in a rather informal way.  There are very few proofs, but hopefully the references given at the end will provide the interested reader sufficient recourse for more details.

This paper makes no attempt to exhaust all aspects of this subject.   In particular, no mention is made of the  work on quasiconformal harmonic maps of the disk and Schoen's conjecture (see \cite{Li-Tam1, Li-Tam2, Rivet, Tam-Wan, Wan}), or of the universal Teichm\"uller space in general.
Other topics that have been covered in great detail in the literature have also been omitted or only briefly touched upon.  For example, there is little discussion of the complex analytic theory of Teichm\"uller space, the Bers embedding, Royden's theorem on automorphisms, etc.
  For the same reason, our summary of Weil-Petersson geometry is rather brief, and instead we refer to Wolpert's recent survey \cite{Wo5}.
  
  Finally, while we have tried to give complete and accurate references to the results stated in this paper, given the expanse of the subject there will inevitably be omissions.  For these we offer our apologies in advance.  
Two useful surveys of earlier results on harmonic maps are \cite{ES2} and \cite{Schoen1}.  Relatively recent general texts on Teichm\"uller theory are \cite{Ab, IT, Nag}. The point of view taken in Tromba's book \cite{Tromba} is especially relevant to the material presented here.   For an interesting account of Teichm\"uller's life and work, see Abikoff \cite{Ab2}.\index{Teichm\"uller}

\section*{Notation}
For simplicity, this paper will deal with connected compact oriented surfaces without boundary and of genus $p\geq 2$.   The notation we shall use is the following:  $S$ will denote the underlying smooth surface, and $j$ will denote a complex structure on $S$.  Hence, a Riemann surface is a pair $(S,j)$.  
The hyperbolic metric on $S$ will be denoted by $\sigma$.  Since it is uniquely determined by and uniquely determines the complex structure, the notation $(S,j)$ and $(S,\sigma)$ will both be understood to represent a Riemann surface structure.
When the complex structure is understood we shall  use letters $S$ and $R$ alone to denote  Riemann surfaces, and hopefully this will not cause confusion.  The following are some of the commonly used symbols in this paper:
\begin{itemize}
\item $\id=$ identity map; 
\item ${\bf I}=$ identity endomorphism;
\item $f\sim f'$ homotopic maps;
\item $\deg(f)=$ the degree of a map between surfaces;
\item $K(f)=$ the dilatation of a quasiconformal map (Section \ref{S:qc});
\item $\Gamma=\pi_1(S)$, or $=\pi_1(M)$ for a manifold $M$;
\item $\widetilde M=$ the universal cover of $M$;
\item $\Omega^p=$ the space of smooth $p$-forms;
\item $\Lambda=$ a Fuchsian group (Section \ref{S:uniformization});
\item $D=$ the unit disk in $\C$;
\item $\HH=$ the upper half plane in $\C$;
\item $\HH^3=$ hyperbolic 3-space $\simeq \SL(2,\C)/\SU(2)$;
\item $i(a,b)=$ the geometric intersection number of  simple closed curves  $a,b$ on $S$;
\item $\ell_c[\sigma]=$ the length of a simple closed curve on $S$ with respect to the hyperbolic metric $\sigma$;
\item ${\mathcal F}=$ a measured foliation on $S$ (Section \ref{S:trees});
\item $i([c],{\mathcal F})=$ the intersection number of an isotopy class of simple closed curves with a measured foliation $\mathcal F$ (see Section \ref{S:trees});
\item $i({\mathcal F}_1, {\mathcal F}_2)=$ the intersection number of a pair of measured foliations (see Section \ref{S:trees});
\item $T_{\mathcal F}=$ the $\R$-tree dual to a measured foliation $\mathcal F$ (Section \ref{S:trees});
\item $\MF$ (resp.\ $\PMF$) $=$  the spaces of measured (resp.\ projective measured) foliations  on $S$ (Section \ref{S:trees});
\item $K_S=$ the canonical line bundle on a Riemann surface $S$;
\item $\chi_S=$ the Euler characteristic of $S$;
\item $\nabla=$ the  covariant derivative, or a  connection on a vector bundle $V$;
\item $d_\nabla=$ the de Rham operator, twisted by a connection $\nabla$;
\item $\nabla_H=$ the Chern connection on a holomorphic bundle with hermitan metric $H$ (Section \ref{S:hs};
\item $F_{\nabla}=$ the curvature of a connection $\nabla$;
\item $\Delta=$ the Laplace-Beltrami operator;
\item $\mu=$ a Beltrami differential (Section \ref{S:qc});
\item $\Vert\mu\Vert_\infty=$ the $L^\infty$ norm of a Beltrami differential $\mu$;
\item $\varphi=$ a holomorphic quadratic differential (Section \ref{S:qd});
\item $\Vert\varphi\Vert_1$ (resp.\ $\Vert \varphi\Vert_2$) $=$ the $L^1$ (resp.\ $L^2$) norms of a quadratic differential $\varphi$ (see eqs.\ \eqref{E:L1} and \eqref{E:L2});
\item $T_\varphi=$ the $\R$-tree dual to the horizontal foliation of $\varphi$ (Section \ref{S:trees});
\item $\QD(S)=$ the space of holomorphic quadratic differentials;
\item $\Fricke=$ the Fricke space (Section \ref{S:uniformization});
\item $\chi(\Gamma,r)$ (resp.\ $\chi(\Gamma)$) $=$ the $\SL(r,\C)$ (resp.\ $\SL(2,\C)$) character varieties of $\Gamma$ (Section \ref{S:hs});
\item $\Teich=$ Teichm\"uller space (Section \ref{S:teich});
\item $\barTeich=$ the Weil-Petersson completion of $\Teich$ (Section \ref{S:wp});
\item $\Diff, \Diffp, \Diffo=$ the diffeomorphisms, orientation preserving diffeomorphisms, and diffeomorphisms connected to the identity of a surface $S$;
\item $\mcg=$ the mapping class group (Section \ref{S:mcg});
\item $\Mod=$ the Riemann moduli space (Section \ref{S:mcg});
\item $\barMod=$ the Deligne-Mumford compactification of $\Mod$ (Section \ref{S:mcg});
\item $\dt=$ the Teichm\"uller metric on $\Teich$  (see eq.\ \eqref{E:teichmetric});
\item $\dwp=$ the Weil-Petersson metric on $\Teich$ (Section \ref{S:wp});
\item $\iso(X)=$ the isometry group of a metric space $(X,d)$;
\item $\partial X=$ the ideal boundary of an NPC space $X$ (Section \ref{S:equivariant});
\item $\tl=$ the translation length function of a representation \eqref{E:tl};
\item $\twp[\phi]=$ the Weil-Petersson translation length of $[\phi]\in \mcg$ (see eq.\ \eqref{E:wptl});
\item $H^1$ (resp.\ $H^1_{loc.}$) $=$ the Sobolev space of square integrable (resp.\ locally square integrable) functions with square integrable (resp.\ locally square integrable) distributional derivatives;
\item $e(f)=$ the energy density of a map $f$ (see eq.\  \eqref{E:ef});
\item $\pi_{\alpha\beta}=$ the directional energy tensor (see eq.\ \eqref{E:pi});
\item $E(f)=$ the energy of a map $f$ (see eq.\ \eqref{E:Ef});
\item $\End(V)$ (resp.\ $\End_0(V)$) $=$ the endomophism (resp.\ traceless endomorphism) bundle of a complex vector bundle $V$ (Section \ref{S:hs});
\item $\ad(V)$ (resp.\ $\ad_0(V)$) $=$ the skew-hermitian (resp.\ traceless skew-hermitian) endomorphism bundle of a hermitian vector bundle $V$ (Section \ref{S:hs});
\item $\Phi=$ a Higgs field (Section \ref{S:hs});
\item ${\mathfrak M}(S,r)=$ the moduli space of polystable Higgs bundles of rank $r$ on $S$ (Section \ref{S:hs}).
\end{itemize}


\section{Teichm\"uller Space and Extremal Maps}

\begin{itemize}
\item 2.1 The Teichm\"uller Theorems
\item 2.2 Harmonic Maps
\item 2.3 Singular Space Targets
\end{itemize}

\vspace{-.25in}


\subsection{The Teichm\"uller Theorems}


This section gives a summary of the basics of Teichm\"uller theory from the point of view of quasiconformal maps.  Section \ref{S:uniformization} reviews the uniformization theorem and the Fricke space.  In Section \ref{S:qc}, we introduce quasiconformal maps, Beltrami differentials, and we state the basic existence theorem for solutions to the Beltrami equation. We also formulate the extremal problem.  In Section \ref{S:qd}, we review quadratic differentials, Teichm\"uller maps, and Teichm\"uller's existence and uniqueness theorems.   In Section \ref{S:teich}, we define the  Teichm\"uller space based on a Riemann surface and discuss the first approach to Teichm\"uller's theorem on the contractibility of Teichm\"uller space.  The proof that we give here is based on the notion of extremal maps, i.e.\ quasiconformal maps that minimize dilatation in their homotopy class.  The connection between extremal and harmonic maps will be explained in Section \ref{S:gr}. Finally, in Section \ref{S:metric}, we provide an alternative definition of Teichm\"uller space via hyperbolic metrics.


\subsubsection{Uniformization and the Fricke space.} \label{S:uniformization}


 The famous \emph{uniformization theorem}\index{uniformization theorem} of Poincar\'e, Klein, and Koebe states that every closed Riemann surface $S$ of genus at least $2$ is biholomorphic to a quotient $\HH/\Lambda$, where $\HH$ denotes the upper half plane and $\Lambda$ is a group of holomorphic automorphisms of $\HH$ acting freely and properly discontinuously.  Such a group can be identified with a  discrete subgroup of $\PSL(2,\R)$, i.e.\ a \emph{Fuchsian group} (cf.\ \cite{For,IT}).\index{Fuchsian group}
 
 On $\HH$ we have the \emph{Poincar\'e metric}\index{Poincar\'e metric}
$$
ds_{\HH}^2=\frac{|dz|^2}{(\imag z)^2}\ .
$$
Under the biholomorphism ${\mathfrak h}:\HH\to D$ given by ${\mathfrak h}(z)=(z-i)/(z+i)$, $ds_{\HH}^2={\mathfrak h}^\ast ds_D^2$, where
$$
ds_{D}^2=\frac{4|dz|^2}{(1-|z|^2)^2}\ .
$$
By a straightforward calculation the curvature of the Poincar\'e metric is constant equal to $-1$, and by Pick's Theorem its isometry group is $\PSL(2,\R)$ (cf.\ \cite{IT}).  In particular, this metric descends to the quotient $\HH/\Lambda$.
Hence, every Riemann surface of genus $\geq 2$ has a hyperbolic metric, and this metric is unique.  On the other hand, \emph{any} Riemannian metric induces a unique complex structure.  This is a consequence of Gauss' theorem on the existence of \emph{isothermal coordinates}:\index{isothermal coordinates} if
$(S,g)$ is an oriented surface with Riemannian metric $g$, then $S$ admits a unique complex structure $j$ such that in local complex coordinates  $g=g(z)|dz|^2$, where $g(z)$ is a smooth, positive (local) function.  Hence, specifying a complex structure on the topological surface $S$ is equivalent to specifying a hyperbolic metric.
We will use Greek letters, e.g. $\sigma=\sigma(z)|dz|^2$, to distinguish the hyperbolic from arbitrary Riemannian metrics $g$.

 Let $\Fricke$\index{Fricke space} denote the Fricke space of conjugacy classes of discrete embeddings $\Gamma=\pi_1(S)\to \PSL(2,\R)$.  Then $\Fricke$ inherits a topology as a \emph{character variety} (cf.\ \cite{CS, Goldman, GMill} and Section \ref{S:hs} below). \index{character variety} The idea is to choose a \emph{marking} of the genus $p$ surface $S$,\index{marking} namely, a presentation 
 $$
 \Gamma=\langle a_1, \ldots, a_p, b_1,\ldots, b_p : \prod_{i=1}^p [a_i,b_i]=\id\rangle\ .
 $$
 where the $a_i$ and $b_i$ are represented by simple closed curves on $S$ with \emph{geometric intersection numbers} satisfying $i(a_i,b_j)=\delta_{ij}$, $i(a_i,a_j)=i(b_i,b_j)=0$.
 A homomorphism $\rho:\Gamma\to \PSL(2,\R)$ is determined by specifying $2p$ elements $A_i, B_i\in \PSL(2,\R)$ satisfying the relation $\prod_{i=1}^p[A_i,B_i]={\bf I}$.  A naive dimension count (which can easily be made precise at irreducible representations $\rho$) suggests that the dimension of the space of such homomorphisms is $6p-3$.  Since  $\PSL(2,\R)$ acts by conjugation, producing a  $3$-dimensional orbit, we have $\dim\Fricke=6p-6$.
  Indeed, since the Fricke space consists of \emph{discrete} embeddings, a more precise analysis can be given which realizes $\Fricke$ as a subset of $\R^{6p-6}$ (cf.\ \cite{Ab}).  
  \begin{proposition} \label{P:fricke}  The Fricke space $\Fricke$ is embedded in $\R^{6p-6}$.
    \end{proposition}
    It is this embedding (the details of which will not be important) that we will use to define the topology on $\Fricke$.  We shall see below that $\Fricke$ is homeomorphic to Teichm\"uller space (Theorem \ref{T:teich2}).


\subsubsection{Quasiconformal maps.} \label{S:qc}


An orientation preserving homeomorphism $f$ of  a domain $\Omega\subset\C$ into $\C$ is called \emph{$K$-quasiconformal}\index{quasiconformal map} (or \emph{$K$-qc}) if
\begin{enumerate}
\item $f$ is of Sobolev class $H^1_{loc.}$, i.e.\ the distributional derivatives $f_z$, $f_{\bar z}$ are locally square integrable on $\Omega$;
\item there exists a constant $0\leq k<1$ such that
$|f_{\bar z}|\leq k|f_z|$, almost everywhere on $\Omega$, where $K=(1+k)/(1-k)$.  
\end{enumerate}
The infimum of $K\geq 1$ such that $f$ is $K$-qc is called the \emph{dilatation} of $f$,\index{dilatation}  \index{quasiconformal map!dilatation of} and it is denoted by $K(f)$.  Clearly, $1$-qc is equivalent to conformal.

An orientation preserving homeomorphism $f:S\to R$ between two Riemann surfaces  is called $K$-qc if its lift to the universal cover $\tilde f:\HH\to\HH$ is $K$-qc.  We define the dilatation $K(f)$ of $f$ to be $K(\tilde f)$.  Given such a map $f$, let 
\begin{equation} \label{E:QC}
\QC[f]=\{ f':S\to R : \text{$f'$ is  a qc homeomorphism homotopic to $f$}\}\ .
\end{equation}

The  main extremal problem in Teichm\"uller theory is a generalization to closed surfaces of Gr\"otzsche's problem for rectangles (see \cite{Ab}):  given a qc map $f:S\to R$, let 
\begin{equation} \label{E:extdilatation}
K^\ast[f]=\inf_{f'\in\QC[f]} K(f') \ .
\end{equation}
\par\indent
{\bf Teichm\"uller's Extremal Problem.}\index{Teichm\"uller!extremal problem}
Is $K^\ast[f]$  realized as the dilatation of a qc map, and if so, what are the properties of the map?

\medskip\noindent
A qc homeomorphism $f$ such that $K(f)=K^\ast[f]$ is called an \emph{extremal map}.\index{quasiconformal map!extremal}\index{extremal map}  The existence of extremal maps is a relatively easy consequence of compactness properties of quasiconformal maps.  The emphasis of this problem is therefore on the uniqueness and  characterization of extremal maps.
  We will give Teichm\"uller's answer to this question in the next section.

Choose coordinates $(U,z)$ on $S$ and $(V,w)$ on $R$ and set $F=w\circ f\circ z^{-1}$.  Define the \emph{Beltrami coefficient of} $f$\index{quasiconformal map!Beltrami coefficient of}\index{Beltrami!coefficient} with respect to the choice of coordinates by
$$
\mu_f=\mu_f(z) d\bar z\otimes (dz)^{-1}=F_{\bar z}/F_z d\bar z\otimes (dz)^{-1}\ .
$$
By (2), $|\mu_f(z)|<1$ almost everywhere.
 The above expression is independent of the choice of coordinates $w$ and  transforms tensorially with respect to coordinate changes in $z$.  More precisely, $\mu_f$ may be regarded as an $L^\infty$-section of the  bundle $\overline K_S\otimes K_S^{-1}$, where $K_S$ is the canonical line bundle of $S$.  Notice, however, that $|\mu_f(z)|$ is independent of a choice of conformal coordinates.  Set $\Vert \mu\Vert_\infty$ to be the essential supremum of $|\mu_f|$ over $S$.

Let ${\mathcal B}(S)$ denote the Banach space of $L^\infty$-sections of $\overline K_S\otimes K_S^{-1}$ with the $L^\infty$-norm.  Set
$$
{\mathcal B}_1(S)=\{ \mu\in{\mathcal B}(S) : \Vert \mu\Vert_\infty<1\}\ .
$$
For any qc map $f:S\to R$ we associate $\mu_f\in {\mathcal B}_1(S)$.  If $S=\HH/\Lambda$ a Beltrami differential\index{Beltrami!differential} on $S$ can be identified with an $L^\infty$ function $\tilde \mu_f$ on $\HH$ satisfying the equations of automorphy
\begin{equation} \label{E:beltrami1}
\tilde \mu_f(\gamma z)\frac{\overline{\gamma'(z)}}{\gamma'(z)}=\tilde \mu_f(z)\ ,\quad z\in \HH\ ,\ \gamma\in\Lambda\ .
\end{equation}
Furthermore, qc homeomorphisms $\tilde f$ of $\HH$ whose Beltrami coefficients satisfy \eqref{E:beltrami1} give deformations of Fuchsian groups\index{Fuchsian group!deformation of} via
\begin{equation} \label{E:fuchsian}
\Lambda\leadsto\Lambda_\mu \ :\ \gamma\in \Lambda\mapsto \tilde f\circ\gamma\circ {\tilde f}^{-1}\in \PSL(2,\R)\ .
\end{equation}
Specifying the Beltrami coefficient and solving for a qc map is called \emph{Beltrami's equation}.\index{Beltrami!equation}
The following is the fundamental existence theorem for solutions to Beltrami's equation.  The seminal reference is Ahlfors \cite{Ah3}.  See also \cite[Chapter 4]{IT}.
\begin{theorem} \label{T:beltrami}
For any Beltrami differential $\mu\in  {\mathcal B}_1(\C)$ there exists a unique qc homeomorphism $f^\mu$ of $\HH$, extending continuously to $\widehat\HH=\HH\cup\{\infty\}$,  whose Beltrami  coefficient is $\mu_{f^\mu}=\mu$, and which fixes the points $0$, $1$, and $\infty$.  Furthermore, $f^\mu$ depends complex analytically on $\mu$.
\end{theorem}

\begin{corollary} \label{C:beltrami}
For any Beltrami differential $\mu\in  {\mathcal B}_1(S)$ there exists a unique qc homeomorphism $f^\mu : S\to R$, for some Riemann surface $R$.  More precisely, if $S=\HH/\Lambda$, then $R=\HH/
{\tilde f}^{\mu}\circ \Lambda\circ (\tilde f^\mu)^{-1}$, where $\tilde f^\mu$ is the solution in Theorem \ref{T:beltrami} for the pullback Beltrami differential.  Furthermore, $f^\mu$ depends complex analytically on $\mu$.
\end{corollary}

Hence, Beltrami differentials can be used to parametrize the Fricke space $\Fricke$.
Of course, there is an infinite dimensional family of Beltrami differentials giving conjugate Fuchsian groups.


\subsubsection{Quadratic differentials and Teichm\"uller maps.} \label{S:qd}


By a \emph{holomorphic quadratic differential}\index{quadratic differential} on a Riemann surface $S$ we mean a holomorphic section of the line bundle $K_S^{ 2}$.  Set $\QD(S)=H^0(S,K_S^{ 2})$.  By the Riemann-Roch Theorem, $\QD(S)$ is a complex vector space of dimension $3p-3$, where $p$ is the genus of $S$.  If $\varphi\in\QD(S)$, then in local conformal coordinates (centered at $z_0$, say)  $\varphi=\varphi(z)dz^2$, where $\varphi(z)$ is a local holomorphic function.  By a  coordinate change we can write $\varphi=z^k dz^2$, where $k=0,1,2,\ldots$ is the order of vanishing of $\varphi$ at $z_0$.  The coordinate system
$$
w(z)=\int_{z_0}^z \sqrt\varphi=\int_{z_0}^z\sqrt{\varphi(z)} dz=\frac{2}{k+2} z^{\frac{k+2}{2}}$$
will be called the $\varphi$-coordinates around $z_0$ (if $m$ is odd, this is multi-valued).  Writing $w=u+iv$, the foliations $v=\text{constant}$ and $u=\text{constant}$ are called the \emph{horizontal} and \emph{vertical foliations} of $\varphi$, respectively. \index{quadratic differential!horiz., vert. foliations}

\setlength{\unitlength}{1cm}
\begin{picture}(14,5.5)
\ifx\pdftexversion\undefined
\put(0,2.2){
{\scalebox{.5}{\includegraphics{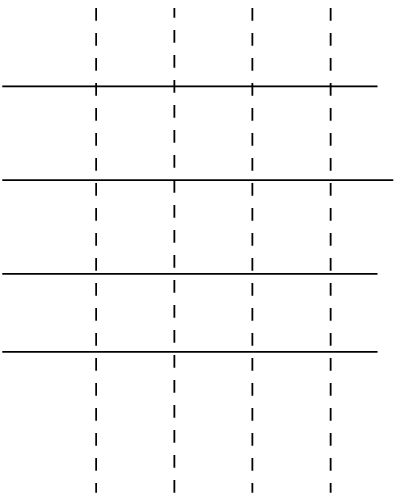}}}}
\put(3.25,2){
{\scalebox{.35}{\includegraphics{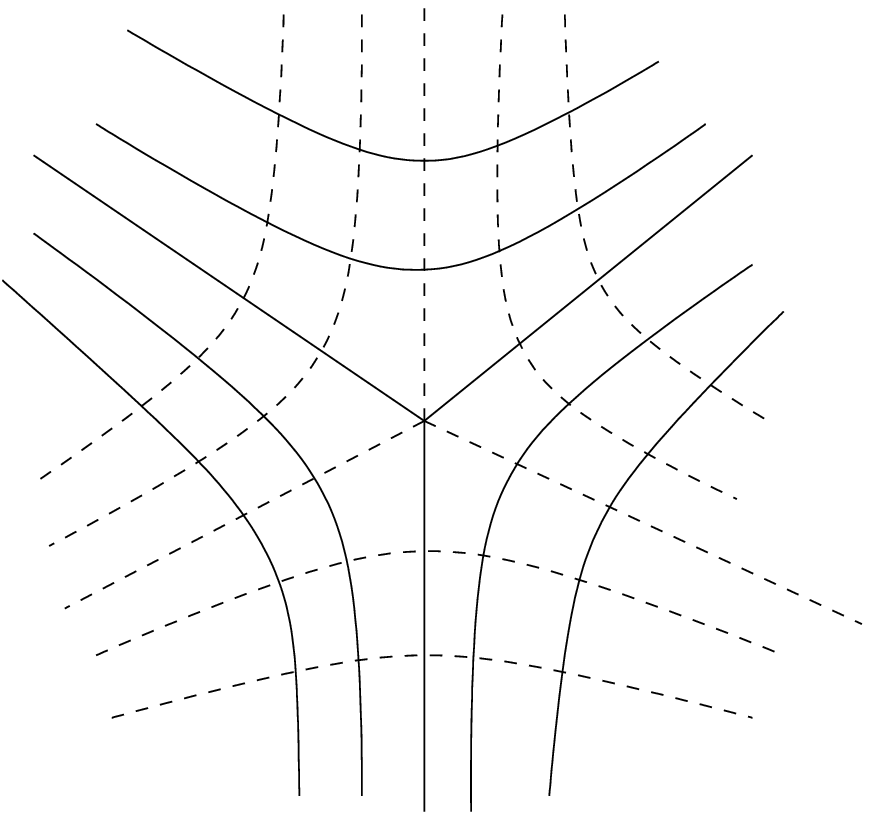}}}}
\put(7.25,2){
{\scalebox{.4}{\includegraphics{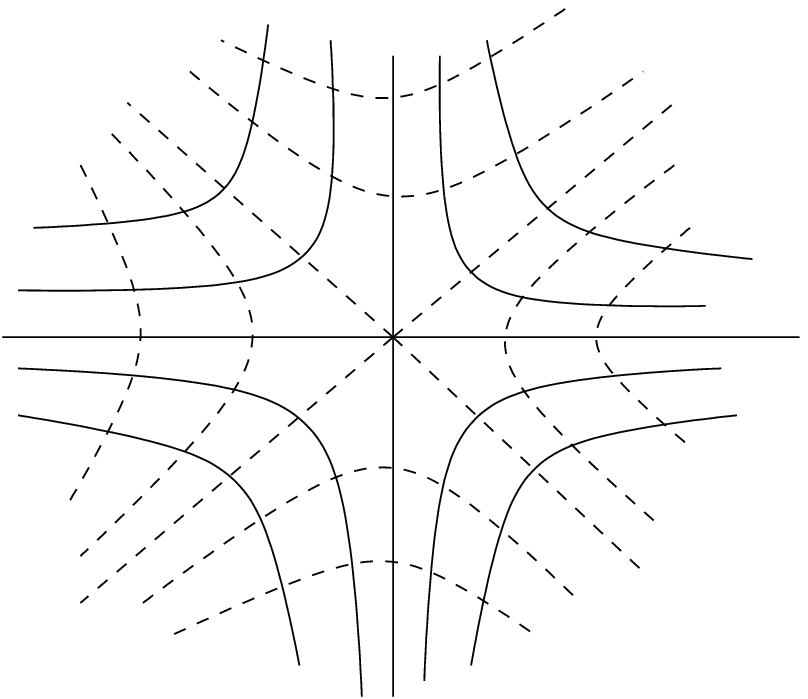}}}}
\else
\put(0,2,2){
{\scalebox{.5}{\includegraphics{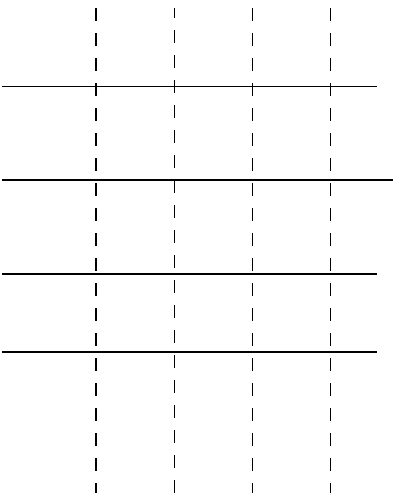}}}}
\put(3.25,2){
{\scalebox{.35}{\includegraphics{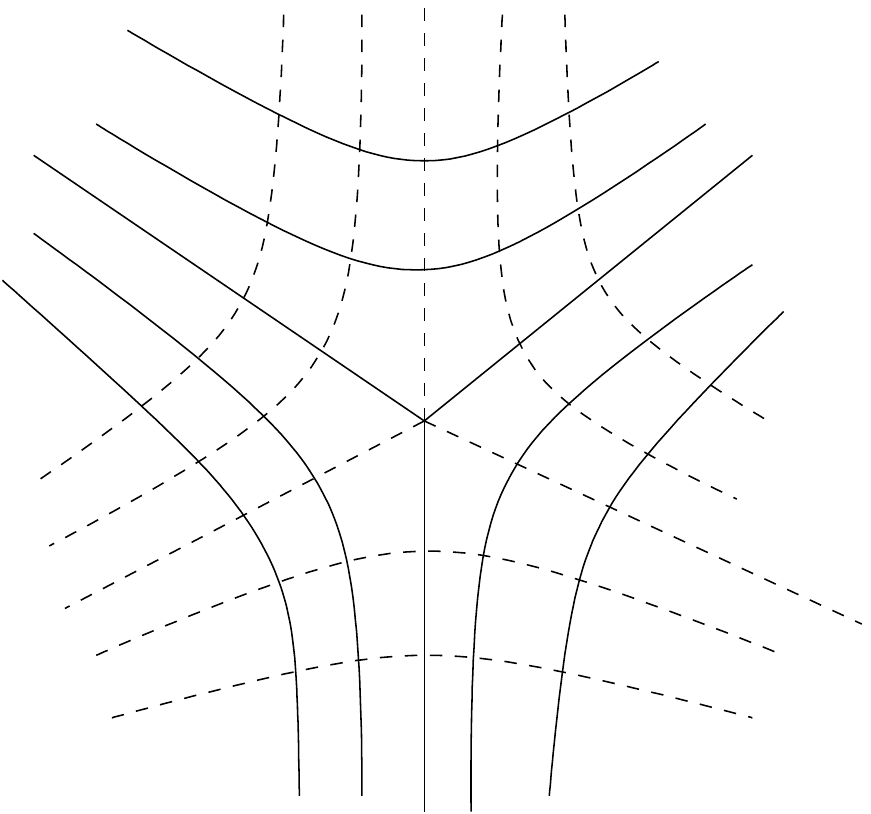}}}}
\put(7.25,2){
{\scalebox{.4}{\includegraphics{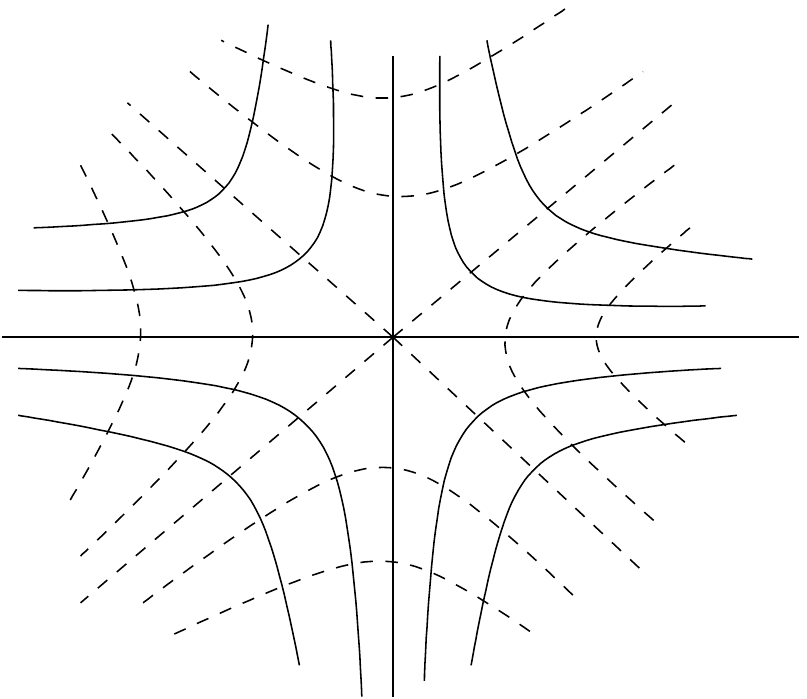}}}}
\fi
\put(.4,1){$\varphi(z_0)\neq0$}
\put(3.5,1){$\varphi(z_0)=0$, $k=1$}
\put(7.5,1){$\varphi(z_0)=0$, $k=2$}
\put(4.5,0){Figure 1.}
\end{picture}

\smallskip
\noindent    For more details, we refer to Strebel's treatise on quadratic differentials \cite{Str}.

A holomorphic quadratic differential on $S=\HH/\Lambda$ is given by $\widetilde\varphi=\tilde \varphi(z) dz^2$, where $\tilde \varphi$ is a holomorphic function on $\HH$ satisfying the equations of automorphy
\begin{equation} \label{E:auto}
\tilde \varphi(\gamma z)\gamma'(z)^2=\tilde \varphi(z)\ ,\quad z\in \HH\ ,\ \gamma\in\Lambda\ .
\end{equation}
Set
$$
\QD_1(S)=\{\varphi\in\QD(S) : \Vert \varphi\Vert_1<1\}
$$
where $\Vert\cdot\Vert_1$ denotes the $L^1$-norm: 
\begin{equation} \label{E:L1}
\Vert\varphi\Vert_1=\int_S|\varphi(z)|dxdy
\end{equation}
  Given $\varphi\in \QD_1(S)\setminus\{0\}$, we say that a qc homeomorphism $f:S\to R$ is a \emph{Teichm\"uller map} \index{Teichm\"uller!map} for $\varphi$ if the Beltrami coefficient of $f$ satisfies
\begin{equation} \label{E:ext}
\mu_f=k\frac{\bar\varphi}{|\varphi|}\ ,\qquad k=\Vert\varphi\Vert_1\ .
\end{equation}

We are now in a position to give Teichm\"uller's solution to the extremal problem stated in the previous section.  First, a \emph{Teichm\"uller map is uniquely extremal}.

\begin{theorem}[Teichm\"uller's Uniqueness Theorem] \label{T:uniqueness}\index{Teichm\"uller!uniqueness theorem}
Let $f:S\to R$ be a Teichm\"uller map.  Then every $f'\in\QC[f]$ satisfies
$$
\Vert \mu_{f'}\Vert_\infty\geq \Vert \mu_f\Vert_\infty\qquad(\text{equivalently,}\ K(f')\geq K(f))\ .
$$
Moreover, the equality holds if and only if $f'=f$.
\end{theorem}

The second result asserts that  \emph{Teichm\"uller maps always exist}.
\begin{theorem}[Teichm\"uller's Existence Theorem] \label{T:existence}\index{Teichm\"uller!existence theorem}
In the homotopy class of every qc homeomorphism $f: S\to R$ there is either a conformal map or  a Teichm\"uller map.
\end{theorem}

 Theorems \ref{T:uniqueness} and \ref{T:existence} were first stated by Teichm\"uller  (see \cite{Te}).  His papers are ``generally considered unreadable" (Abikoff, \cite[p.\ 36]{Ab}).  Subsequent proofs were given in \cite{Ah1} and \cite{Bers1} (see also \cite{GL, Ham1}).     Below we outline a proof of these two fundamental results  based on harmonic maps to singular spaces (see Section \ref{S:gr}).

 Teichm\"uller maps are essentially affine with respect to a natural choice of coordinates (see \cite{IT}):

\begin{theorem} \label{T:affine}
Fix $\varphi\in \QD_1(S)\setminus\{0\}$, $k=\Vert\varphi\Vert_1<1$, and let $f:S\to R$ be a Teichm\"uller map for $\varphi$.  Then there exists a unique holomorphic quadratic differential $\psi$ on $R$ satisfying the following conditions
\begin{enumerate}
\item  If $z$ is a zero of $\varphi$  then $f(z)$ is a zero of $\psi$ of the same order;
\item If $z$ is not a zero of $\varphi$ and $\zeta$ is a $\varphi$-coordinate about $z$, then there exists a $\psi$-coordinate $w$ at $f(z)$ such that
\begin{equation} \label{E:dilatation}
w\circ f=\frac{\zeta +k\bar\zeta}{1-k}\ .
\end{equation}
\end{enumerate}
\end{theorem}
The quadratic differentials $\varphi$ and $\psi$ are called the \emph{initial} and \emph{terminal} differentials of the Teichm\"uller map $f$, respectively.\index{quadratic differential!initial, terminal}


\subsubsection{The Teichm\"uller space.} \label{S:teich}


We now come to the definition of Teichm\"uller space.\index{Teichm\"uller!space}  Let $S$ be a closed Riemann surface of genus $p\geq 2$.  Consider triples $(S,f, R)$, where $R$ is a Riemann surface  and $f:S\to R$  is an orientation preserving diffeomorphism.  Triples $(S,f_1, R_1)$ and $(S,f_2,R_2)$ are said to be equivalent if $f_2\circ f_1^{-1} : R_1\to R_2$ is homotopic to a biholomorphism.   The set of all equivalence classes $[S,f, R]$ of triples $(S,f,R)$  is denoted $\Teich$ and is called the \emph{Teichm\"uller space based on $S$}.  The definition of $\Teich$ turns out to be independent of the complex structure on $S$ (see Theorem \ref{T:teichmetric} below).
Since any homeomorphism (in particular quasiconformal ones) is homotopic to a diffeomorphism, one obtains the same space if one considers pairs $(S,f, R)$ where $f$ is quasiconformal.  This is a point of subtlety when dealing with Riemann surfaces with punctures.

Restricting as we are to the case of closed surfaces,
Teichm\"uller space may be regarded as parametrizing complex structures up to biholomorphisms connected to the identity.  Indeed, if $S_0=(S,j_0)$ denotes the basepoint
 and $(S,j)$ is another complex structure on the underlying surface $S$, then by choosing $f=\id$ and $R=(S,j)$  there is an associated  point $[j]=[S_0,\id, R]\in \Teich$.  Two points $[S_0,\id,R_1]$ and $[S_0,\id,R_2]$ obtained in this way are equivalent if and only if $(S,j_1)$ and $(S,j_2)$ are biholomorphic via a map connected to the identity.  Conversely, given any triple $(S_0,f,R)$, let $j$ denote the pullback by $f$ of the complex structure  on $R$ to the underlying surface $S$.  Then by definition $f:(S,j)\to R$ is a biholomorphism; hence, $(S_0,\id,(S,j))$ is equivalent to $(S_0,f,R)$.  With this understood, we sometimes represent points in $\Teich$ by equivalence classes $[j]$.

Given $[j_1], [j_2]\in \Teich$, recall that  $\QC[\id]$ is the set of all qc homeomorphisms $(S,j_1)\to (S,j_2)$ homotopic to the identity.  The \emph{Teichm\"uller metric} is defined
\begin{equation} \label{E:teichmetric} \index{Teichm\"uller!metric}
\dt([j_1],[j_2])=\inf_{f\in{\QC[\id]}} \log K(f)\ .
\end{equation}
For the next result we refer to \cite[\S 5.1]{IT}.
\begin{theorem} \label{T:teichmetric}
$\Teich$ is a complete metric space with respect to the Teichm\"uller metric $\dt$.  Furthermore, given $[S,f, R]\in \Teich$, the map $[f]_\ast : {\mathcal T}(S)\to {\mathcal T}(R)$ given by 
$
[S,f', R']\mapsto [R, f'\circ f^{-1}, R']
$
is an isometry.
\end{theorem}
Henceforth, the topology on $\Teich$ is that given by the metric $\dt$.  Also, in light of the theorem we identify all Teichm\"uller spaces independent of the choice of base point. 
 Now we are ready for Teichm\"uller's third result.

\begin{theorem}[Teichm\"uller's Theorem] \label{T:teich}\index{Teichm\"uller!theorem}
$\Teich$ is homeomorphic to a cell of dimension $6p-6$.
\end{theorem}

By Corollary \ref{C:beltrami} on solutions to Beltrami's equation, Teichm\"uller maps with initial differential  $\varphi$ exist for any $\varphi\in\QD_1(S)$.  Hence, we may define a map
\begin{equation} \label{E:tau}
\tau : \QD_1(S)\lra\Teich\ : \ \tau(\varphi)=[S,f, R]\ ,
\end{equation}
where $f$ is a Teichm\"uller map for $\varphi\neq 0$, and $f=\id$, $R=S$, for $\varphi=0$.  
Theorem \ref{T:teich}  follows from
\begin{theorem} \label{T:teich2}
The map $\tau$ in \eqref{E:tau} is a homeomorphism.  Moreover, $\Teich$ is homeomorphic to $\Fricke$.
\end{theorem}
\begin{proof}
First, note that there is a natural bijection $F:\Teich\to \Fricke$ defined as follows: 
given $[S,f,R]\in \Teich$, by the uniformization theorem applied to the Riemann surface $R$ there is a discrete embedding  $\rho_R:\pi_1(R)\to \PSL(2,\R)$, determined up to conjugation. Since the diffeomorphism $f$ induces an isomorphism $f_\ast:\Gamma=\pi_1(S)\simrightarrow \pi_1(R)$, we obtain a  discrete embedding $\rho=\rho_R\circ f_\ast:\Gamma\to \PSL(2,\R)$.  Notice that if $[S,f_1,R_1]=[S,f_2,R_2]$, then the corresponding homomorphisms are conjugate.  Hence, there is a well-defined  point  $F[S,f,R]\in \Fricke$.  Conversely, given a discrete embedding $\rho:\Gamma\to \PSL(2,\R)$, consider the Riemann surface $R=\HH/\rho(\Gamma)$.   The Poincar\'e polygon theorem  realizes the boundary of a fundamental domain for the action of $\rho(\Gamma)$ as the lift of simple closed curves  $\alpha_i$, $\beta_i$  on $R$ satisfying the relations 
$i(\alpha_i,\beta_j)=\delta_{ij}$, $i(\alpha_i,\alpha_j)=i(\beta_i,\beta_j)=0$ (cf.\ \cite{Beardon}).
 The identification of $a_i$, $b_i$ with $\alpha_i$, $\beta_i$ fixes a homotopy class of diffeomorphisms $f:S\to R$, and it is clear that $F[S,f,R]=[\rho]$.  Hence, $F$ is a bijection.
Moreover, $F$ is continuous  by Corollary \ref{C:beltrami}, since a qc map of small dilatation is close to the identity, hence the corresponding deformation of the Fuchsian groups is small.
Consider the following diagram
\begin{displaymath}
 \xymatrix{ \QD_1(S) \ar[dr]_{G} \ar[r]^{\,\tau} & \ar[d]^F \Teich   \\
  & \Fricke}
 \end{displaymath}
 where $G=F\circ\tau$.  By Theorems \ref{T:uniqueness} and \ref{T:existence}, $\tau$ is a bijection.  It is also continuous.  Indeed, $\dt(\tau(0), \tau(\varphi))=\log((1+k)/(1-k))$, where $k=\Vert\varphi\Vert_1$ (recall that $\tau(0)=S$), so $\tau$ is clearly  continuous at the origin.  Continuity at general points follows from the change of basepoints in Theorem \ref{T:teichmetric}.  It follows that $G$ is a continuous bijection.  By the embedding $\Fricke\hookrightarrow\R^{6p-6}$ (Proposition \ref{P:fricke}) and  Invariance of Domain, $G$ is a homeomorphism;  hence, so are $F$ and $\tau$.
\end{proof}

We have proven Teichm\"uller's Theorem via his existence and uniqueness results (Theorems \ref{T:uniqueness} and \ref{T:existence}).  The proof uses the Fricke space $\Fricke$ and the finite dimensionality of the space of holomorphic quadratic differentials.  
 In  Sections \ref{S:wolf}, \ref{S:higgsproof}, and  \ref{S:proper} we shall give three alternative
  proofs of Theorem \ref{T:teich} using harmonic maps and the metric description of Teichm\"uller space.  


\subsubsection{Metric definition of Teichm\"uller space.}   \label{S:metric}


Let $S$ be an oriented surface of genus $p\geq 2$.  Let  $\hypmet$ be the space  of metrics with constant curvature $-1$.   This has a smooth structure inherited as a smooth submanifold of the space $\met$ of all smooth metrics on $S$.
As discussed in Section \ref{S:uniformization}, a hyperbolic metric defines a complex structure on $S$ via Gauss' theorem, and conversely, in every conformal class of metrics compatible with a given complex structure there is a unique hyperbolic metric.  
The group $\Diffo$ of diffeomorphisms isotopic to the identity acts on $\hypmet$ by pullback.  Define

 \begin{equation} \label{E:metteich}\index{Teichm\"uller!space}
 \metTeich=\hypmet/\Diffo\ ,
 \end{equation}
 with the quotient topology.   
 By constructing a slice for the action of $\Diffo$ on
 $\hypmet$ it is not hard to prove (see \cite{Ebin, FT1, Palais, Tromba}) 
 \begin{proposition} \label{P:smooth}
 $\metTeich $
is a smooth manifold of dimension $6p-6$.
\end{proposition}

To elaborate on this statement, we review the description of the tangent and cotangent spaces to $\Teich$ (for this approach, cf.\ \cite{DP, FT1}).  Let $\nabla$ denote the covariant derivative for a metric $g$ on $S$.  On the tangent space $T_g\met$ there is a natural $L^2$-pairing:
\begin{equation} \label{E:met}
\langle\delta g, \delta g'\rangle=\int_S(g^{\alpha\beta}g^{\mu\nu}\delta g_{\alpha\mu}\delta g'_{\beta\nu})
dvol_{(S,g)}
\end{equation}
where the metrics and variations are expressed with respect to local coordinates $\{x^\alpha\}$, $z=x^1+ix^2=x+iy$,  and repeated indices are summed.  For $\sigma$ a hyperbolic metric, the condition that $\delta\sigma$ be tangent to $\hypmet$ is
\begin{equation} \label{E:hyptangent}
0=(-\Delta+1)\Tr(\delta\sigma)+\nabla^\alpha\nabla^\beta(\delta\sigma_{\alpha\beta})\ ,
\end{equation}
where $\Delta$ is the Laplace-Beltrami operator associated to $\sigma$.  Finally, the tangent space to the orbit $\Diffo\cdot\sigma$  at $\sigma$ consists of Lie derivatives of the metric:
\begin{equation} \label{E:ds}
\delta\sigma_{\alpha\beta}=(L_v\sigma)_{\alpha\beta}=\nabla_\alpha v_\beta+\nabla_\beta v_\alpha
\end{equation}
for smooth vector fields $\{v^\alpha\}$.  From \eqref{E:met} and \eqref{E:ds}, the $L^2$-orthogonal in $T_\sigma\met$ to $T_\sigma(\Diffo\cdot\sigma)$ consists of variations satisfying 
\begin{equation}\label{E:nabla}
\nabla^\alpha\delta\sigma_{\alpha\beta}=0\ .
\end{equation}
  Restricting to hyperbolic metrics, it then follows from \eqref{E:hyptangent} that these variations must also be traceless.  Hence, we have an identification of $T^\ast_{[\sigma]}\metTeich$ with the space of traceless symmetric $2$-tensors satisfying \eqref{E:nabla}.
  Now the bundle of traceless symmetric $2$-tensors is real isomorphic to $K_S^2$ via
  $
 2 \varphi(z)=\delta\sigma_{11}-i\delta\sigma_{12}
  $.
  Moreover, \eqref{E:nabla} is precisely the statement that the corresponding quadratic differential $\varphi=\varphi(z)dz^2$ is holomorphic.  Hence,  $T^\ast_{[\sigma]}\Teich\simeq\QD(S)$.\index{Teichm\"uller!space!cotangent space}

This description makes contact with the Kodaira-Spencer theory of deformations of a complex structure
(cf.\ \cite{Kodaira}).  Indeed, infinitesimal deformations of a complex structure  are parametrized by smooth sections  $\mu$ of $\overline K_S\otimes K^{-1}_S$.  These are just (smooth) Beltrami differentials. 
Note that
there is a natural pairing between Beltrami differentials and holomorphic quadratic differentials on a Riemann surface $S$ obtained by raising indices in \eqref{E:met}:
\begin{equation} \label{E:pairing}
\langle \mu, \varphi\rangle= \int_S \mu(z)\varphi(z) |dz|^2\ ,
\end{equation}
where $\mu=\mu(z) d\bar z\otimes (dz)^{-1}$ and $\varphi=\varphi(z) dz^2$.
Let ${\mathcal H}{\mathcal B}(S)$  denote the space of \emph{harmonic} Beltrami differentials,\index{Beltrami!differential!harmonic} i.e.
$$
 {\mathcal H}{\mathcal B}(S)=\{\mu\in {\mathcal B}(S) : \bar\partial^\ast\mu=0\}\ ,
$$
where the adjoint $\bar\partial^\ast$ is defined with respect to the hyperbolic metric.
For any holomorphic quadratic differential $\varphi$, the Beltrami differential $\mu=\sigma^{-1}\bar\varphi$ is harmonic.  Moreover,  $\langle \mu,\varphi\rangle=\Vert\varphi\Vert^2_{2}$, where $\Vert\cdot\Vert_2$ denotes the $L^2$-norm with respect to the metric $\sigma$:
\begin{equation} \label{E:L2}
\Vert\varphi\Vert^2_{2}=\int_S |\varphi(z)|^2\sigma(z)^{-1} dxdy\ .
\end{equation}
  It follows that the pairing 
\begin{equation} \label{E:nondegenerate}
\langle\cdot,\cdot\rangle : {\mathcal H}{\mathcal B}(S)\times \QD(S)\lra \C
\end{equation}
is nondegenerate and that the tangent space is given by $T_{[\sigma]}\Teich={\mathcal H}{\mathcal B}(S)$.\index{Teichm\"uller!space!tangent space}

To complete this circle of ideas, one can directly compute the Beltrami differential associated to $\delta\sigma$. Let $\sigma_t$ be a differentiable family of hyperbolic metrics with $\sigma_0=\sigma$, $(d\sigma_t/dt)|_{t=0}=\delta\sigma$, and let $\nu_t$ be the Beltrami differentials associated to the identity map $(S,\sigma)\to (S,\sigma_t)$, $\nu_0=0$, $(d\nu_t/dt)|_{t=0}=\mu$.
 If $w=w_t$, $w_0=z$ is a family of conformal coordinates such that 
 $$
 ds_{\sigma_t}^2=\sigma_t(w)|dw|^2=\sigma_t(w)|w_z|^2|dz+\nu_t d\bar z|^2\ ,
 $$
 then since $\delta\sigma_{\alpha\beta}$ is traceless,
 $$
 \frac{d}{dt}\biggr|_{t=0} \sigma_t(w)|w_z|^2=0\ .
 $$
 It then follows that $2\sigma\mu=\delta\sigma_{11}+i\delta\sigma_{12}$, in agreement with the correspondence above.

There is a canonical map ${\mathfrak c}: \metTeich\to \Teich$ obtained by associating to an equivalence class of  hyperbolic metrics the corresponding equivalence classes of complex structures obtained via Gauss' theorem (see Section \ref{S:uniformization}). This map is continuous, for if two hyperbolic metrics are close in the smooth topology, then the identity has small dilatation. Furthermore, $\mathfrak c$ is a bijection by the uniformization theorem.   
With this understood, we now see that the two definitions of Teichm\"uller space are equivalent.
\begin{theorem} \label{T:metteich}
The canonical map ${\mathfrak c}: \metTeich\to \Teich$ obtained by associating to the hyperbolic metric its conformal class is a homeomorphism.
\end{theorem}

\begin{proof}
Recall from the proof of Theorem \ref{T:teich2} that the map $F:\Teich\to\Fricke$ is also a continuous bijection.  Since $\Fricke\subset\R^{6p-6}$, it follows by Proposition \ref{P:smooth} and Invariance of Domain that the composition 
$$F\circ{\mathfrak c}:\metTeich\to\Fricke\hookrightarrow\R^{6p-6}$$
 is a homeomorphism; hence, both $F$ and $\mathfrak c$ are as well.
\end{proof}

\begin{remark} \label{R:equivalence}
\begin{enumerate}
\item We emphasize that the proof of the homeomorphism $\metTeich\simeq \Teich$ given above is independent of  the Teichm\"uller Theorems \ref{T:existence}, \ref{T:uniqueness}, and \ref{T:teich}.
\item
By Theorem \ref{T:metteich}, we may regard the topological space $\Teich$ either as equivalence classes of marked Riemann surfaces or as the moduli space of hyperbolic metrics.    In particular, for the alternative proofs of Theorem \ref{T:teich} given below, it suffices to prove that $\metTeich$ is homeomorphic to a cell.
\item The $L^2$-metric \eqref{E:L2} is the \emph{Weil-Petersson cometric} on $\Teich$ (see Section \ref{S:wp} below).  In this description, it is easy to see that the Teichm\"uller metric \eqref{E:teichmetric} is a Finsler metric defined by the $L^1$-norm \eqref{E:L1}. 
\end{enumerate}
\end{remark}


\subsection{Harmonic Maps}


This section is a brief summary of the theory of harmonic maps with an emphasis on those aspects that relate to Teichm\"uller theory.  In Section \ref{S:def}, we give the basic definitions and present the variational formulation along with some examples.  In Section \ref{S:es}, we state the  existence and uniqueness theorem of Eells-Sampson-Hartman for nonpositively curved targets, and we indicate the importance of the Bochner formula.  In Section \ref{S:2d}, we specialize to the case of surface domains.  We discuss conformal invariance, the Hopf differential, and some applications.  In Section \ref{S:wolf}, we present another proof that Teichm\"uller space is a cell using harmonic maps. 


\subsubsection{Definitions.} \label{S:def}


Let $(M, g)$ and $(N, h)$ be Riemannian manifolds.  With respect to coordinates $\{x^\alpha\}$ on $M$ and $\{y^i\}$ on $N$, write $g=(g_{\alpha\beta})$, $h=(h_{ij})$.  Given a smooth map $f:M\to N$, its differential 
$$(df)^k_\alpha=(\partial f^k/\partial x^\alpha) dx^\alpha\otimes (\partial/\partial y^k)\ ,$$
 is a section of the bundle $T^\ast M\otimes f^\ast TN$ with the induced metric and connection. 
Define the \emph{energy density}\index{energy!density} and \emph{energy of} $f$ by\index{energy}
\begin{align} 
e(f)&=\tfrac{1}{2}\langle df, df\rangle_{T^\ast M\otimes f^\ast TN}=\frac{1}{2}\frac{\partial f^i}{\partial x^\alpha}\frac{\partial f^j}{\partial x^\beta} g^{\alpha\beta}h_{ij}\circ f\ ,
\label{E:ef} \\
E(f)&=\int_M e(f) dvol_M\ ,\label{E:Ef}
\end{align}
respectively (repeated indices are summed).  The energy can be viewed as a functional on the space of smooth maps between $M$ and $N$. 

The second  extremal problem,  analogous to the Teichm\"uller  problem in Section \ref{S:qc}, may now be formulated as follows:  given a smooth map $f:(M,g)\to(N,h)$, let 
\begin{equation} \label{E:extenergy}
E^\ast[f]=\inf\{ E(f') : f' \ \text{smooth}\ ,\ f'\sim f\}
\end{equation}

\medskip\indent
{\bf Energy Extremal Problem.} \index{energy!extremal problem}
Is $E^\ast[f]$ is realized as the energy  of a smooth map, and if so, what are the properties of the map?

\medskip
A smooth map $f$ such that $E(f)=E^\ast[f]$ is called an \emph{energy minimizer}.\index{energy!minimizer}
Unlike the problem for quasiconformal maps, existence of energy minimizers is not obvious.  We will discuss this at greater length in the next section.

   The covariant derivative  $\nabla df$
is a section of $\Sym^2(T^\ast M)\otimes f^\ast TN$, where $\Sym^2$ denotes symmetric $2$-tensors.
 The trace  $\tau(f)=\Tr_g \nabla df$ is called  the \emph{tension field of $f$}.  \index{tension field}
Let $\Delta$ denote the Laplace-Beltrami operator on $(M,g)$.  Then
$$
\tau(f)^k=\Delta f^k+(\Gamma^k_{ij}\circ f)\frac{\partial f^i}{\partial x^\alpha}\frac{\partial f^j}{\partial x^\beta}g^{\alpha\beta}\ .
$$
 Here, $\Gamma^k_{ij}$ denotes the Christoffel symbols of  $N$.
A smooth map $f:M\to N$ is called \emph{harmonic}\index{harmonic map} if $\tau(f)\equiv 0$.  Let
$$
d_\nabla : \Omega^p(f^\ast TN)\lra \Omega^{p+1}(f^\ast TN)
$$
denote the exterior derivative coupled with pulled-back Levi-Civit\`a connection on $N$.  It is easily seen that $d_\nabla(df)=0$ for all differentiable maps. The equations for harmonic maps are then equivalent to 
\begin{equation} \label{E:harmonic}
d_\nabla(\ast df)=0\ ,
\end{equation}
i.e.\ $df$ is a harmonic form (cf.\ \cite{ES, ES2}).  Here are some examples:

\begin{itemize}
\item Harmonic maps $S^1\to N$ are closed geodesics in $N$;
\item When $N=\R^n$ the harmonic map equations are equivalent to the harmonicity of the coordinate functions.
\item Totally geodesic maps satisfy $\nabla df=0$, and so are harmonic.
\item Holomorphic or anti-holomorphic maps between K\"ahler manifolds are harmonic. \index{energy!and holomorphic maps}
\item Minimal isometric immersions are harmonic.

\end{itemize}

 Now let us consider variational formulas for the energy $E(f)$. 
A smooth vector field $v$ along $f$, i.e. $v\in C^\infty(f^\ast TN)$, defines a variation of $f$ by $f_t(x)=\exp_{f(x)}(tv(x))$.  Since $N$ is assumed to be complete, this defines a smooth map $M\times\R\to N$ with $f_0=f$.  The \emph{first variational formula} is\index{energy!first variation of}
\begin{equation} \label{E:firstvar}
\delta_v E(f)=\frac{dE(f_t)}{dt}\biggr|_{t=0}=-\int_M \langle\tau(f),v\rangle_h dvol_M\ .
\end{equation}
It follows that the Euler-Lagrange equations for $E$ are precisely the harmonic map equations \eqref{E:harmonic}.  

In general there is a distinction between energy minimizers, smooth minimizers of $E$ which then necessarily satisfy \eqref{E:harmonic}, and smooth solutions to \eqref{E:harmonic} which may represent higher critical points of the energy functional.  We shall see below that this distinction vanishes when the target manifold $N$ has nonpositive curvature.  Another case where minimizers can be detected is the following:
let $S$ be a compact Riemann surface and
 $N$  a compact K\"ahler manifold. \index{energy!and holomorphic maps}
 \begin{proposition} \label{P:holo}
  If
 $f: S\to N$ is holomorphic or anti-holomorphic, then for any conformal metric on $S$,  $f$ is harmonic and is energy minimizing in its homotopy class. 
 \end{proposition}
  Indeed, a computation in local coordinates as above shows that for any smooth map $f:S\to N$, 
 \begin{align} 
 E(f)&=\int_S f^\ast\omega + 2 \int_S|\bar\partial f|^2 dvol_S\label{E:antiholo} \\
   &=-\int_S f^\ast\omega + 2 \int_S|\partial f|^2 dvol_S \label{E:holo}
 \end{align}
 where $\omega$ is the K\"ahler form on $N$.  Since the first terms on the right hand sides depend only on the homotopy class of $f$, the result follows.

Now let $v,w\in C^\infty(f^\ast TN)$ and $f_{s,t}$ be a two-parameter family of maps such that $f_{0,0}=f$, $v=(\partial f_{s,t}/\partial s)|_{s=t=0}$,  $w=(\partial f_{s,t}/\partial t)|_{(s,t)=(0,0)}$, where $f$ is harmonic.  Then\index{energy!second variation of}
\begin{equation} \label{E:secondvar}
H_f(v,w)=\frac{\partial^2E(f_{s,t})}{\partial s\partial t}\bigr|_{s=t=0}=-\int_M \langle J_fv,w\rangle_h dvol_M\ ,
\end{equation}
where
\begin{equation} \label{E:J}
J_f(v)=\Tr_g(\nabla^2 v+\Riem^N(df,v)df)\,
\end{equation}
is the Jacobi operator, and $\Riem^N$ is the Riemannian curvature of $(N,h)$.
  In particular, if $N$ has nonpositive Riemannian sectional curvature, then 
$$
H_f(v,v)\geq \int_M|\nabla v|^2 dvol_M\geq 0\ ,
$$
and hence every harmonic map is a local minimum of the energy.

Given smooth maps $f:M\to N$ and $\psi:N\to P$, one has the composition formula
$$
\nabla d(\psi\circ f)=d\psi\circ \nabla df+\nabla d\psi(df,df)\ .
$$
Taking traces we obtain the formula for the tension (cf.\ \cite{ES})
\begin{equation} \label{E:composition}
\tau(\psi\circ f)=d\psi\circ\tau(f)+\Tr_g\nabla d\psi(df,df)\ .
\end{equation}
In particular, if $f$ is harmonic and $\psi$ is totally geodesic then $\psi\circ f$ is also harmonic.  If $P=\R$ and $f$ is harmonic, then \eqref{E:composition} becomes 
$$\Delta(\psi\circ f)=\Tr_g\nabla d\psi(df,df)\ ,
$$
 and therefore a harmonic map pulls back germs of convex functions to germs of subharmonic functions.  The converse is also true:
\begin{theorem}[Ishihara \cite{Ish}] \label{T:ishihara}\index{Ishihara's theorem}
A map is harmonic if and only if it pulls back germs of convex functions to germs of subharmonic functions.
\end{theorem}


\subsubsection{Existence and uniqueness.} \label{S:es}


In the case of nonpositively curved targets the energy extremal problem has a solution.  
The basic existence result  is the following
\begin{theorem}[Eells-Sampson  \cite{ES}]\label{T:es} \index{harmonic map!existence of} \index{Eells-Sampson theorem}
Let $f:M\to N$ be a continuous map between compact Riemannian manifolds, and suppose that $N$ has nonpositive sectional curvature.  Then there exists an energy minimizing  harmonic map  homotopic to $f$.
\end{theorem}

The proof is based on the heat equation method  to deform a map to a harmonic one (cf.\ \cite{Ham2}).  Namely, one solves the initial value problem for a nonlinear parabolic equation
\begin{equation} \label{E:heat}
\frac{\partial f}{\partial t}(x,t)= \tau(f)(x,t)\ ,\ f(x,0)=f(x)\ .
\end{equation}
Stationary solutions to \eqref{E:heat} satisfy the harmonic map equations.  Furthermore, by taking the inner product on both sides in \eqref{E:heat} with $\tau(f)(x,t)$ and integrating over $M$, one observes, using \eqref{E:firstvar}, that the energy of the map $x\mapsto f(x,t)$ is decreasing in $t$.  Hence, one hopes that as $t\to \infty$, $f(\cdot, t)$ converges to a harmonic map. Unfortunately, it turns out that this procedure does not always work.   In general, even existence of a solution to \eqref{E:heat} for all $t\geq 0$ is not guaranteed (cf.\ \cite{CDY, Coron}). However, we have

\begin{theorem} \label{T:convergence}
Assume $M, N$ are compact Riemannian manifolds and $N$ has nonpositive sectional curvature.  Given a smooth map $f:M\to N$, then the solution to \eqref{E:heat} exists for all $t\in [0,\infty)$ and converges as $t\to\infty$ uniformly to a harmonic map homotopic to $f$.
\end{theorem}

The key to this theorem is following parabolic \emph{Bochner formula}.\index{Bochner formula}  Suppose $f(x,t)$ is a solution to \eqref{E:heat} for $0\leq t<T$, and let $e(f)(x,t)$ denote the energy density of the map $x\mapsto f(x,t)$.
Then for any orthonormal frame $\{u_\alpha\}$ at a point $x\in M$ we have the following pointwise identity:
 \begin{align} 
 -\frac{\partial e(f)}{\partial t}+ \Delta e(f)&=|\nabla df|^2+\langle df \Ric^M(u_\alpha), df(u_\alpha)\rangle
 \label{E:bochner} \\
 &\qquad -\langle \Riem^N(df(u_\alpha), df(u_\beta)) df(u_\beta), df(u_\alpha)\rangle\ . \notag
 \end{align}
 where $\Ric^M$ is the Ricci curvature of $(M,g)$.  In particular, if $M$ is compact and $N$ is nonpositively curved, then
 \begin{equation} \label{E:key}
 \frac{\partial e(f)}{\partial t}\leq \Delta e(f)+ C e(f)
 \end{equation}
 for some constant $C\geq 0$.
 If the solution $f(x,t)$ exists for $0\leq t< T$, then it follows easily from \eqref{E:key} that $e(f)$ is uniformly bounded in $x$ and $t<T$.  The bound on the energy density means that the maps $f(\cdot, t)$ form an equicontinuous family from which convergence as $t\to T$ can be deduced.  Resolving the initial problem at $t=T$ then allows one to extend the solution for some small time $t>T$.
 This is the rough idea behind the existence for all $0\leq t<\infty$.  In fact, more sophisticated methods show that $e(f)$ is bounded for all time (cf.\ \cite{Moser}).
 
 The Bochner formula for harmonic maps, i.e.\ stationary solutions of \eqref{E:heat}, \index{Bochner formula} is
  \begin{align} 
 \Delta e(f)&=|\nabla df|^2+\langle df \Ric^M(u_\alpha), df(u_\alpha)\rangle \label{E:bochner2}
  \\
 &\qquad -\langle \Riem^N(df(u_\alpha), df(u_\beta)) df(u_\beta), df(u_\alpha)\rangle\ . \notag
 \end{align}
As before, this implies
 \begin{equation} \label{E:key2}
  \Delta e(f)\geq - C e(f)\ .
 \end{equation}
Inequality \eqref{E:key2} is the key to regularity of weakly harmonic maps to nonpositively curved spaces (cf.\ \cite{Schoen1}).
To state the result precisely, we note that if $\Omega\subset M$ is a domain with smooth boundary, one can solve the Dirichlet problem for an energy minimizing map $f:\Omega\to N$ with prescribed boundary conditions.  If $f:M\to N$ is energy minimizing then it is automatically energy minimizing with respect to its boundary values for any $\Omega\subset M$.  This is what is meant by \emph{locally energy minimizing}.  The following Lipschitz bound  follows from \eqref{E:key2} by  iterating the Sobolev embedding.  
\begin{proposition} \label{P:lipschitz}
If $f:\Omega\to N$ is harmonic with energy $E(f)$ and $N$ has nonpositive curvature, then for any $U\subset\subset\Omega$,
$$\sup_{x\in U} e(f)(x)\leq C(U) E(f)$$
for some constant $C(U)$ independent of $f$.
\end{proposition}

 Next, we have the following result on uniqueness.

\begin{theorem}[Hartman \cite{Ht}] \label{T:homotopy} \index{harmonic map!uniqueness of}
Assume $M, N$ are compact Riemannian manifolds and $N$ has nonpositive sectional curvature.   Let $f_0, f_1 : M\to N$ be homotopic harmonic maps, and let $f_s:M\to N$ be a geodesic homotopy where $s\in [0,1]$  is proportional to arc length.  Then:
\begin{enumerate}
\item  for every $s$, $f_s$ is a harmonic map with $E(f_s)=E(f_0)=E(f_1)$; and
\item  the length of the geodesic $s\mapsto f_s(x)$ is independent of $x$.
\end{enumerate}
In case $N$ has \emph{negative} sectional curvature,
 any nonconstant harmonic map $f:M\to N$ is unique in its homotopy class \emph{unless} $f$ maps onto a geodesic,  in which case all homotopic harmonic maps are translations of $f$ along the geodesic.\end{theorem}

\noindent
A good reference for these results  is Jost's book \cite{J1}.

The  theorems above  fail to hold if the curvature assumption on $N$ is relaxed.  In this case the analytic complexity increases substantially, and there is no satisfactory existence result   in general.  There is something of an exception in the case of surface domains (see \cite{J3, SU1, SU2, SY2}), where the conformal invariance with respect to the domain metric leads to bubbling phenomena.
 We will not attempt to present any results for the case of higher dimensional domains, since the relation with Teichm\"uller theory is less important.

 \subsubsection{Two dimensional domains.} \label{S:2d}
 
 We now specialize to the case where the domain  is a Riemann surface.
   Here the salient feature, as we have just mentioned above, is that the energy functional is
invariant under conformal changes of metric on $S$, i.e.\ $g\mapsto e^\phi g$.  Hence, the harmonic map equations for surface domains depend only on the complex structure on $S$.

Let $f:(S,\sigma)\to (N,h)$ be a smooth map, where $N$ is an arbitrary Riemannian manifold.  Then $\varphi=(f^\ast h)^{2,0}=\Hopf(f)$ is a quadratic differential, called the \emph{Hopf differential of $f$}.\index{Hopf differential}  A key fact is that \emph{$\varphi$ is holomorphic  if $f$ is harmonic}.  Indeed, in local coordinates, $\varphi=\varphi(z)dz^2$, where
\begin{equation} \label{E:hopf}
\varphi(z)=\left\langle f_z, f_z\right\rangle=\tfrac{1}{4}\left(
|f_x|^2-|f_y|^2-2i\left\langle f_x, f_y\right\rangle
\right)\ ,
\end{equation}
Notice that $\varphi\equiv 0$ if and only if $f$ is conformal.
In normal coordinates at $f(z)$, the harmonic map equations are $\Delta f^k=0$, for all $k$.  Together with the vanishing of the derivatives of the metric, this implies
$$
\varphi_{\bar z}(z)=\left\langle f_z,f_z\right\rangle_{\bar z}=\left(h_{ij}(f(z)) f^i_z f^j_z\right)_{\bar z}
=2h_{ij}(f(z)) f^i_z f^j_{z\bar z}=0\ .
$$
  Several results in this article depend on the holomorphicity of the Hopf differential.  In Section \ref{S:npc}, we will present a different argument due to Schoen \cite{Schoen2} which works for a more general class of metric space targets.

To see how holomorphicity can have topological consequences, take for example the case where the target  is also a Riemann surface $R$.  Writing the metric $h$ on $R$ in local conformal coordinates $w$,
the energy  of a map $f$  is then
\begin{equation}\label{E:energy3}
E_h(f)
= \int_S h(f(z))(|f_z|^2+|f_{\bar z}|^2)dxdy\ , \
z=x+iy
\end{equation}
where we have confused the notation $f$ and $w\circ f$.  When the metric $h$ is understood, we shall simply write $E(f)$.
The harmonic map equations are  (cf.\ \cite[Ch.\ 1]{SY1})
\begin{equation} \label{E:surface}
f_{z\bar z}+\frac{h_w}{h} f_z f_{\bar z}=0\ .
\end{equation}
As an immediate application of \eqref{E:surface} it follows that if $f:S\to R$ is harmonic, then $|\partial f|$ and $|\bar \partial f|$ are either identically zero or have a well-defined order.  Indeed, if $H=f_z$ and $G=-(h_w/h)f_{\bar z}$, and $\zeta$ satisfies the equation $\zeta_{\bar z}=-G$, then it is easily checked that $He^\zeta$ is holomorphic.  By setting $n_p=\ord_p He^\zeta$ we obtain
$|\partial f|=|z|^{n_p}k(z)$, where $k(z)$ is a smooth strictly positive function.  We call $n_p$ the \emph{order} of $|\partial f|$ at $p$.  This leads to

\begin{theorem}[Eells-Wood \cite{EW}] \label{T:eellswood}
Let $f:S\to R$ be a harmonic map between surfaces.  If $|\partial f|$ is not identically zero, then
$$
\sum_{|\partial f|(p)=0}n_p=\deg(f)\chi_R-\chi_S\ .
$$
If $|\bar\partial u|$ is not identically zero, then
$$
\sum_{|\partial f|(p)=0}m_p=-\deg(f)\chi_R-\chi_S\ .
$$
\end{theorem}
Here, $n_p$ and $m_p$ are the orders of $|\partial f|$, $|\bar\partial f|$ at $p$, respectively.
An immediate consequence of this is  Kneser's Theorem:
\begin{corollary}[Kneser, \cite{Kn}] \label{C:kneser}
Let $f: S\to R$ be a continuous map between surfaces, $\chi_R <0$.  Then $|\deg(f)|\chi_R\geq \chi_S$.
\end{corollary}

Pushing these ideas further, Schoen-Yau and Sampson proved
\begin{theorem}[Schoen-Yau \cite{SY}, Sampson \cite{Sa1}, see also Jost-Schoen \cite{JS}] \label{T:jostschoen}
Suppose $f:S\to R$ is a harmonic map between surfaces of the same genus.  If $\deg f=1$ and $R$ has negative curvature, then $f$ is a diffeomorphism.
\end{theorem}

Theorems \ref{T:eellswood} and \ref{T:jostschoen} depend on the following formulas for a harmonic map between surfaces.
\begin{equation}
\Delta \log|\partial f| = -K_R J(f)+K_S\, ,\
\Delta \log|\bar\partial f| = K_R J(f)+K_S\, , \label{E:du}
\end{equation}
where $K_S$, $K_R$ are the Gaussian curvatures of $S$ and $R$, and
$J(f)=|\partial f|^2-|\bar\partial f|^2$ is the Jacobian of $f$.
The equations \eqref{E:du}  are related to the
Bochner formula \eqref{E:bochner2}.  The proof is a simple
calculation which can be found, for example, in \cite[Ch.\ 1]{SY1}.

As an application, the next theorem regarding the quotient \eqref{E:metteich} is due to Earle and Eells (cf.\ \cite{EE} and also \cite{Ebin}).
 \begin{theorem} \label{T:ee}
 The bundle given by the quotient map $p:\hypmet\to \metTeich$ is trivial, i.e.\ there exists a homeomorphism $H$ with the property that the diagram
 \begin{displaymath}
 \xymatrix{ \hypmet \ar[dr]_{p} \ar[r]^{H\qquad } & \ar[d]^{\pi} \metTeich \times\Diffo   \\
  & \metTeich  }
 \end{displaymath}
 commutes, where $\pi$ is the projection onto the first factor.
 \end{theorem}
 The map $H$ can be constructed as follows:  fix a metric
 $\sigma_0\in\hypmet$.  For any other $\sigma\in \hypmet$, let $f_\sigma :
 (S,\sigma_0)\to (S,\sigma)$ be the harmonic diffeomorphism $\sim \id$ from
 Theorem \ref{T:jostschoen}.   Then $H$ is defined by $F(\sigma)=(p(\sigma),
 f_\sigma^{-1})$.

\subsubsection{A second proof of Teichm\"uller's theorem.} \label{S:wolf}

We now give a second proof that Teichm\"uller space is a cell (Theorem \ref{T:teich}) using harmonic maps and Hopf differentials as opposed to Teichm\"uller maps.  \index{Teichm\"uller!theorem}
  Let $\varphi_\sigma=\Hopf(f_\sigma)$ be the Hopf differential of the map $f_\sigma$ defined above.  By uniqueness of the harmonic diffeomorphism in its homotopy class (Theorem \ref{T:homotopy}) we obtain a well-defined map
\begin{equation} \label{E:wolfmap}
{\mathcal H} : \metTeich\lra\QD(S) : [\sigma]\mapsto {\mathcal H}[\sigma]=\varphi_\sigma\ .
\end{equation}
Then we have
\begin{theorem}[Wolf \cite{Wf1}] \label{T:wolf}
The map $\mathcal H$ is a diffeomorphism.
\end{theorem}

The fact that $\mathcal H$ is 1-1 is due to Sampson \cite{Sa1}.  The smooth dependence of $\mathcal H$ follows easily as in \cite{ES}.  This seems to have been first observed also by Sampson.  That $\mathcal H$ is proper is due to Wolf.  
The idea is based on the following energy bound (see also \cite{Minsky}):
\begin{equation} \label{E:wolf}
E(f_\sigma)\leq 2\int_S |\varphi_\sigma| - 2\pi\chi_S\ .
\end{equation}
To see this, let $f: (S,\sigma_0)\to (S,\sigma)$ be \emph{any} quasiconformal map
with Beltrami coefficient $\mu$ and Hopf differential $\varphi$.  Then
$$
|\bar\partial f|^2 dvol=\sigma f_{\bar z} \bar f_z|dz|^2= \sigma  f_z \bar f_z\frac{ f_{\bar z}}{ f_z}|dz|^2=\varphi\mu|dz|^2\leq |\varphi| \ ,
$$
since $|\mu|<1$.
Then since $f_\sigma$ has degree $1$,  \eqref{E:wolf} is a consequence of the  above inequality and \eqref{E:antiholo}.   
Similarly, using \eqref{E:holo}, one has by the same argument
\begin{equation} \label{E:wolf2}
 2\int_S |\varphi_\sigma| + 2\pi\chi_S\leq E(f_\sigma)\ .
\end{equation}
 Properness of $\mathcal H$ now follows from \eqref{E:wolf} and properness of the energy.  The latter is due to Schoen and Yau \cite{SY2}. In 
Section \ref{S:proper},  we will sketch the proof.  Finally, that $\mathcal H$ is onto follows from the properness and the fact that both the domain and target are smooth manifolds of dimension $6p-6$.
Hence, Theorem \ref{T:teich} is a consequence of the theorem above, along with Theorem \ref{T:metteich}.


\subsection{Singular Space Targets}


Harmonic maps to singular spaces were first introduced in a systematic 
way in the paper of Gromov and Schoen \cite{GS}  in connection 
with arithmetic superrigidity.  Since then the subject has played an 
important role in Teichm\"uller theory and is one of the main themes of 
this review.  In Section \ref{S:gr}, we will indicate how singular 
space targets make a  connection between the extremal maps discussed in 
Section \ref{S:teich} and the harmonic maps of Section \ref{S:def}.  The highlight is the proof of Teichm\"uller's existence and uniqueness theorems.  The idea, going back 
to Gerstenhaber and Rauch, provides a clear motivation for the use 
of singular targets from the point of view of Teichm\"uller theory.  We will 
defer the technical aspects of the general theory to Section \ref{S:npc}.  In Section
\ref{S:trees}, we 
discuss the notion of $\R$-trees and their connection to measured foliations 
and quadratic differentials.  We also state the famous Hubbard-Masur Theorem.  Section \ref{S:npc} contains all of the 
technical results on harmonic maps to metric spaces that we will 
need in this article.  There we give an outline of the main results of 
\cite{GS, KS1, KS2}.  In addition, we describe several results that are special to harmonic maps to trees.

\subsubsection[The Gerstenhaber-Rauch approach.]{The Gerstenhaber-Rauch approach to Teichm\"uller's extremal problem.} \label{S:gr}


Teichm\"uller's extremal problem (Section \ref{S:qc}) and the energy extremal problem (Section \ref{S:def}) bear obvious similarities; hence,  the natural

\medskip\noindent {\bf Question.}
Are Teichm\"uller maps  harmonic
for some metric?
\medskip

This leads to the notion of energy minimizing maps
to singular space targets, which is the subject of this section.
We begin with a simple example.
Given a holomorphic quadratic differential $\psi$ on a Riemann surface $R$,  $|\psi|$ defines a singular flat metric  with conical singularities at the zeros of $\psi$ (cf.\ \cite{Str}).  Indeed,
away from the zeros we may write $|\psi|=|dw|^2$ for some conformal
coordinate $w$, whereas at  a zero of order $m\geq 1$,
$|\psi|=|w|^{m+2/2}|dw|^2$. Notice that for $h(w)=|w|^{m+2/2}$,
the Gauss curvature
\begin{equation}\label{E:curvature}
K=-\frac{1}{2h}\Delta \log h\leq 0\ ,
\end{equation}
in the sense of distributions.  We say that $S$ with the metric
$h=|\psi|$ is a \emph{nonpositively curved space}.  \index{NPC}

Let $S$ be another Riemann surface.
Given a map $f:S\to R$ one can define the Sobolev class $H^1$ and the
energy of $f$ with respect to the singular conformal metric $|\psi|$ on $R$
by \eqref{E:energy3}.  Following the definitions of Section \ref{S:def} we call such a map harmonic if it is an energy minimizer (see also Theorem \ref{T:homotopy}).
This is a special case of the general theory of
Gromov, Korevaar, and Schoen that we will describe below;
in particular, such minimizers always exist and are Lipschitz by
Theorem \ref{T:ks1}.  The following result builds on 
earlier, weaker versions  due to Miyahara \cite{Miy} and Leite \cite{Leite}.

\begin{theorem}[Kuwert \cite{Ku}] \label{T:k1}
A Teichm\"uller map $f_0:S\to R$ is the unique harmonic map in its
homotopy class when $R$ is endowed with the singular flat metric
$h=|\psi|$ defined by the terminal quadratic differential of $f_0$.
\end{theorem}

Let us show how this gives a 
\begin{proof}[Proof of Teichm\"uller's Uniqueness Theorem \ref{T:uniqueness}]\index{Teichm\"uller!uniqueness theorem}
Let $f:S\to R$ be any quasiconformal map with Beltrami differential $\mu_f$.  Then by \eqref{E:energy3} we have
\begin{align*}
E(f)&=\int_S(|f_z|^2+|f_{\bar z}|^2)|\psi(f(z))|dxdy\\
&=\int_S(1+|\mu_f|^2)|\psi(f(z))||f_z|^2 dxdy\\
&\leq (1+\Vert \mu_f\Vert_\infty^2)\int_S|\partial f|^2 dvol_S\ .
\end{align*}
Now by \eqref{E:holo}, which continues to hold for the singular metric, 
\begin{align}
E(f)&\leq \frac{1}{2}(1+\Vert \mu_f\Vert_\infty^2)(E(f)+ C[f])\notag \\
E(f)&\leq \frac{1+\Vert \mu_f\Vert_\infty^2}{1-\Vert \mu_f\Vert_\infty^2}C[f]\ , \label{k1}
\end{align}
where $C[f]$ is a constant depending only on the homotopy class of $[f]$ and the area of the metric $|\psi|$.  On the other hand, for the Teichm\"uller map $f_0$ we have by the same computation
\begin{equation}
E(f_0)= \frac{1+\Vert \mu_{f_0}\Vert_\infty^2}{1-\Vert \mu_{f_0}\Vert_\infty^2}C[f_0]\ . \label{k2}
\end{equation}
If $f\sim f_0$, then $C[f]=C[f_0]$.
By Theorem \ref{T:k1}, $E(f_0)\leq E(f)$, which by \eqref{k1} and \eqref{k2} implies $\Vert \mu_{f_0}\Vert_\infty\leq \Vert \mu_{f}\Vert_\infty$, with equality if and only if $f=f_0$.
\end{proof}

This result does not answer the question of existence of extremal
maps by harmonic map methods.  In their 1954 paper, Gerstenhaber and
Rauch proposed a minimax method of finding a Teichm\"uller map \cite{GRauch}.  Let $\CM(R)$
denote the space of conformal metrics on $R$ with unit area and with
at most conical singularities (see below for more details).  For
each $h\in {\CM(R)}$, let $E_h(f)$ be defined as in \eqref{E:energy3}, where $f:S\to R$ is in $H^1$.
Gerstenhaber-Rauch conjectured that
\begin{equation} \label{E:gr} \index{Gerstenhaber-Rauch conjecture}
\sup_{h\in{\CM(R)}}\inf_{f\sim f_1} E_h(f)=\frac{1}{2}\left(
K^\ast[f_1]+\frac{1}{K^\ast[f_1]}\right)\ ,
\end{equation}
and that the sup-inf in \eqref{E:gr} is realized by the
Teichm\"uller map homotopic to $f_1$.  The problem was investigated further by Reich
and Reich-Strebel in the case where $S,R$ are both the disk
\cite{Re,ReSt}.  Kuwert, assuming the existence of the Teichm\"uller
map, proved 
\begin{theorem}[Kuwert \cite{Ku}] \label{T:k2}
The Teichm\"uller map $f_0$ and the singular metric $h_0=|\psi|$
defined by its terminal differential realize the sup-inf in
\eqref{E:gr}.
\end{theorem}
The full Gerstenhaber-Rauch conjecture was recently proved by Mese (cf.\
\cite{Me3,Me4}) using the harmonic map theory of Gromov, Korevaar,
and Schoen. Before we state Mese's theorem we need to set up some
notation and terminology.
Let $(X,d)$ be a metric space which is also a \emph{length space},\index{length space}
i.e.\ for all pairs $p, q\in X$ there exists a rectifiable curve
$\gamma_{pq}$ whose length equals $d(p,q)$ (which we sometimes write
$d_{pq}$).  We call $\gamma_{pq}$ a geodesic from $p$ to $q$.  Then
$X$ is \emph{NPC} (= nonpositively curved)\index{NPC} if every point of $X$ is
contained in a neighborhood $U$ so that for all $p,q,r\in U$,
$$
d^2_{pq_\tau}\leq (1-\tau)d^2_{pq}+\tau d^2_{pr}
-\tau(1-\tau)d^2_{qr}\ ,
$$
where $q_\tau$ is the point on $\gamma_{qr}$ so that
$d_{qq_\tau}=\tau d_{qr}$.  Note that equality is achieved for every
triple $p,q,r\in \R^2$.  More generally, one defines a length space with curvature bounded above by $\kappa$ by making comparisons with geodesic triangles in surfaces of constant curvature $\kappa$ (cf.\ \cite{Me2}).  It follows
from \eqref{E:curvature} that if $h(w)|dw|^2$ is a conformal
metric on $R$ with 
\begin{equation} \label{E:curvbound}
\Delta\log h \geq -2\kappa h
\end{equation}
then the induced metric space has curvature bounded above by
$\kappa$ (cf. \cite{Me1}).  We will use this fact when we give a
harmonic map construction of the Teichm\"uller map.

Let ${\CM}_{\alpha,\kappa}(R)$ denote the set of
metrics $h=h(w)|dw|^2$ on $R$ where $h\geq 0$ is bounded
of Sobolev class $H^1$, satisfies \eqref{E:curvbound} weakly, and has area
$=\alpha$.  Let $d_h$ denote the distance function associated to the
above metric.  As we have discussed before it is not hard to see
that $(R, d_h)$ has curvature bounded above by $\kappa$.  The key
result  is the following
\begin{theorem}[Mese \cite{Me4}] \label{T:mese}
Let $h_i\in {\CM}_{\alpha,\kappa}(R)$, $\kappa>0$, and $f_i:S\to (R,h_i)$ be
such that 
\begin{enumerate}
\item $f_i$ is harmonic;
\item $\displaystyle \lim_{i\to \infty} E_{h_i}(f_i)=\sup_{h\in {\CM}_{\alpha,\kappa}(R)}\inf_{f\sim f_1} E_h(f)$.
\end{enumerate}
Then  the $f_i$ converge in the pullback sense to the
Teichm\"uller map $f_0$.
\end{theorem}

Convergence in the pullback sense is essentially Gromov-Hausdorff convergence.   This will be explained in greater detail below (see Section \ref{S:ms}).  Theorem \ref{T:mese}, along with earlier work, gives a proof of Teichm\"uller's Existence Theorem \ref{T:existence}.\index{Teichm\"uller!existence theorem}

\subsubsection{$\R$-trees.} \label{S:trees}

The use of singular metrics to prove the Teichm\"uller theorems is motivation to study energy minimizing maps for other metric space targets.  Here we discuss another ubiquitous example.
An  \emph{$\R$-tree} \index{tree}is a length space such that any
two points can be joined by a \emph{unique} path parametrized by arc
length.  This path is called \emph{the geodesic} between the points, say
$p,q$, and it is denoted $pq$.  An equivalent definition is that an
$\R$-tree is a simply connected length space with curvature bounded above by $\kappa$
for any $\kappa\in \R$ (cf.\ \cite{Sun}).

\begin{example}
Let $T$ be a simplicial tree, i.e. a simply connected 1-dimensional 
simplicial complex.  Then $T$ can be thought of as an $\R$-tree by
assigning to each edge a unit length.  An $\R$-tree is called
simplicial if it is obtained from a simplicial tree in this way.
Note that we do not assume the simplicial tree is locally finite, although the set of vertices clearly is.
\end{example}

\begin{example}
Take $T=\R^2$ and define $d(p,q)=|p-q|$ if $p,q$ lie on some ray
from the origin, and $d(p,q)=|p|+|q|$, otherwise.  Clearly, $T$
with this metric is not locally compact, though it is simplicial.
\end{example}

\begin{example}
A slight modification of the above yields a non-simplicial tree.
Again take $T=\R^2$ and define $d(p,q)=|p-q|$ if $p$ and $q$ lie on
the same vertical line.  In all other cases, let
$d(p,q)=d(p,p')+d(p',q')+d(q,q')$, where $p', q'$ are the
projections of $p,q$ to the $x$-axis. Then every point on the $x$-axis becomes a vertex.
\end{example}

$\R$-trees appear in Teichm\"uller theory in several ways.  The primary example is the leaf space of the horizontal and vertical foliations of a
holomorphic quadratic differential.  First recall that a
\emph{measured foliation}\index{measured foliation} $\mathcal F$ on a surface $S$ with
singularities at the points $z_1,\ldots, z_\ell$ and multiplicities
$k_1,\ldots, k_\ell$ is described by the following (cf.\ \cite{FLP}):  an open cover
$\{U_i\}$ of $S\setminus\{z_1,\ldots, z_\ell\}$ and open sets
$V_1,\ldots, V_\ell$ about $z_1,\ldots, z_\ell$ along with smooth
real valued functions $u_i$ defined on $U_i$ such that 
\begin{enumerate}
\item $|du_i|=|du_j|$ on $U_i\cap U_j$;
\item $|du_i|=|\imag(z-z_j)^{k_j/2} dz| $ on $U_i\cap V_j$.
\end{enumerate}
Clearly, $\ker du_i$ defines a vector field on $S$ which integrates
to give a foliation away from $\{z_1,\ldots, z_\ell\}$, with
$(k_j+2)$-pronged singularities at $z_j$ (see Figure 1).  A leaf containing a singularity is called a \emph{critical trajectory}\index{critical trajectory}, whereas the other leaves are called \emph{noncritical}.
An
important  attribute of measured foliations is that they carry a 
\emph{transverse measure}.  More precisely, if $c$ is a
rectifiable path then we denote by $\nu(c)$,
 the number $\nu(c)=\int_c |du|$, where $|du|$ is defined
by $|du|_{U_i}=|du_i|$.  An important feature of this measure is its translation invariance
along the leaves.  Namely, if $c_0$ is a path transverse to the foliation $\mathcal F$, and 
if we deform $c_0$ to $c_1$ via an isotopy that maintains the transversality to the foliation at every time, then $\nu(c_1)=\nu(c_0)$.  For the free homotopy class $[c]$ of a simple closed curve we define
\begin{equation} \label{E:intersection}
i([c],{\mathcal F})=\inf\{ \nu(c) : c\in[c]\}\ .
\end{equation}
Two measured foliations $({\mathcal F}, \nu)$ and $({\mathcal F}',\nu')$ are called equivalent if $i([c], {\mathcal F})=i([c], {\mathcal F}')$ for all free homotopy classes of simple closed curves.  We denote the space of equivalence classes of measured foliations on $S$ by $\MF$. Then the collection of intersection numbers \eqref{E:intersection}, as $c$ ranges over isotopy classes of simple closed curves, endows $\MF$ with a topology.
 We call ${\mathcal F}$ and ${\mathcal F}'$ \emph{projectively equivalent}\index{measured foliation!projective equivalence} if there is $b>0$ such that $i([c], {\mathcal F}')=b\, i([c], {\mathcal F})$ for all free homotopy classes of simple closed curves. In this case, we write ${\mathcal F}'=b{\mathcal F}$.  The space of projective equivalence classes will be denoted $\PMF$.

 Given a measured foliation $({\mathcal F},\nu)$  we can associate a \emph{dual tree} $T_{\mathcal F}$ \index{tree!dual to a foliation}to the foliation with an isometric action of $\Gamma=\pi_1(S)$.    Explicitly, let $(\widetilde{\mathcal F}, \tilde\nu)$ denote the pullback of $({\mathcal F},\nu)$ to the universal cover $\HH$ of $S$.  On $\HH$ we define a pseudodistance $\tilde d$ via 
$$
\tilde d(p,q)=\inf\{ \tilde\nu(c) : c \ \text{a rectifiable path between}\ p, q\}\ .
$$
It follows by \cite[Corollary 2.6]{Bow} that the Hausdorffication of $(\HH, \tilde d)$ is an $\R$-tree with an isometric action of $\Gamma$.  Strictly speaking, the setup in \cite{Bow} works for measured foliations on arbitrary $2$-complexes.  The approach is useful in that it avoids introducing the notion of a geodesic lamination.  For a proof using laminations, see \cite{MS2, Ot}.

For a holomorphic quadratic differential $\varphi\neq 0$ on $S$ we have seen in Section \ref{S:qd} how to define horizontal and vertical foliations.  If the $\varphi$-coordinate is locally given by $w=u+iv$, then transverse measures may be defined by $|du|$ and $|dv|$, respectively.  In other words, a nonzero quadratic differential defines a measured foliation via its horizontal foliation. \index{measured foliation!and quadratic differentials}  We denote the corresponding dual tree by $T_\varphi$.

The following fundamental theorem, due to Hubbard-Masur and also announced by Thurston, asserts that every measured foliation on $S$ arises in this way:

\begin{theorem}[Hubbard-Masur \cite{HM}] \label{T:hm} \index{Hubbard-Masur theorem} 
Given a measured foliation $({\mathcal F},\nu)$ on a closed Riemann surface $S$ of genus $p\geq 2$ there is a unique holomorphic quadratic differential whose horizontal foliation is equivalent to $({\mathcal F},\nu)$.  In particular, $\MF$ is homeomorphic to $\R^{6p-6}\setminus\{0\}$, and $\PMF\simeq S^{6p-6}$.
\end{theorem}
In Section \ref{S:skora}, we will sketch how we can interpret the Hubbard-Masur theorem via harmonic maps to trees (see \cite{Wf3, Wf4}).  

There is a particular class of quadratic differentials on $S$ called \emph{Jenkins-Strebel differentials} (cf.\ \cite{Str}). \index{Jenkins-Strebel differential}They are characterized by the property that the noncritical trajectories are all closed and they partition the complement of the critical trajectories in $S$ into cylinders with the standard foliations (see Figure 2).  Notice that in this case the dual tree $T_\varphi$ is a simplicial tree with a $\Gamma$ action.  The quotient $T_\varphi/\Gamma$ is a graph $G_\varphi$, and the quotient map $\tilde p:\HH\to T_\varphi$ descends to a map $p:S\to G_\varphi$, as indicated in the Figure 2.

\setlength{\unitlength}{1cm}
\begin{picture}(14,5.5)
\put(5,2.3){$ p$}
\put(1.5,4.1){$S\ =$}
\put(1.5,1.2){$G_\varphi\ =$}
\ifx\pdftexversion\undefined
\put(2.6,1){
{\scalebox{.4}{\includegraphics{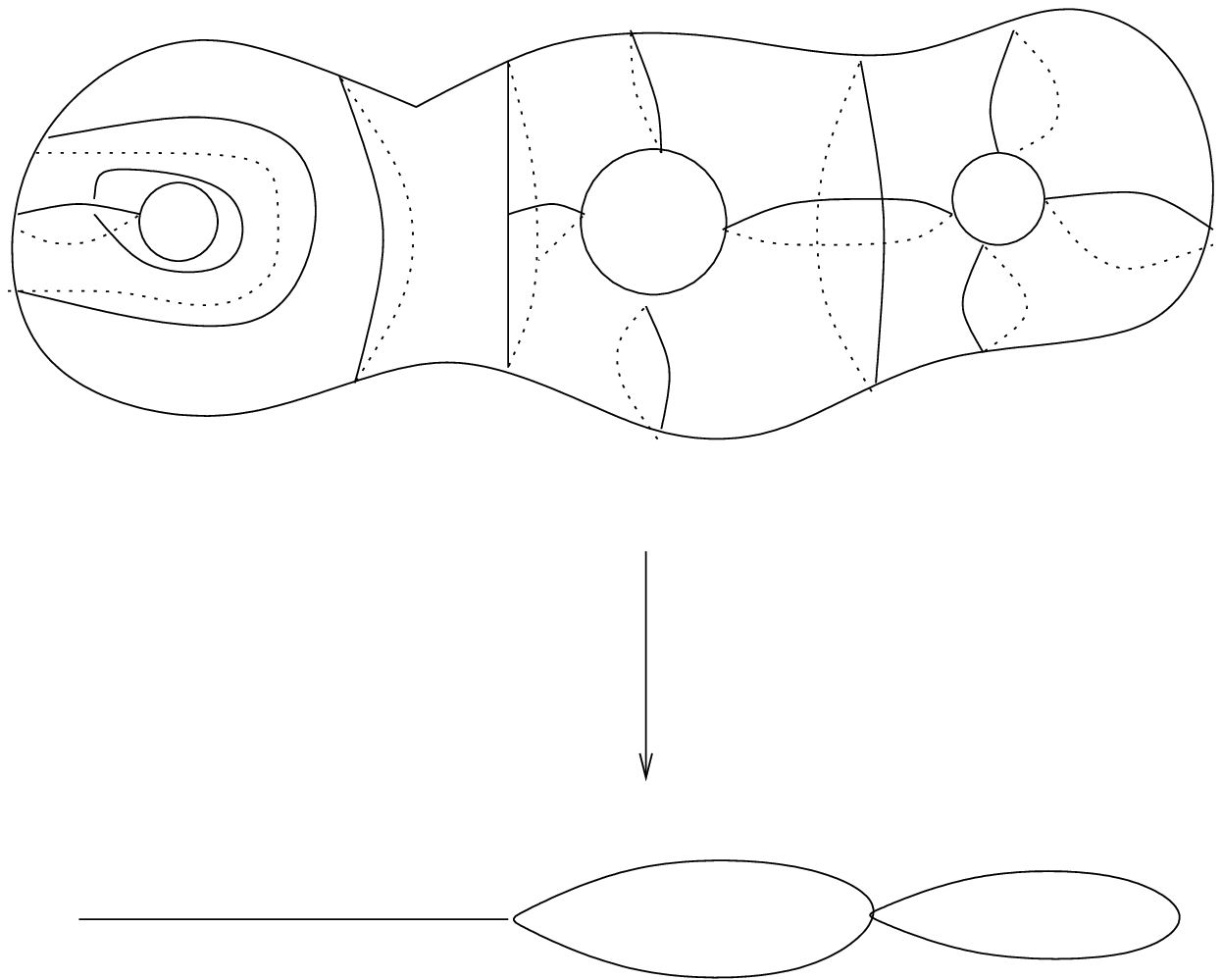}}}}
\else
\put(2.6,1){
{\scalebox{.4}{\includegraphics{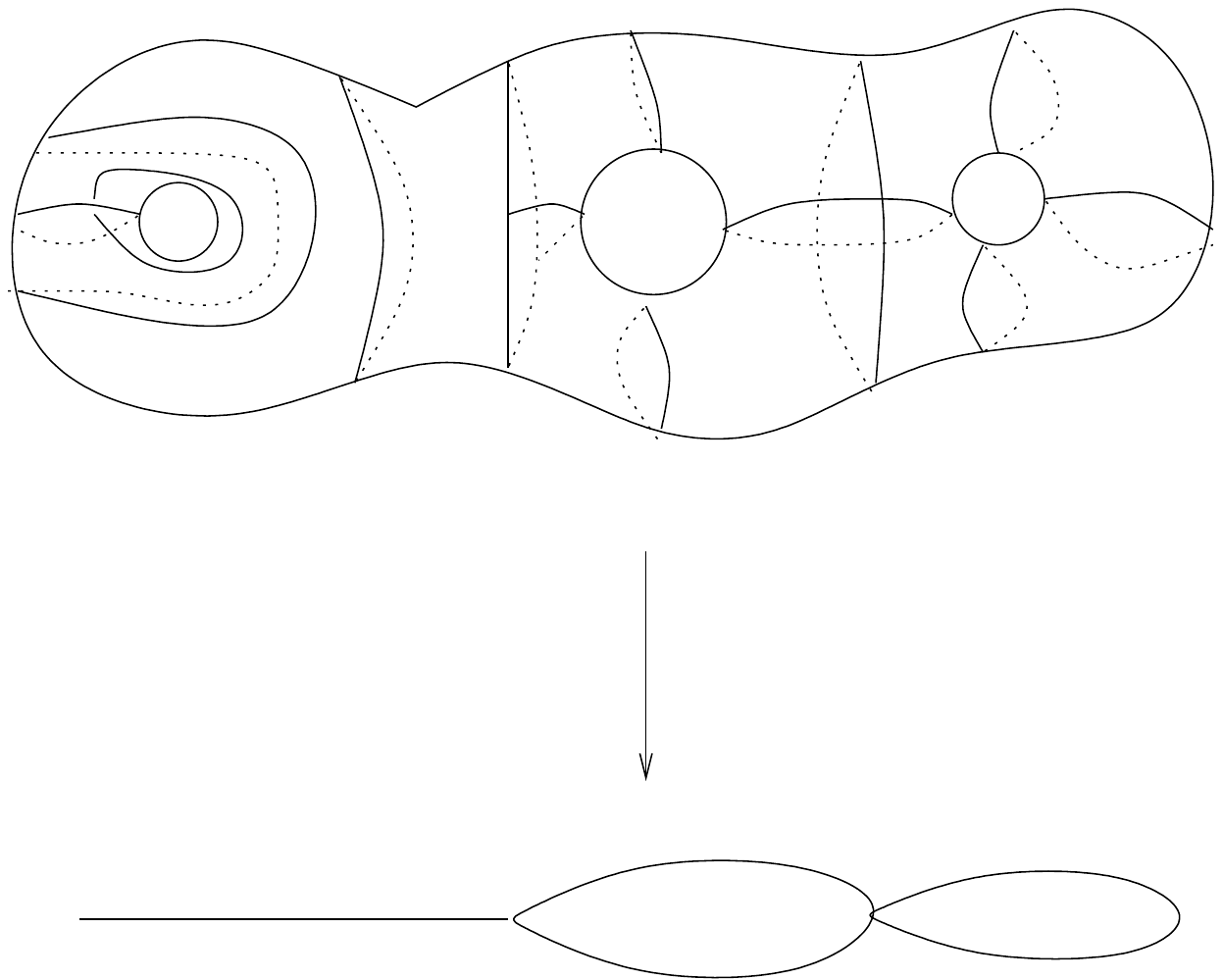}}}}
\fi
\put(4.5,0){Figure 2.}
\end{picture}
\smallskip

Hence, the intersection number \eqref{E:intersection} number is a generalization of the geometric intersection number of simple closed curves.
Let us point out two facts (cf.\ \cite{FLP, Str}).
\begin{itemize}
\item  There are examples of measured foliations where all the noncritical leaves are noncompact.  For example, the fixed points in $\PMF$ of pseudo-Anosov mapping classes.  These give rise to $\R$-trees which are not simplicial.
\item  However, the measured foliations whose associated trees are simplicial are dense in $\MF$.
Furthermore, the intersection number \eqref{E:intersection} extends continuously to $\MF\times \MF$.\index{measured foliation!intersection number}
\end{itemize}

Prior to a rigorous definition of energy minimizers to NPC spaces, we first  introduce the notion of a harmonic map to a tree.  This definition is due to Wolf \cite{Wf2} and  is motivated by Ishihara's Theorem \ref{T:ishihara}.  As we shall see in the next subsection, it turns out that for the case of trees it is equivalent to the definition of energy minimizers due to Korevaar-Schoen.

Let $S$ be a Riemann surface and let $(T,d)$ be a minimal $\R$-tree with an isometric action of $\Gamma=\pi_1(S)$.  Let $f:\HH\to T$ be a $\Gamma$-equivariant, continuous map.  We say that $f$ is \emph{harmonic} if it pulls back germs of convex functions on $T$ to germs of subharmonic functions on $\HH$.  Notice that a function $f:U\to\R$, where $U$ is a convex open subset of an $\R$-tree, is called convex if for any segment $pq\subset U$ and $r\in pq$ we have
$$
(f(r)-f(p))d(q,r)\leq (f(q)-f(r))d(p,r)\ .
$$
A basic example of a harmonic map to a tree is the projection
\begin{equation} \label{E:projection}
p:\HH\lra T_{\varphi}
\end{equation}
where $T_{\varphi}$ is the dual tree to the horizontal foliation of a holomorphic quadratic differential.
  It is not hard to see by direct observation that $p$ is harmonic (cf.\ \cite{Wf2}).

\subsubsection{Harmonic maps to NPC spaces.} \label{S:npc}

For the purpose of this subsection $(\Omega, g)$ will be  a bounded Riemannian domain of dimension $m$ with Lipschitz boundary and $(X,d)$ any complete NPC space.  References for the following are \cite{GS, J5, KS1, KS2}.  The generalization to the case where $X$ is assumed only to have curvature bounded from above can be found in \cite{Me2}.

A Borel measurable map $f:\Omega\to X$ is said to be in $L^2(\Omega,X)$ if for $p\in X$,
$$
\int_\Omega d^2(p, f(x)) dvol_\Omega(x) <\infty\ .
$$
By the triangle inequality, the condition is independent of the choice of point $p$.  For $f\in L^2(\Omega,X)$ we construct an $\varepsilon$-approximate energy function $e_\varepsilon(f) : \Omega_\varepsilon\to \R$, where $\Omega_\varepsilon=\{ x\in \Omega : d(x,\partial\Omega)>\varepsilon\}$ by 
$$
e_\varepsilon(f)(x)=\frac{1}{2\omega_m}\int_{\partial B_\varepsilon(x)} \frac{d^2(f(x),f(y))}{\varepsilon^2}\frac{d\sigma(y)}{\varepsilon^{m-1}}\ ,
$$
where $\omega_m$ is the volume of the unit sphere in $\R^m$ and $d\sigma$ is the induced volume on the sphere $\partial B_\varepsilon(x)\subset\Omega$ of radius $\varepsilon$ about $x$.  Setting $e_\varepsilon(f)(x)=0$ for $x\in \Omega\setminus\Omega_\varepsilon$, we can consider $e_\varepsilon(f)$ to be an $L^1$ function on $\Omega$.  In particular, it defines a linear functional $E_\varepsilon : C_c(\Omega)\to \R$.  We say that $f$ has \emph{finite energy} (or that $f\in H^1(\Omega, X)$) if
$$
E(f)\equiv \sup_{0\leq \varphi\leq 1}\limsup_{\varepsilon\to 0} E_\varepsilon(\varphi) <\infty\ .
$$
It can be shown that if $f$ has finite energy, the measures $e_\varepsilon(f)(x)dvol_\Omega(x)$ converge weakly to a measure that is absolutely continuous with respect to Lebesgue measure on $\Omega$.  Therefore, there is a well-defined integrable function $e(f)(x)$, which we call the \emph{energy density}, \index{energy!density}so that for each $\varphi\in C_c(\Omega)$,
$$\lim_{\varepsilon\to 0}\int_\Omega e_\varepsilon(f)(x)\varphi(x) dvol_\Omega(x)=\int_\Omega e(f)(x)\varphi(x) dvol_\Omega(x)\ .
$$
By analogy with the case of smooth maps we write $e(f)(x)=\frac{1}{2}|\nabla f|^2(x)$ with total energy
$$
E(f)=\frac{1}{2}\int_\Omega |\nabla f|^2 dvol_\Omega\ .
$$
Similarly, the directional energy measures $|f_\ast (Z)|^2 dvol_\Omega$ for $Z\in \Gamma( T\Omega)$ is a Lipschitz tangent vector field can also be defined as a weak-$\ast$ limit of measures ${}^Ze_\varepsilon(f) dvol_\Omega$.  Here,
$$
{}^Ze_\varepsilon(f)(x)=\frac{d^2 (f(x), f(\bar x(x,\varepsilon))}{\varepsilon^{2}}\ ,
$$
where $\bar x(x,\varepsilon)$ denotes the flow along $Z$ at time $\varepsilon$, starting at $x$.  For almost all $x\in \Omega$,
$$
|\nabla f|^2(x)=\frac{1}{\omega_m}\int_{S^{m-1}}|f_\ast(v)|^2 d\sigma(v)\ ,
$$
where $S^{m-1}$ is the unit sphere in $T_x\Omega$.
This definition of the  Sobolev space $H^1(\Omega, X)$ is consistent with the usual definition when $X$ is a Riemannian manifold.  

For any map $f\in H^1(\Omega, X)$ we can also make sense of the notion of the pullback metric
\begin{equation} \label{E:pi}
\pi : \Gamma(T\Omega)\times \Gamma(T\Omega)\lra L^1(\Omega)
\end{equation}
defined by 
$$
\pi(V,W)=\tfrac{1}{4} |f_\ast(V+W)|^2-\tfrac{1}{4}|f_\ast(V-W)|^2\ ,\ V,W\in \Gamma(T\Omega)\ .
$$
If the tangent space to $(\Omega, g)$ has a local frame $(u_1,\ldots, u_m)$, we write $\pi_{\alpha\beta}=\pi(u_\alpha,u_\beta)$, and 
\begin{equation} \label{E:density}
e(f)=\tfrac{1}{2}|\nabla f|^2=\tfrac{1}{2}g^{\alpha\beta}\pi_{\alpha\beta}\ .
\end{equation}
The $L^1$-tensor will be used in the next section to define the Hopf differential.

A finite energy  map $f:\Omega\to X$ is said to be \emph{harmonic} \index{harmonic map}if it is locally energy minimizing. In other words, for each point $x\in \Omega$  and each neighborhood of $x$, all comparison maps agreeing with $f$ outside this neighborhood have total  energy no less than $f$.   The following are the basic existence and regularity results.  For an alternative approach, see \cite{J5}.

\begin{theorem}[Korevaar-Schoen \cite{KS1}, see also \cite{Serb}] \label{T:ks1}
Let $(X,d)$ be an NPC space.  If $f:\Omega\to X$ is harmonic, then $f$ is locally Lipschitz continuous.  The Lipschitz constant on $U\subset\subset \Omega$ is of the form $C(U)\sqrt{E(f)}$, where $C(U)$ is independent of the map $f$ (cf.\ Proposition \ref{P:lipschitz}).
\end{theorem}

\begin{theorem}[Korevaar-Schoen \cite{KS1}] \label{T:kscompact}
Let $(X,d)$ be compact and NPC.  Let $M$ be a compact Riemannian manifold without boundary, and $f:M\to X$ a continuous map.  Then there exists a Lipschitz harmonic map homotopic to $f$.
\end{theorem}

Note that Theorem \ref{T:kscompact} is a generalization of the Eells-Sampson Theorem \ref{T:es}.  The uniqueness result in the singular case is due to Mese.

\begin{theorem}[Mese \cite{Me5}] \label{T:meseunique}
Let $M$ be a compact Riemannian manifold and $X$ a compact metric space with curvature bounded above by a constant $\kappa<0$.  If $f:M\to X$ is a nonconstant harmonic map, then $f$ is unique in its homotopy class unless it maps onto a geodesic.
\end{theorem}

An important tool in understanding the structure of harmonic maps is the \emph{monotonicity formula}\index{harmonic map!monotonicity formula} for energy minimizers.  The idea goes back to Almgren \cite{Alm}.  The statement is that for nonconstant energy minimizers, the quantity
\begin{equation} \label{E:monotonicity}
e^{C\varepsilon}\frac{\varepsilon\int_{B_\varepsilon(x)}|\nabla f|^2 dvol_\Omega}{\int_{\partial B_\varepsilon(x)} d^2(f(x), f(y))ds(y)}\ ,
\end{equation}
is monotone increasing in $\varepsilon$, for some constant $C$.  The extension of this to singular space targets was obtained in \cite{GS, Schoen2}, and further developed in \cite{Me2}.
The basic idea is that since the derivation of the formula depends only on domain variations, and not on any differentiability of the target space, it continues to hold for maps to metric space targets.

The monotonicity \eqref{E:monotonicity} can be used to construct linear approximations to harmonic maps, and in some cases further regularity can be derived.
A key quantity is  the \emph{order function}. \index{harmonic map!order of}Roughly speaking, the order of a harmonic map $f:\Omega\to X$ at a point $x$ measures the degree of the dominant homogeneous harmonic polynomial which approximates $f-f(x)$.  This is precisely true when $X$ is a smooth manifold.  In the general case, it is defined as follows.   Define
\begin{equation} \label{E:order}
\ord_x(f)=\lim_{\varepsilon\downarrow 0}\frac{\varepsilon\int_{B_\varepsilon(x)}|\nabla f|^2 dvol_\Omega}{\int_{\partial B_\varepsilon(x)} d^2(f(x), f(y))ds(y)}\ .
\end{equation}

It follows from the monotonicity formula \eqref{E:monotonicity}  that the above limit exists and is $\geq 1$ for nonconstant maps.  We call this limit the \emph{order of $f$ at $x$}.  It is not an integer in general.  For example, let $p:\HH\to T_\varphi$ be the projection map \eqref{E:projection}.  If $x$ is not a zero of $\varphi$, then $p$ is locally a harmonic function and $\ord_x(p)$ is the order of vanishing.  If $x$ is a zero of order $k$, then $\ord_x(p)=(k+2)/2$.  
The order is related to the eigenvalues of the Laplacian on subdomains of $\partial B_r(x)$, as explained in \cite[Theorem 5.5]{GS}.  In Figure 3, $\ord_x(p)$ is equal to the first Dirichlet eigenvalue of the domain $D_i$ in the circle around $x$.  It is clear in this case that it is equal to $3/2$. 

\setlength{\unitlength}{1cm}
\begin{picture}(14,4.75)
\put(3.66,2.65){$x$}
\put(2.72,3.2){$y$}
\put(2.5, 3.05){$\bullet$}
\put(3.54, 2.43){$\bullet$}
\put(6.1,2.75){$p$}
\put(8.5,2.5){$p(x)$}
\put(8.11, 2.43){$\bullet$}
\put(3.45,4.2){$D_1$}
\put(5,1.4){$D_2$}
\put(1.8,1.4){$D_3$}
\ifx\pdftexversion\undefined
\put(2,1){
{\scalebox{.45}{\includegraphics{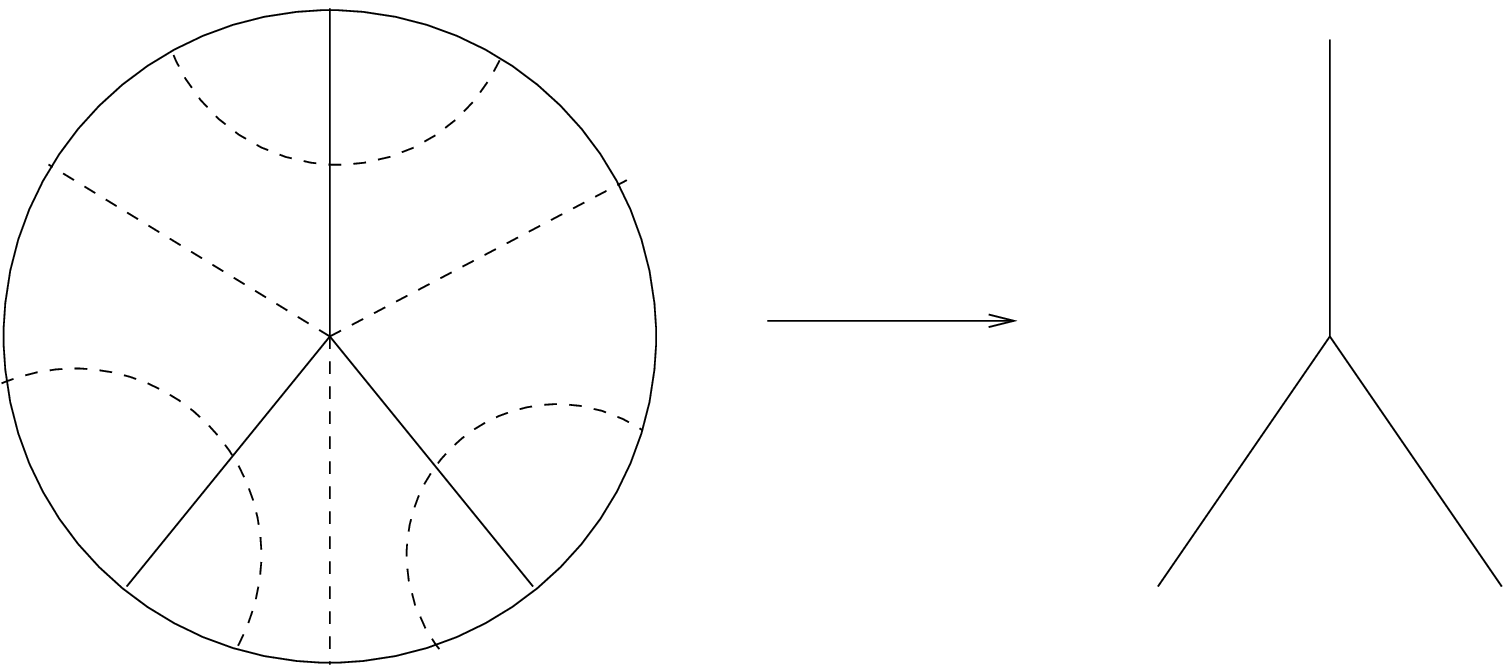}}}}
\else
\put(2,1){
{\scalebox{.45}{\includegraphics{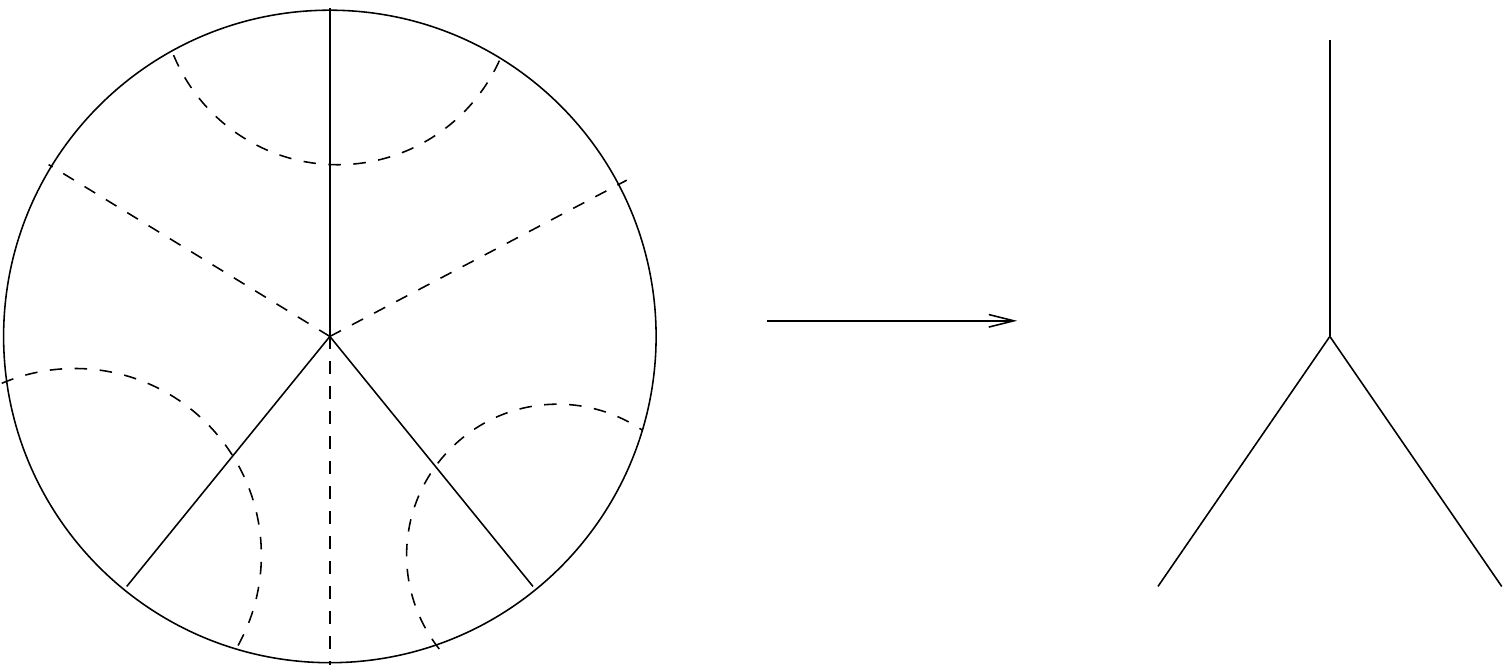}}}}
\fi
\put(4.5,0){Figure 3.}
\end{picture}
\smallskip

On the other hand, if $y\in p^{-1}(p(x))$ is not a zero, then $\ord_y(p)=1$, and indeed locally near $y$, $p$ maps to an interval.

This can be generalized.
Let $f:\Omega\to T$ be a harmonic map to an $\R$-tree.  A point $x\in \Omega$ is called \emph{regular}\index{regular point} if there exists $r>0$ such that $f(B_r(x))$ is an embedded arc.  In particular, $f$ restricted to $B_r(x)$ is then a harmonic function.  Nonregular points are called \emph{singular}.\index{singular point}
In the case of two-dimensional domains, the harmonic map $p: \HH\to T_\varphi$ has singularities precisely at the zeros of $\varphi$.  In particular, they are of codimension $2$.  The next result  was proven in \cite{GS} for simplicial trees and in \cite{Sun} for $\R$-trees.

\begin{theorem} \label{T:regularity}
Let $f:\Omega\to T$ be a harmonic map to an $\R$-tree.  Then $x\in \Omega$ is regular if $\ord_x(f)=1$.  Moreover,
   the Hausdorff codimension of the singular set is at least $2$.
\end{theorem}

X. Sun also proved the following useful fact.
\begin{theorem}[Sun \cite{Sun}] \label{T:sun}
Let $f:\Omega\to T$ be a harmonic map to an $\R$-tree.  Then for any point $x\in \Omega$ there is $r>0$ such that $f(B_r(x))$ lies in a locally finite subtree.
\end{theorem}


\section{Harmonic Maps and Representations}

\begin{itemize}
\item 3.1 Equivariant Harmonic Maps
\item 3.2 Higgs Bundles and Character Varieties
\end{itemize}

\vspace{-.25in}



\subsection{Equivariant Harmonic Maps} \label{S:equivariant}


In this section we describe the equivariant harmonic map problem and its applications.  In Section \ref{S:reductive}, we introduce the notion of reductivity (or semisimplicity) in different contexts and indicate how it is related to the existence problem for equivariant harmonic maps.  In Section \ref{S:skora}, we discuss the holomorphicity of the Hopf differential for harmonic maps and show how it can be used to simplify the proofs of the Hubbard-Masur and Skora theorems.  We also give   the first variation harmonic maps with respect to the domain metric and apply this to derive Gardiner's formula.

\subsubsection{Reductive representations.} \label{S:reductive}
Throughout this section, unless otherwise noted, $(M,g)$ is  a closed Riemannian manifold with $\Gamma=\pi_1(M)$, and  $(X,d)$ is  a simply connected NPC space. Let $\widetilde M$ denote the universal cover of $M$.  We assume that $\Gamma$ acts on $X$ via isometries, i.e.\ that there is a homomorphism $\rho: \Gamma\to \iso(X)$. Associated to $\rho$ is a \emph{translation length function}\index{translation length!of a representation} 
\begin{equation} \label{E:tl}
\tl:\Gamma\lra\R^+: \gamma\mapsto\inf_{x\in X} d(x,\rho(\gamma)x)\ .
\end{equation}
 Let $f:\widetilde M\to X$ be a $\rho$-equivariant map.  Provided that $f$ is a locally in $H^1$, the energy density $|\nabla f|^2$ is $\Gamma$-invariant, and therefore we can define the energy by
\begin{equation}\label{E:totalenergy}
E(f)=\frac{1}{2}\int_{M=\widetilde M/\Gamma}|\nabla f|^2 dvol_M\ .
\end{equation}
Finite energy maps always exist, and indeed energy minimizing sequences  can be taken to be uniformly Lipschitz \cite{KS2}. 
Under conditions that will be made precise below and which we will always assume,  there exist maps with finite  energy.
A $\rho$-equivariant map $f:\widetilde M\to X$ which is locally in $H^1$ is called \emph{harmonic}\index{harmonic map!equivariant} if it minimizes the  energy \eqref{E:totalenergy} among all other equivariant maps in $H^1_{loc.}$.

It follows from the trace theory in \cite{KS1} that equivariant harmonic maps are locally energy minimizers.  Therefore, in the case where $X$ is a smooth manifold the first variational formula \eqref{E:firstvar} implies that a $\rho$-equivariant harmonic map is equivalent to a smooth $\rho$-equivariant map that satisfies the harmonic map equations \eqref{E:harmonic}.  For general NPC targets it follows from Theorem \ref{T:ks1} that $\rho$-equivariant harmonic maps are Lipschitz.

The existence of equivariant harmonic maps is more complicated than in the case of compact targets.  The reason for this is that in the process of choosing an energy minimizing sequence, e.g.\ using the heat flow as in the Eells-Sampson theory, the map can ``escape to infinity," and fail to converge.  An example of this phenomenon can be found in  \cite{ES}.   
One naturally looks for a condition on the homomorphism $\rho$ which rules out this kind of behavior.  For example, it is reasonable to rule out the existence of a sequence of points escaping to infinity whose translates by fixed elements in the image of $\rho$ remain bounded.  This is the notion of a \emph{proper} action (see below).

Before making this more precise we introduce the notion of the \emph{ideal boundary}\index{ideal boundary} of an NPC space.  By a \emph{ray}\index{ray} in $X$ we mean a geodesic $\alpha$ parametrized by arc length on the interval $[0,\infty)$.  Two rays $\alpha_1$, $\alpha_2$,  are said to be equivalent if the Hausdorff distance between them is finite.  Denote by $\partial X$ the set of equivalence classes of rays.  Notice that since $\Gamma$ acts by isometries, $\Gamma$ also acts on $\partial X$.  We have the following facts:
\begin{enumerate} 
\item (cf.\ \cite{BH})  If $X$ is locally compact then $\overline X=X\cup\partial X$ can be topologized so that it becomes a compact metric space.
\item (cf.\ \cite{CM})  If $(X,d)$ is an $\R$-tree (not necessarily locally compact) then two rays $\alpha_1$ and $\alpha_2$ are equivalent if and only if $\alpha_1\cap\alpha_2$ is another ray.  
\item (cf.\ \cite{CM})   If $(X,d)$ is an $\R$-tree with $\Gamma$ action, then $\Gamma$ fixes a point on $\partial X$ if and only if $\tl(\gamma)=|r(\gamma)|$ where $r:\Gamma\to \R$ is a homomorphism.
\end{enumerate}

We now state 
\begin{theorem}[Korevaar-Schoen \cite{KS2}] \label{T:ks2} \index{harmonic map!existence of}
 Suppose $\rho:\Gamma\to \iso(X)$ is a homomorphism  that  does \emph{not} fix a point of $\partial X$.  If either (i) $X$ is locally compact, or (ii) $X$ has curvature bounded above by $\kappa<0$, then there exists a $\rho$-equivariant harmonic map $f:\widetilde M\to X$.
\end{theorem}

The equivariant version of Theorem \ref{T:meseunique} also holds:

\begin{theorem}[Mese \cite{Me5}] \label{T:meseunique2}
If $X$ has curvature bounded above by a constant $\kappa<0$, and  if $f:\widetilde M\to X$ is a nonconstant equivariant harmonic map, then $f$ is unique in its equivariant homotopy class unless it maps onto a geodesic.
\end{theorem}

Special cases of Theorem \ref{T:ks2} had been proven earlier:

\medskip\noindent $\bullet$ \emph{The Corlette-Donaldson Theorem}.
\begin{theorem}[cf.\ \cite{C1, C2, D1}]  \label{T:dc} \index{Corlette-Donaldson theorem}
Let $X$ be  a Riemannian symmetric space of noncompact type $X=G/K$, where $G$ is a semisimple Lie group and $K$ a maximal compact subgroup.  Let $\rho:\Gamma\to G$ be a homomorphism with Zariski dense image.  Then there is a $\rho$-equivariant (smooth) harmonic map $f:\widetilde M\to X$.
\end{theorem}
This theorem is implied by Theorem \ref{T:ks2}, since if $\rho(\Gamma)$ fixes a point $[\alpha]\in\partial X$, then $\overline {\rho(\Gamma)}$ would be closed subgroup contained in the stabilizer of $[\alpha]$, which is a proper subgroup of $G$.  See also \cite{JY1}.

\medskip\noindent$\bullet$  \emph{Labourie's Theorem}.
In the Riemannian case, the criterion for existence in terms of fixing a point in the ideal boundary was conjectured in \cite{C1} and proved in \cite{Labourie} (see also \cite{J4}).   A homomorphism $\rho:\Gamma\to \iso(X)$ is called \emph{semisimple}\index{semisimple!representation} (or \emph{reductive}) if either $\rho(\Gamma)$ does not fix a point in $\partial X$ or it fixes a geodesic.  Then we have the following
\begin{theorem}[Labourie  \cite{Labourie}]  \label{T:labourie} 
Let  $X$ be a Riemannian manifold with  negative sectional curvature.  Then there exists a $\rho$-equivariant harmonic map $f:\widetilde M\to X$ if and only if $\rho$ is semisimple.
\end{theorem}

\medskip\noindent $\bullet$ \emph{$\R$-trees}.   Let $(T,d)$ be an $\R$-tree and $\rho:\Gamma\to \iso(T)$ a homomorphism.  We assume (without loss of generality) that the action of $\Gamma$ on $T$ is minimal.  
\begin{theorem}[Culler-Morgan \cite{CM}]
Let $\rho_1$, $\rho_2$ be nontrivial semisimple actions on $\R$-trees $T_1$, $T_2$ with the same translation length functions.  Then there exists an equivariant isometry $T_1\simeq T_2$.  If either action is not isometric to an action on $\R$, then the equivariant isometry is unique.
\end{theorem}

Then we have the following generalization of Theorem \ref{T:labourie} to trees.

\begin{theorem} \label{T:semisimple}
Let $(T,d)$ be a minimal $\R$-tree and $\rho:\Gamma\to \iso(T)$.  Then there exists a $\rho$-equivariant harmonic map $u:\widetilde M\to T$ if an only if $\rho$ is semisimple.
\end{theorem}

\begin{proof}
The sufficiency follows from Theorem \ref{T:ks2} (see Section \ref{S:trees}).  For the converse, suppose $\Gamma$ fixes a point in $\partial T$. If there is a $\rho$-equivariant harmonic map there would necessarily be a family of distinct such maps (see \cite{DDW1}).  By the uniqueness Theorem \ref{T:meseunique2} and the minimality of $T$, it follows that $T$ in this case is equivariantly isometric to $\R$.
\end{proof}

In the case where $X$ is not locally compact, the condition of not fixing a point at infinity does not seem to be sufficient to guarantee existence.  Korevaar and Schoen developed a slightly stronger condition to cover this case.
  Let  $\rho:\Gamma\to\iso(X)$ be a
homomorphism.  To each set of generators $\mathcal G$ of $\Gamma$ we associate a function on $X$:
$$
D_\rho(x) = \max\left\{ d(x,\rho(\gamma) x) : \gamma\in{\mathcal G}\right\}\ .
$$ 
A homomorphism $\rho:\Gamma\to\iso(X)$  is called \emph{proper}\index{proper action}
 if  for every
$B\geq 0$, the set $\{x\in X : D_\rho(x)\leq B\}$ is bounded.
 Clearly, this condition  is independent of the choice of generating set $\mathcal G$.   For complete manifolds
of nonpositive curvature, the existence of two  hyperbolic isometries in the image of $\rho$ with nonasymptotic axes is sufficient to prove properness.   More generally,  $\rho$ being proper implies that $\rho$ has no fixed end, for if $R$ is a fixed ray then $D_\rho$ is bounded along $R$.

\begin{theorem}[Korevaar-Schoen \cite{KS2}] \label{T:ksproper}
  Suppose $\rho:\Gamma\to \iso(X)$ is proper.   Then there exists a $\rho$-equivariant harmonic map $f:\widetilde M\to X$.
\end{theorem}

\noindent   In case $X$ is locally compact this is implied by Theorem \ref{T:ks2}, but for nonlocally compact spaces it is not.   Yet another sufficient condition is introduced in \cite{KS3}.

To end this section, we connect the definition of harmonicity given in this section with that at the end of Section \ref{S:trees}.

\begin{theorem} \label{T:convex}
Let $S$ be a Riemann surface, $(T,d)$ an $\R$-tree, and $\rho:\Gamma=\pi_1(S)\to\iso(T)$ a reductive action.  A $\rho$-equivariant map $f:\widetilde M\to T$ is harmonic if and only if $f$ pulls back germs of convex functions to germs of subharmonic functions.
\end{theorem}

\begin{proof}
The fact that harmonic maps pull back functions to subharmonic ones is the content of \cite[Prop.\ 3.2]{FW} (see also \cite{KS1}).  For the converse, we argue as follows:  suppose $f:\HH\to T$ is a $\rho$-equivariant map that pulls back germs of convex functions to subharmonic ones.  Let $f':\HH\to T$ be a $\rho$-equivariant harmonic map.  Since both $f$, $f'$ pull back germs of convex functions to subharmonic functions, it follows that the same is true for $f\times f':\HH\to T\times T$.  Hence, $d(f,f')$ is $\Gamma$-equivariant and subharmonic, hence constant.  But because of the $1$-dimensionality of trees it is easy to see that the energy densities of $f$ and $f'$ must be equal, so that $f$ is energy minimizing.
\end{proof}
So far as we know, this result for general NPC targets is open.


\subsubsection{Measured foliations and Hopf differentials.} \label{S:skora}


 Recall from Section \ref{S:npc} that if $X$ is a metric space target and $f:(M, g)\to (X,d)$ is a finite energy map, then one can associate an integrable symmetric $2$-tensor $\pi_{\alpha\beta}$ on $S$ with the property that the energy density $|\nabla f|^2=g^{\alpha\beta}\pi_{\alpha\beta}$.  Hence, while the energy density may not be the square of the norm of a derivative, it is a trace of directional energies.    Let us specialize to the case where where the domain is a Riemann surface, and let $f$ be an energy minimizer.  By varying among finite energy maps  obtained from pulling $f$ back by a local diffeomorphism defined by a vector field $v$, we arrive at
$$
0=\int_M\langle \pi, L_v g-(1/2)\Tr_g(L_v g)\rangle_g dvol_M\ .
$$
Note that the integrand is well-defined since $\pi$ is integrable.  By a particular choice of $v$, and using Weyl's lemma on integrable weakly holomorphic functions, we obtain \cite{Schoen2} that
$$\varphi(z)dz^2=\tfrac{1}{4}(\pi_{11}-\pi_{22}-2i\pi_{12})dz^2$$
is a holomorphic quadratic differential on $S$ (cf.\ Section \ref{S:es}).  We call $\varphi$ the \emph{Hopf differential} of $f$.  Since these computations are local, they apply as well to the case of equivariant harmonic maps.

Thus far we have seen that a measured foliation $\mathcal F$ on a surface $S$ gives rise to an $\R$-tree $T_{\mathcal F}$ with an isometric action of $\Gamma$.  This action has the following properties:
\begin{enumerate}
\item  the action is \emph{minimal}\index{tree!minimal action}\index{minimal action} in the sense that no proper subtree is invariant under $\Gamma$ (cf.\ \cite{Ot} -- strictly speaking, the proof there uses geodesic laminations but it can be easily adapted to the case of measured foliations);
\item the action is \emph{small}\index{tree!small action}\index{small action} in the sense that the edge stabilizer subgroups do not contain free groups on $2$-generators (cf.\ \cite{MO} -- more precisely, the stabilizers are cyclic, since leaves on the quotient surface are either lines or circles).
\end{enumerate}

Shalen conjectured \cite{Sh} that every minimal, small action of a surface group on an $\R$-tree is dual to a measured foliation.  This conjecture, which plays an important role in Thurston's hyperbolization theorem for fibered $3$-manifolds (see \cite{Ot}),
was proved by R. Skora, building upon previous work of  Morgan-Otal \cite{MO}.

\begin{theorem}[Skora \cite{Sk}] \label{T:skora} \index{Skora's theorem}
Let $S$ be a surface of genus at least $2$.  Then if $(T,d)$ is an $\R$-tree with a minimal, small isometric action of $\Gamma=\pi_1(S)$, there is a measured foliation $\mathcal F$ on $S$ such that $(T,d)$ is equivariantly isometric to $T_{\mathcal F}$.
\end{theorem}


   For example, we have seen that if $(T,d)$ is dual to a measured foliation on $S$ then the action is small.  It is also a simple matter to see that  a small action is semisimple.  Indeed, choose $\gamma_i\in \Gamma$, $i=1,\ldots, 4$ such that the commutators $[\gamma_1,\gamma_2]$ and $[\gamma_3,\gamma_4]$ generate a group $G$ containing a free group on $2$-generators.  Then if $\rho$ had a fixed end, then $\rho(\gamma_i)$ would act by translations along a common ray.  In particular, $\rho(G)$ would stabilize this ray, contradicting the assumption of smallness.
Hence, small actions are semisimple, and by Theorem \ref{T:semisimple} there exists a $\rho$-equivariant harmonic map $f:\HH\to T$.  

In general, 
let $\tilde\varphi$ be the Hopf differential of an equivariant harmonic map.  Then $\tilde\varphi$ is the lift of a holomorphic quadratic differential $\varphi$ on $S$.
  Let $T_\varphi$ denote the dual tree to the vertical foliation of $\tilde \varphi$.  It is not hard to see (cf.\ \cite{DDW2,FW}) that there is a $\Gamma$-equivariant map $F:T_\varphi\to T$ such that the following diagram commutes
 \begin{equation} \label{E:factorization}
 \xymatrix{ \HH \ar[dr]_{f} \ar[r]^{p} & \ar[d]^{F} T_\varphi    \\
  & T }
 \end{equation}
where $p:\HH\to T_\varphi$ is the natural projection.  Moreover, this is a \emph{morphism of trees},\index{tree!morphism of} meaning that any segment $xy\in T_\varphi$ decomposes into a finite union of subsegments  along which $p$ is an isometry.
By \cite{MO}, it follows that $F$ is either an equivariant isometry, or $F$ \emph{folds} at some point.
This means there is an identification of two or more segments $z'y'$ and $z'y''$ in $T_\varphi$ to a single segment $zy$ in $T$.  An example of a folding\index{tree!folding of} is shown in Figure 4.   

\setlength{\unitlength}{1cm}
\begin{picture}(14,4)
\put(1.4,2){$x'$}
\put(2.31, 2.27){$\bullet$}
\put(3,2){$x''$}
\put(2.5,1.4){$y'$}
\put(2.5,2.9){$y''$}
\put(5,2){$F$}
\put(7.4,2.27){$x$}
\put(8.6, 2.27){$\bullet$}
\put(10,2.27){$y$}

\ifx\pdftexversion\undefined
\put(1,1){
{\scalebox{.5}{\includegraphics{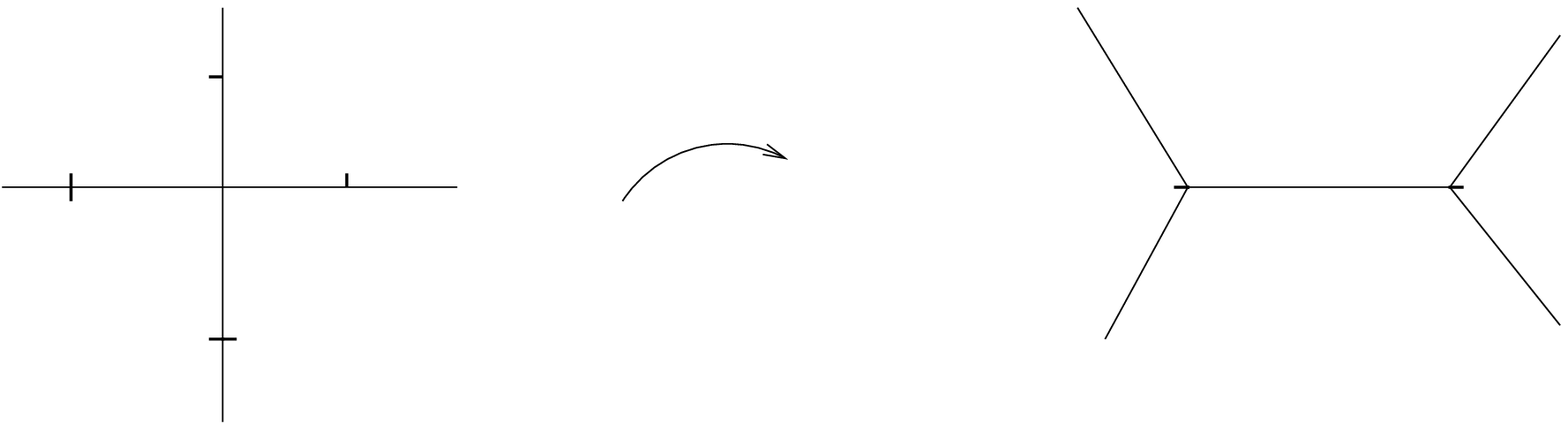}}}}
\else
\put(1,1){
{\scalebox{.5}{\includegraphics{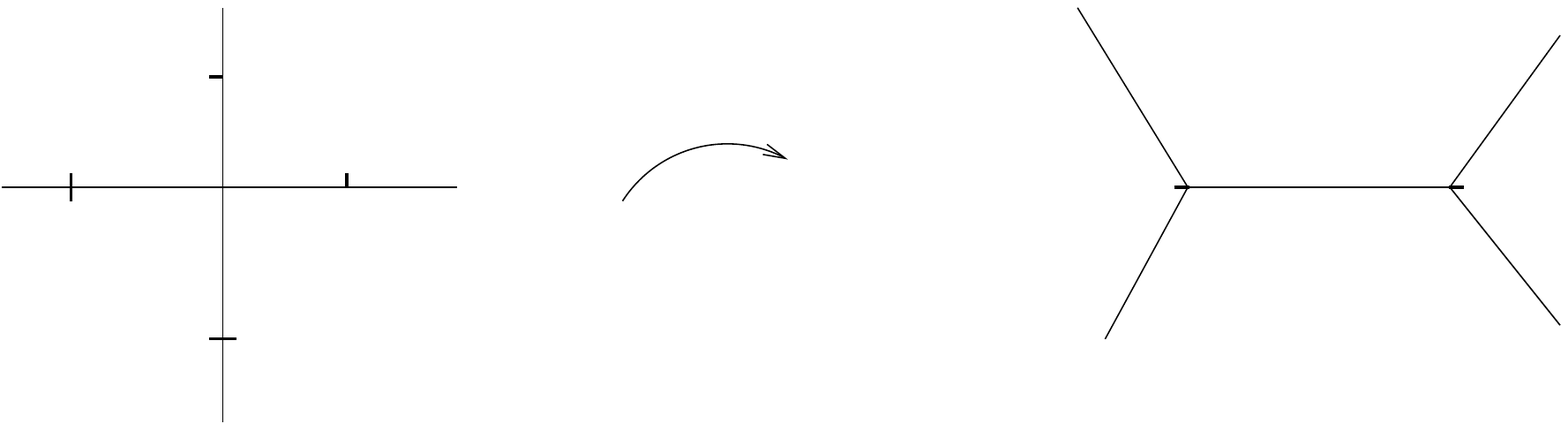}}}}
\fi
\put(4.5,0){Figure 4.}
\end{picture}
\smallskip

Interestingly, the harmonicity of the map $f$ precludes a whole class of undesirable foldings.  For example,  the  following is a consequence of the maximum principle.
\begin{proposition}[\cite{DDW2, FW}]  Suppose  $T_\varphi\to T$  arises from an equivariant harmonic map to $T$, as in \eqref{E:factorization}.  Then folding occurs only at points in $T_\varphi$ corresponding to zeros of $\tilde\varphi$ of multiplicity at least two. Moreover, adjacent  edges may not be identified under such a folding.
\end{proposition}
This type of resolution of the tree $T$ by the dual tree $T_\varphi$ to a measured foliation, with the folding properties of the proposition, had been obtained by Morgan-Otal in \cite{MO}, and it is the first step in proving Skora's Theorem \ref{T:skora}.

By an ingenious counting argument using interval exchanges, Skora went on to show that provided the action of $\Gamma$ is small, folding at vertices  cannot occur either, and in fact $F$ is an isometry.  This completes the proof  of Theorems \ref{T:hm} and \ref{T:skora}.  An alternative source for the counting argument is \cite[\S 8.4]{Ot}.  The reader may also consult \cite{FW, Wf3, Wf4}.

As a second application, consider a  measured foliation $\mathcal F$ on a Riemann surface $(S,\sigma)$.   We have seen above that there is a unique holomorphic quadratic differential $\varphi_{\mathcal F}=\varphi(\sigma,{\mathcal F})$ whose horizontal foliation is measure equivalent to $\mathcal F$.  The \emph{extremal length} of $\mathcal F$ is defined by
\begin{equation} \label{E:height} \index{measured foliation!extremal length}
\ext_{\mathcal F}[\sigma]= \int_S |\varphi_{\mathcal F}| 
\end{equation}
and is a well-defined function on $\Teich$.
It is a generalization of the extremal length of a simple closed curve to the case of arbitrary measured foliations. 
In \cite{G1}, Gardiner gave a formula for the first variation of $\ext_{\mathcal F}[\sigma]$.  Here we show how this formula  arises naturally as the variation of the energy of harmonic maps.

Let $X$ be an NPC space, and suppose $\rho:\Gamma\to \iso(X)$ is proper in the sense of Section \ref{S:reductive}.  For simplicity, assume $X$ has curvature bounded above by some $\kappa<0$. Then for each complex structure $j$ on $S$, Theorems \ref{T:ksproper} and \ref{T:meseunique2} guarantee the existence of a unique $\rho$-equivariant harmonic map $f:\widetilde S\to X$. The energy of $f=f_{[j],\rho}$   gives a well-defined function  depending upon $[j]$ and $\rho$:

\begin{equation}  \label{E:eminus}
{\mathcal E}^-_\rho : \Teich\to \R^+ : [j]\mapsto E(f_{[j],\rho})\ .
\end{equation}

\begin{theorem} \label{T:diff}
The function  ${\mathcal E}^-_\rho$ is differentiable on $\Teich$.  If $\sigma_t$, $-1\leq t\leq 1$, is a differentiable family of metrics on $S$ with Beltrami differential $\mu$ at $t=0$, and $\varphi$ is the Hopf differential of a $\rho$-equivariant energy minimizer $(\widetilde S, \sigma_0)\to X$, then
\begin{equation} \label{E:firstvariation}
\frac{d}{dt}\biggr|_{t=0} {\mathcal E}^-_\rho[\sigma_t]= -4\real\langle\mu,\varphi\rangle\ ,
\end{equation}
where the pairing is as in \eqref{E:pairing}.
\end{theorem}

In the case where $X$ is a smooth Riemannian manifold, this formula has some history.  Wolf \cite{Wf5} provides a derivation and refers to earlier notes of Schoen, as well as \cite{Tromba5, Tromba, J3}.  The  earliest computation of this sort may be due to Douglas (cf.\ \cite[eq.\ (12.29)]{Doug}).
Formally, the proof of \eqref{E:firstvariation} goes as follows.  The total energy is the contraction of the energy density tensor $\pi_{ij}$ with the metric on $S$.  Hence, the first variation involves varying first $\pi_{ij}$, i.e.\ the harmonic map, and then the metric.  But the term associated to the variation of the map is necessarily zero, since the map is energy minimizing.  It follows that the only contribution comes from variations with respect to the metric.  Formula \eqref{E:firstvariation} then follows easily.  Some care must be taken to justify this in the case of metric space targets (see \cite{W2}).

Now consider a measured foliation $\mathcal F$ on $S$ with associated dual tree $T_{\mathcal F}$.  The energy of the unique equivariant harmonic map $f:\widetilde S\to T_{\mathcal F}$ is precisely the extremal length $\ext_{\mathcal F}$.  From \eqref{E:firstvariation} we have

\begin{theorem}[Gardiner \cite{G1,G2}] \label{T:gardiner} \index{Gardiner's formula}
 For any measured foliation $\mathcal F$,  $\ext_{\mathcal F}$ is differentiable on $\Teich$ with derivative
 $$
\frac{d}{dt}\biggr|_{t=0} \ext_{\mathcal F}[\sigma_t]= 2\real\langle\mu,\varphi_{\mathcal F}\rangle\ .
$$
Here,  $\varphi_{\mathcal F}$ is the Hubbard-Masur differential for $\mathcal F$ at $\sigma_0$.
\end{theorem}


\subsection{Higgs Bundles and Character Varieties}


This section discusses the relationship between character varieties and certain special metrics on holomorphic vector bundles.  The link between these two comes via the equivariant harmonic map problem of the previous section.
 In Section \ref{S:hs}, we introduce the notion of a Higgs bundle and discuss the correspondence between stable Higgs bundles, the self-duality equations, and flat $\SL(2,\C)$ connections.  In Section \ref{S:higgsproof}, we give a Higgs bundle interpretation of the Teichm\"uller space and another proof of Theorem \ref{T:teich} using the self-duality equations.  Finally, in Section \ref{S:ms}, we discuss the notion of convergence in the pullback sense and give a harmonic maps interpretation of the Morgan-Shalen-Thurston compactification of character varieties.


\subsubsection{Stability and the Hitchin-Simpson Theorem.} \label{S:hs}


By a \emph{Higgs bundle}\index{Higgs bundle} on a Riemann surface $S$ we mean a pair $(V,\Phi)$, where $V\to S$ is a holomorphic vector bundle and $\Phi$ is a holomorphic section of the associated bundle $\End(V)\otimes K_S$.  Two Higgs bundles $(V,\Phi)$, $(V', \Phi')$ are isomorphic if there exists an isomorphism $\imath : V\to V'$ of holomorphic structures such $\Phi'\circ \imath=\imath\circ\Phi$.  

Recall that  a complex bundle  has a well-defined \emph{degree}\index{degree}, 
$$\deg(V)=\int_S c_1(V)\ ,$$
 where $c_1(V)$ denotes the first Chern class.
The \emph{slope}\index{slope} is defined by
$$
\slope(V)=\deg(V)/\rk(V)\ ,
$$
where $\rk(V)$ is the rank of $V$.  A Higgs bundle $(V,\Phi)$ is called \emph{stable}\index{Higgs bundle!stability}\index{stability} if $\slope(V')<\slope(V)$ for all nontrivial $\Phi$-invariant proper subbundles $V'\subset V$, i.e.\ $V'\neq 0, V$ and $\Phi(V')\subset V'\otimes K_S$.  A Higgs bundle is called \emph{polystable}\index{Higgs bundle!polystable} \index{polystable} if can be written as a direct sum of stable Higgs bundles.

Given a hermitian metric $H$ on a holomorphic  bundle $V$ we will denote the Chern connection by $\nabla _H$, i.e.\ $\nabla _H$ is the unique connection compatible with $H$ and the holomorphic structure (cf.\ \cite{Chern}).  The curvature $F_{\nabla _H}$ takes values in $\ad(V)\otimes \Omega^2(S)$, where $\ad(V)\subset\End(V)$ is the bundle of skew-hermitian endomorphisms.
Let $\omega$ be a K\"ahler form on $S$ normalized so that $\int_S\omega=1$.
 The following result  is due to Hitchin, who first introduced Higgs bundles in this form \cite{Hitchin1, Hitchin3}.  The result for  higher dimensional K\"ahler manifolds is due to Simpson \cite{Si1}.  The case $\Phi\equiv 0$ corresponds to stable bundles on Riemann surfaces and was proved first by Narasimhan-Seshadri \cite{NS} and later, using very different methods, by Donaldson \cite{D}.  Higher dimensional versions of the Narasimhan-Seshadri Theorem were obtained
  by Donaldson \cite{D1} and Uhlenbeck-Yau \cite{UY}.

\begin{theorem}[Hitchin, Simpson]  \label{T:hitchin} \index{Hitchin-Simpson theorem}
Let $ (V,\Phi)$ be a Higgs bundle on a closed Riemann surface $S$.  Then $(V,\Phi)$ is polystable if and only if there exists a hermitian metric $H$ on $V$ solving the self-duality equations\index{self-duality equations}
\begin{equation} \label{E:he}
\frac{i}{2\pi} F_{\nabla _H} + [\Phi, \Phi^{\ast_H}]=s\, {\bf I}\otimes \omega\ ,
\end{equation}
where $ \Phi^{\ast_H}$ is the adjoint of $\Phi$ with respect to $H$, and $s=\slope(V)$.
Furthermore, $H$ is unique up to scalars.
\end{theorem}

From both the algebro-geometric and topological points of view, it is preferable to fix determinants.  In other words, fix a holomorphic line bundle $L\to S$ with hermitian metric $h$ such that $\det(V)=L$, and let $\Phi\in H^0(\End_0 V\otimes K_S)$, where $\End_0(V)$ is the bundle of traceless endomorphisms.  
We shall call $(V,\Phi)$ a Higgs bundle of fixed determinant $L$.

\begin{corollary} \label{C:hitchin}
A Higgs bundle $(V,\Phi)$ of fixed determinant $L$ is polystable if and only if there exists a hermitian metric $H$ on $V$ with $\det H=h$ and such that \eqref{E:he} holds.  In the case such an $H$ exists, it is unique.
\end{corollary}

Following Corlette \cite{C1} we call a flat $\SL(r,\C)$ connection $\nabla  $ on the trivial rank $r$ bundle on $S$ \emph{reductive}\index{reductive} if any $\nabla  $-invariant subbundle has a $\nabla  $-invariant complement.  Clearly, a reductive flat $\SL(r,\C)$ connection is a direct sum of irreducible flat $\SL(r', \C)$ connections for values $r'<r$.

Define ${\mathfrak M}(S,r)$ to be the moduli space of isomorphism classes of polystable Higgs bundles on $S$ of rank $r$ and fixed trivial determinant.\index{Higgs bundle!moduli space} We denote  the space of equivalence classes of reductive flat $\SL(r,\C)$ connections on the trivial rank $r$ bundle $V\to S$ by $\chi(\Gamma, r)$.  We have the following

\begin{theorem}[Corlette, Donaldson]  \label{T:dc2} \index{Corlette-Donaldson theorem}
The map 
$$
\Psi :{\mathfrak M}(S,r)\to \chi(\Gamma, r)\ :\
(V,\Phi)\mapsto \nabla_H+\Phi+\Phi^{\ast_H}\ ,
$$
is a bijection, where $H$ satisfies \eqref{E:he}.  
\end{theorem}
That $\Psi$ is well-defined follows from Corollary \ref{C:hitchin}, and the injectivity is a consequence of the uniqueness of the solution $H$.  The surjectivity part was first conjectured by Hitchin in \cite{Hitchin1} and was subsequently proven for rank $2$ by Donaldson \cite{D2} and in general by Corlette \cite{C1}.
It is equivalent to the Corlette-Donaldson  Theorem  \ref{T:dc}  on equivariant harmonic maps discussed above.

Indeed, given a reductive flat connection $\nabla  $, let $\rho:\Gamma\to\SL(r,\C)$ denote its holonomy representation.  Since reductive representations split into irreducible factors, we may assume without loss of generality that $\rho$ is irreducible.  By Theorem \ref{T:dc} there exists an equivariant harmonic map $H:\HH\to \SL(r,\C)/\SU(r)$.  Equivalently, we can view $H$ as a section of the ``twisted bundle"
$$
\HH\times_\rho \SL(r,\C)/\SU(r)\lra S\ ,
$$
i.e.\ $H$ is  nothing but the choice of a hermitian metric on $V$.  Therefore, we can split $\nabla  =\nabla_H+\Phi+\Phi^{\ast_H}$, where $\nabla_H$ is a hermitian connection with respect to $H$, and $\Phi$ is a smooth section of $\End_0(V)\otimes K_S$.  Clearly, the flatness of $\nabla  $ is equivalent to the equations \eqref{E:he} together with the Bianchi identity
\begin{equation} \label{E:bianchi}
d_{\nabla_H}(\Phi+\Phi^{\ast_H})=0\ .
\end{equation}
 The \emph{harmonicity} of $H$ is equivalent to the condition \cite{C1, D2}
\begin{equation} \label{E:harmonicity}
(d_{\nabla_H})^{\ast_H}(\Phi+\Phi^{\ast_H})=0\ .
\end{equation}
Conditions \eqref{E:bianchi} and \eqref{E:harmonicity} are together equivalent to the holomorphicity of $\Phi$.  Hence, $V$ with the induced holomorphic structure from $\nabla_H$ and $\Phi$ define a Higgs bundle with $\Psi(V,\Phi)=\nabla  $.

Notice that in the above argument we indicated that the existence of a $\rho$-equivariant harmonic map $\HH\to \SL(r,\C)/\SU(r)$  was equivalent to the reductivity of the flat connection $\nabla  $.  Therefore, by Labourie's Theorem \ref{T:labourie}, it is also equivalent to the reductivity of the holonomy representation $\rho:\Gamma\to \SL(2,\C)$,  in the case of $\HH^3=\SL(2,\C)/\SU(2)$.   

The question of the complex structure on the spaces ${\mathfrak M}(S,r)\simeq\chi(\Gamma, r)$, originally addressed by Hitchin \cite{Hitchin1}, is an extremely interesting one.  As a \emph{character variety},\index{character variety}  $\chi(\Gamma, r)$ is an affine algebraic variety.  For example, given any $\gamma\in \Gamma$ we define a regular function $\tau_\gamma:\chi(\Gamma,2)\to \C$ by $\tau_\gamma[\rho]=\Tr\rho(\gamma)$.  Here, $[\rho]$ denotes the conjugacy class of representations containing $\rho$.
By \cite{CS} the ring generated by all elements $\tau_\gamma$, $\gamma\in \Gamma$, is finitely generated.  Fix a generating set associated to $\{\gamma_1,\ldots,\gamma_m\}$, and define
$$
t:\chi(\Gamma,2)\lra \C^m\ :\ [\rho]\mapsto (\tau_{\gamma_1}(\rho),\ldots,\tau_{\gamma_m}(\rho))\ .
$$
Then $t$ is a bijection onto its image and gives $\chi(\Gamma,2)$ the structure of an affine variety.
For higher rank, one needs to consider other invariant polynomials in addition to  traces.
On the other hand, Nitsure and Simpson have shown that ${\mathfrak M}(S,r)$ with its complex structure induced as a moduli space over the Riemann surface $S$ has the structure of a quasiprojective algebraic scheme \cite{Nitsure, Si2, Si3}.  The bijection ${\mathfrak M}(S,r)\simeq\chi(\Gamma,r)$ is \emph{not} complex analytic. On the contrary, Hitchin shows that the two complex structures are part of a hyperk\"ahler family.  For more details, we refer to \cite{Hitchin1, Si2}.

A consequence of the realization of $\chi(\Gamma, r)$ as a moduli space of Higgs bundles is that there is a natural $\C^\ast$-action.  Indeed, if $(V,\Phi)$ is a polystable Higgs bundle then so is $(V,t\Phi)$, $t\in \C^\ast$.  This defines a holomorphic action on ${\mathfrak M}(S,r)$, and therefore also an action (not holomorphic) on $\chi(\Gamma, r)$.  This action depends on the complex structure on $S$ and is not apparent from the point of view of representations.  Nevertheless, we shall see in the next section that it has some connection with Teichm\"uller theory.

\subsubsection{Higgs bundle proof of Teichm\"uller's theorem.} \label{S:higgsproof}


For the purposes of this section we specialize to the case $r=2$ and set ${\mathfrak M}(S)={\mathfrak M}(S,2)$ and $\chi(\Gamma)=\chi(\Gamma,2)$.  Define the \emph{Hitchin map}\index{Hitchin map}
\begin{equation} \label{E:hitchinmap}
\det : {\mathfrak M}(S)\lra \QD(S)\ :\ [V,\Phi]\mapsto\det\Phi=-\frac{1}{2}\Tr\Phi^2\ .
\end{equation}
Hitchin proved that $\det$ is a proper, surjective map with generic fibers being half-dimensional tori.  This last property in fact realizes ${\mathfrak M}(S)$ as a completely integrable system (see \cite{Hitchin2}).  More importantly for us, notice that under the Corlette-Donaldson correspondence $\Psi:{\mathfrak M}(S)\to \chi(\Gamma)$, $\det\circ\Psi^{-1}$ is just the Hopf differential of the associated harmonic map (cf.\ \cite{DDW2}).  Indeed, for $[\rho]\in \chi(\Gamma)$ with an associated equivariant harmonic map $f_\rho$, $\varphi_\rho=\Hopf(f_\rho)$ is given by
$$
\varphi_\rho=\langle \nabla f_\rho^{1,0}, \nabla f_\rho^{1,0}\rangle=-\Tr(\nabla f_\rho^{1,0})^2=2\det\circ\Psi^{-1}[\rho]\ .
$$
In order to realize the Teichm\"uller space inside ${\mathfrak M}(S)$, let $\imath : {\mathfrak M}(S)\to{\mathfrak M}(S)$ denote the involution $\imath(V,\Phi)=(V,-\Phi)$.  Notice that $\imath$ is a restriction of the full $\C^\ast$-action on ${\mathfrak M}(S)$ described at the end of the previous section.  Also notice that under the Corlette-Donaldson correspondence $\Psi$, $\imath$ corresponds to complex conjugation.  Hence, the fixed points of $\imath$ are either $\SU(2)$ or $\SL(2,\R)$ representations.  The former correspond under the  Narasimhan-Seshadri Theorem to the Higgs pair $(V,0)$, i.e.\ $\Phi\equiv 0$.  If $(V,\Phi)$ is a fixed point of $\imath$ with $\Phi\not\equiv 0$, Hitchin shows that $V$ must be a split holomorphic bundle $L\oplus L^\ast$, and with respect to this splitting $\Phi$ is of the form
$$
\Phi=\left(\begin{matrix} 0 & a\\ b& 0\end{matrix}\right)\ ,
$$
where $a\in H^0(S, L^2\otimes K_S)$, and $b\in H^0(S, L^{-2}\otimes K_S)$.  Stability implies $b\neq 0$, and hence by  vanishing of cohomology, $\deg L\leq p-1$, where $p$ is the genus of $S$.  This fact, as pointed out in \cite{Hitchin1}, turns out to be equivalent to the Milnor-Wood inequality which states that the Euler class of any $\PSL(2,\R)$ bundle on $S$ is $\leq 2p-2$ (cf.\ \cite{Milnor, Wood}).\index{Higgs bundle!and the Milnor-Wood inequality} 

We next restrict ourselves to the components of the fixed point set of $\imath$ corresponding to line bundles $L$ of maximal degree $p-1$.  In this case, $L$ must be  a spin structure, i.e.\ $L^2=K_S$, for otherwise $b=0$, contradicting stability.  We denote this moduli space by ${\mathcal N}_L(S)$.  After normalizing by  automorphisms of $L\oplus L^{-1}$, we can write
$$
\Phi=\left(\begin{matrix} 0 & a\\ 1& 0\end{matrix}\right)\ ,
$$
for some quadratic differential $a\in \QD(S)$.  It follows that the restriction of the Hitchin map  to ${\mathcal N}_L(S)$ defines a homeomorphism $\det:{\mathcal N}_L(S)\simrightarrow \QD(S)$.  

The following gives another proof of Theorem \ref{T:teich}. \index{Teichm\"uller!theorem}

\begin{theorem}[Hitchin] \label{T:hitchinteich} \index{Higgs bundle!and Teichm\"uller's theorem}
Given a Higgs bundle 
$$\left(L\oplus L^{-1}, \left(\begin{matrix} 0 & a\\ 1& 0\end{matrix}\right)\right)$$
 in ${\mathcal N}_L(S)$, let $H$ denote the metric on $L\oplus L^{-1}$ solving the self-duality equations, and let $h$ be the induced metric on $K_S^{-1}=L^{-2}=T^{1,0}S$.  Then 
\begin{enumerate}
\item  the tensor 
$$
\hat h= a+(h+h^{-1}a\bar a)+\bar a\in \Omega^0(S, \Sym^2(T^\ast S)\otimes\C)
$$
is a Riemannian metric on $S$ of constant curvature $-4$.
\item  any metric of constant curvature $-4$ on $S$ is isometric to one of this form for some $a\in \QD(S)$.
\end{enumerate}
\end{theorem}
The new ingredient in this theorem is the use of the existence of solutions to the self-duality equations \eqref{E:he}.  Notice that in the reducible case described in Theorem \ref{T:hitchin} the self-duality equations reduce to the \emph{abelian vortex equations}\index{vortex equations}
$F_h=-2(1-\Vert a\Vert_{L^2}^2)\omega$ (cf.\ \cite{JT}).  The relation between the vortex equations and  curvature of metrics on surfaces had been noted previously in the  work of Kazdan and Warner \cite{KW}.

Notice that the definition of ${\mathcal N}_L(S)$ depends on a choice of spin structure $L$, and there are
$\# H^1(S,\Z_2)=2^{2p}$ such choices.   This reflects the fact that on $\chi(\Gamma)$ there is an action of $\Z_2^{2p}$, and the quotient is
$$
\chi(\Gamma)/\Z_{2}^{2p}=\Hom(\Gamma, \PSL(2,\C))// \PSL(2,\C)\ ,
$$
the character variety of $\PSL(2,\C)$, of which the Fricke space $\Fricke$ is a natural subset.  The preimage of $\Fricke$ in $\chi(\Gamma)$ is the disjoint union  of the ${\mathcal N}_L(S)$, and each of these is homeomorphic to Teichm\"uller space.



\subsubsection{The Thurston-Morgan-Shalen compactification.}  \label{S:ms}

Let us first explain the notion of convergence in the pullback sense,\index{convergence in the pullback sense}  due to Korevaar-Schoen, that appears in the statement of Theorem \ref{T:mese}.  
Let $\Omega$ be a set and $f:\Omega\to X$ a map into a simply connected NPC space $(X, d)$.  Use $f$ to define a pseudometric on $\Omega$, $d_f(x,y)=d(f(x),f(y))$, $x,y\in \Omega$.  To obtain convergence in an NPC setting, some convexity is needed.  This is achieved by enlarging $\Omega$ to a space $\Omega_\infty$, defined recursively by: 
\begin{align*}
\Omega_0&=\Omega\ ,\\
\Omega_{k+1}&=\Omega_k\times \Omega_k\times [0,1]\ ,\\
\Omega_\infty&=\bigsqcup_{k=0}^\infty\Omega_k\bigr/\sim\ ,
\end{align*}
where the identification $\sim$ is generated by an inclusion $\Omega_k\hookrightarrow\Omega_{k+1}$, $x\mapsto (x,x,0)$.  The map $f$ extends to $\Omega_\infty$ recursively by setting $f(x,y,t)$, where $x,y\in\Omega_{k+1}$, equal to the point on the geodesic $t$ of the way from $f(x)$ to $f(y)$.
Let $d_\infty$ denote the pullback pseudometric on $\Omega_\infty$.  After identifying points of zero pseudodistance in $(\Omega_\infty, d_\infty)$ and completing, one obtains a metric space $(Z, d_Z)$ isometric to the closed convex hull $C(f(\Omega))\subset X$ (see \cite{KS2}).  

Given a sequence $f_i:\Omega\to X_i$ of maps into simply connected NPC spaces $X_i$, we say that $f_i\to f$ \emph{in the pullback sense}\index{convergence in the pullback sense} if the pullback pseudodistances $d_{i,\infty}$ on $\Omega_\infty$ converge locally uniformly to a pseudometric $d_\infty$, and if the map $f$ is the quotient $\Omega\hookrightarrow\Omega_\infty\to (Z,d_Z)$.   

This notion is equivalent to Gromov-Hausdorff convergence  (cf.\ \cite{Wf2}).\index{Gromov-Hausdorff convergence}  Indeed, (uniform) Gromov-Hausdorff convergence $(Z_i, d_i)\to (Z, d_Z)$ means that for any $\varepsilon>0$ there are   relations $R_i\subset Z_i\times Z$ whose projections surject onto $Z_i$ and $Z$, and such that  if $(z_i,z), (z_i^\prime,z')\in R_i$, then 
$$
| d_i(z_i, z_i^\prime)-d_Z(z,z')|<\varepsilon\ .
$$
Convergence of the maps $f_i:\Omega\to Z_i$ to $f:\Omega\to Z$ imposes the additional requirement that $(f_i(x), f(x))\in R_i$ for all $x\in \Omega$.    It is easy to see that $f_i:\Omega\to X_i$ converges in the pullback sense if and only if the convex hulls  $Z_i=C(f_i(\Omega))$ and the maps $f_i$ converge in the the Gromov-Hausdorff sense.  Indeed, pulling everything back to $\Omega_\infty$, the relations $R_i$ can be taken to be  the diagonal.  We also point out that it is easy to extend these notions equivariantly in the presence of isometric group actions.

We have the following compactness property:
\begin{proposition} \label{P:ks}
Let $\Omega$ be a metric space, and let $f_i:\Omega\to X_i$ be a sequence of maps into NPC spaces such there is a uniform modulus of continuity: i.e. for each $x\in\Omega$ there is a monotone function $\omega_x$ so that $\lim_{R\to 0}\omega_x(R)=0$, and $\max_{y\in B(x,R)} d_{f_i}(x,y)\leq\omega_x(R)$.  Then $f_i$ converges (after passing to a subsequence) in the pullback sense to a map $f:\Omega\to Z$, where $Z$ is an NPC space.
\end{proposition}
We call the NPC space $(Z,d_Z)$ a \emph{Korevaar-Schoen limit}.\index{Korevaar-Schoen limits} Strictly speaking,  the target surfaces $(R,h_i)$ in Theorem \ref{T:mese} are not simply connected and are not NPC. To deal with the former,  consider equivariant convergence of the lifts to the universal covers as mentioned above.  For the latter, one shows that under the assumption that curvature is bounded from above, geodesics are locally unique, so the the construction of $\Omega_\infty$ above works at a local level.

As usual, we denote by  $\Gamma$ the fundamental group of a hyperbolic surface.  Let $\chi(\Gamma)$ be the $\SL(2,\C)$-character variety of $\Gamma$.  As we have seen, $\chi(\Gamma)$ is a noncompact algebraic variety.  In this section we describe a construction, introduced by Thurston in the case of $\SL(2,\R)$ representations, to compactify $\chi(\Gamma)$.  It is important to note that this is not a compactification in an  algebro-geometric sense, and indeed $\overline {\chi(\Gamma)}$ will not be an complex analytic space.

Let ${\mathcal C}$ denote the set of conjugacy classes of $\Gamma$, and let
$$
\PP({\mathcal C})=\left\{ [0,\infty)^{\mathcal C}\setminus\{0\}\right\}\bigr/ \R^+\ ,
$$
where $\R^+$ acts by homotheties. Topologize $\PP({\mathcal C})$ with the product topology. We define a map
$$
\vartheta : \chi(\Gamma) \lra \PP({\mathcal C}) \ :\
[\rho]\mapsto \left\{ \log(|\Tr\rho(\gamma)| + 2)\right\}_{\gamma\in {\mathcal C}}
$$
The purpose of the ``$+2$" in the formula is to truncate the logarithm so that it goes to infinity only when the trace goes to infinity.  It is easy to see (cf.\ \cite{Cooper}) that $ \log(|\Tr\rho(\gamma)| + 2)$ is asymptotic to $\tl(\gamma)$, where $\tl$ is
the translation length function  of $\rho$ acting on hyperbolic space $\HH^3$ (see \eqref{E:tl} and recall that $\iso(\HH^3)=\PSL(2,\C)$).
In case $\rho$ is a discrete faithful $\SL(2,\R)$ representation, hence defining an element of Teichm\"uller space, $\tl(\gamma)$ is just the length of the closed geodesic in the hyperbolic surface $S=\HH/\rho(\Gamma)$ in the free homotopy class of $\gamma$.  

Next, recall from Section \ref{S:hs} that by definition of the affine variety structure on $\chi(\Gamma)$, coordinate functions are of the form $\tau_\gamma$, where $\tau_\gamma(\rho)=\Tr\rho(\gamma)$.  Hence, $\tau_\gamma$, $\gamma\in {\mathcal C}$ generate the coordinate ring of $\chi(\Gamma)$ as a $\C$-algebra, and it follows that $\vartheta$ is a continuous injection.
  Define $\overline{\chi(\Gamma)}$ to be the closure of the image of $\vartheta$ as a subset of $\PP({\mathcal C})$.  It follows, essentially from the finite generation of the coordinate ring of $\chi(\Gamma)$, that $\overline{\chi(\Gamma)}$ is compact (cf.\ \cite{MS1}).  We call $\overline{\chi(\Gamma)}$ the \emph{Morgan-Shalen compactification of $\chi(\Gamma)$},\index{Morgan-Shalen compactification} and set $\partial \chi(\Gamma)=\overline{\chi(\Gamma)}\setminus \chi(\Gamma)$ to be the set of ideal points.  The really useful ingredient in this construction is that the ideal points are not arbitrary but are translation length functions for isometric actions of $\Gamma$ on $\R$-trees.
Another important property  is the following:  the group $\Aut(\Gamma)$  of automorphisms of $\Gamma$ clearly acts continuously on $\chi(\Gamma)$   and this action  admits a continuous extension to $\overline{\chi(\Gamma)}$ (this is essentially the action of the \emph{mapping class group} to be discussed in Section \ref{S:mcg} below).

Let $\Teich\subset \chi(\Gamma)$ denote the Teichm\"uller space, viewed as a component of the discrete faithful representations.  The closure $\thurston$ of $\Teich$ in $\overline{\chi(\Gamma)}$  is called the \emph{Thurston compactification of $\Teich$}\index{Thurston compactification} and $\bthurston=\thurston\setminus \Teich$  is called the \emph{Thurston boundary of $\Teich$}. \index{Thurston boundary}The action of  $\Aut(\Gamma)$ extends continuously to $\Teich\subset\thurston$, and indeed this was part of the motivation for Thurston's compactification.

In terms of a finite set $\gamma_1,\ldots, \gamma_m\in\Gamma$, where $\{ \tau_{\gamma_i}\}$ generate the coordinate ring of $\chi(\Gamma)$,  we can rephrase the compactness of $\overline{\chi(\Gamma)}$  as follows.  Given a sequence of representations $\rho_i:\Gamma\to\SL(2,\C)$, only one of the following can occur:
\begin{enumerate}
\item For some subsequence $\{i'\}$, all traces $\rho_{i'}(\gamma_j)$, $j=1,\ldots, m$, are bounded (in this case, we call the sequence $\rho_{i'}$ \emph{bounded}).  Then $[\rho_{i'}]$ converges (after possibly passing to a further subsequence) in $\chi(\Gamma)$.  
\item For any subsequence $\{i'\}$  there is some $s=1,\ldots, m$ such that $\Tr\rho_{i'}(\gamma_s)\to \infty$ as $i'\to \infty$.  Then there is a function $\ell:{\mathcal C}\to \R^+$, $\ell\neq 0$, such that (after possibly passing to a further subsequence) $\ell_{\rho_{i'}}\to \ell$, projectively.
\end{enumerate}

In terms of the relationship between representations and equivariant harmonic maps we have the following simple but important observation (cf.\ \cite{DDW1}):
\begin{proposition}
A sequence of representations $\rho_i: \Gamma\to \SL(2,\C)$ with associated $\rho_i$-equivariant harmonic maps $f_i:\HH\to \HH^3$ is bounded (up to conjugation) if and only if the energy of the harmonic maps $f_i$ is uniformly bounded.
\end{proposition}

We now assume that $\rho_i$ is an unbounded sequence of representations with $f_i$ as above. Consider the sequence of $\rho_i$-equivariant harmonic maps 
\begin{equation} \label{E:map}
\hat f_i :\HH\lra (\HH^3, d_i)\ ,
\end{equation}
where the hyperbolic metric $d$ on  $\HH^3$ is scaled by the square-root of the energy: $d_i(x,y)=d(x,y)/E^{1/2}(f_i)$, and $\hat f_i=f_i$.  Then because of the scaling the $\hat f_i$ have uniform modulus of continuity.  Furthermore, by  properties  of thin triangles in $\HH^3$ and the fact that $E(f_i)\to \infty$, one can see that geodesic triangles in the convex hull of the image of $\hat f_i$ become infinitely thin (cf.\ \cite{Bestvina, Paulin}).  
\setlength{\unitlength}{1cm}
\begin{picture}(14,3.5)
\ifx\pdftexversion\undefined
\put(1.75,1){
{\scalebox{.45}{\includegraphics{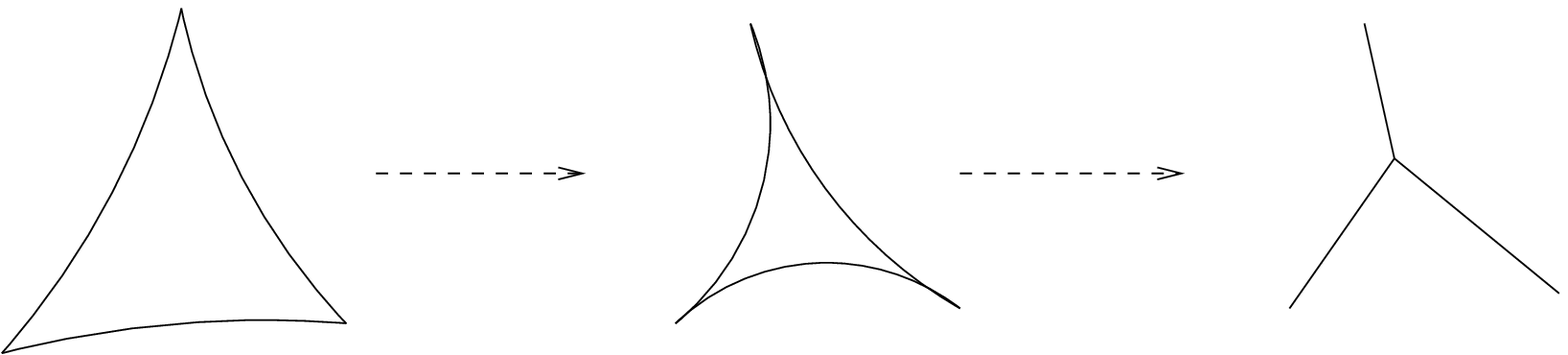}}}}
\else
\put(1.75,1){
{\scalebox{.45}{\includegraphics{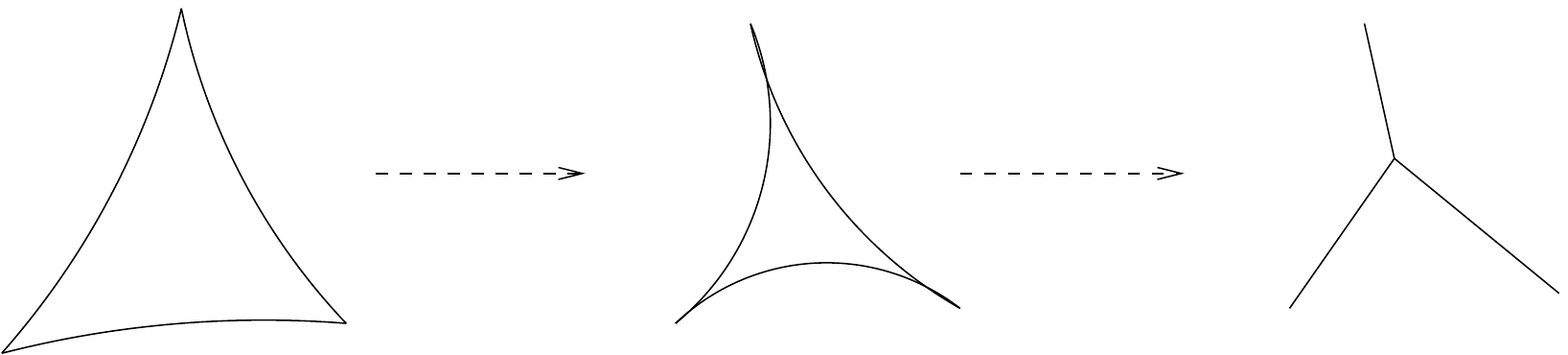}}}}
\fi
\put(4.5,0){Figure 5.}
\end{picture}
\smallskip

\noindent Using these ideas we have
\begin{theorem}[Daskalopoulos-Dostoglou-Wentworth \cite{DDW1}] \label{T:ddw}
For an unbounded sequence of irreducible $\SL(2,\C)$ representations $\rho_i$ the corresponding harmonic maps $\hat f_i$ in \eqref{E:map} converge (after possibly passing to a subsequence) in the pullback sense to a $\Gamma$-equivariant harmonic map $\hat f:\HH\to X$, where $X$ is an $\R$-tree with isometric $\Gamma$ action such that
\begin{enumerate}
\item $\Gamma$ acts on $X$ without fixed points;
\item the length function of the action of $\Gamma$ on $X$ is in the projective class of the Morgan-Shalen limit of the sequence $\rho_i$;
\item the image of $\hat f$ is a minimal tree.
\end{enumerate}
\end{theorem}

Let $\chi_{df}(\Gamma)\subset\chi(\Gamma)$ denote the subspace of discrete faithful representations.  It is a consequence of Jorgenson type inequalities   (cf.\ \cite{MS1}) that the Morgan-Shalen limit of a sequence of discrete faithful representations is the length function of a \emph{small} action on an $\R$-tree.  By Skora's Theorem \ref{T:skora}, the tree is dual to a measured foliation, and therefore $\partial\chi_{df}(\Gamma)\subset\PMF$.  We actually have
\begin{align} 
\partial\chi_{df}(\Gamma)&=\PMF\simeq\QD(S)\ , \label{E:df} \\
 \bthurston &=\PMF\simeq \QD(S)\ .\label{E:compact}
\end{align}
The second equality \eqref{E:compact}, first proven by Thurston using the density of Jenkins-Strebel differentials,  was also proven by Wolf using harmonic maps \cite{Wf1}.  We show how this result follows from the discussion above.  Recall from Theorem \ref{T:wolf}
that the map ${\mathcal H} : \metTeich\to \QD(S)$ defined in \eqref{E:wolfmap} is a homeomorphism. 
Choose  $t_j\to \infty$ and a sequence $\{\varphi_j\}\in \QD(S)$, $\Vert\varphi_j\Vert_1=1$.  Without loss of generality, we may assume $\{\varphi_j\}$ converges to some nonzero $\varphi\in \QD(S)$.  Let $[\sigma_j]={\mathcal H}^{-1}(t_j\varphi_j)$, and let $f_j$ be the associated harmonic maps.  By definition, the Hopf differentials of the rescaled maps $\hat f_j$ converge to $(1/2)\varphi$.  Indeed, $\Hopf(\hat f_j)=t_j\varphi_j/E(f_j)$, and by 
\eqref{E:wolf} and \eqref{E:wolf2}, $E(f_j)\sim 2t_j$.
On the other hand, the Hopf differentials of $\hat f_j$  converge to the Hopf differential of the limiting equivariant map $\HH\to X$.  By the smallness of the action of $\Gamma$ on $X$, Skora's theorem implies that $X$ is dual to a measured foliation.  This measured foliation must coincide with the horizontal foliation of $(1/2)\varphi$.  Hence, we have shown the equality \eqref{E:compact} and that the map $\mathcal H$ defined in \eqref{E:wolfmap} extends continuously as a map from $\partial\Teich$ to the sphere at infinity in $\QD(S)$. Equality \eqref{E:df} follows from \eqref{E:compact} and the fact that $\bthurston\subset\partial\chi_{df}(\Gamma)\subset\PMF$.

\section{Weil-Petersson Geometry and Mapping Class Groups}

\begin{itemize}
\item 4.1 Weil-Petersson Geodesics and Isometries
\item 4.2 Energy of Harmonic Maps
\end{itemize}

\vspace{-.25in}

\subsection{Weil-Petersson Geodesics and Isometries}

Teichm\"uller space has a length space structure given by the Teichm\"uller distance \eqref{E:teichmetric}.  An alternative Riemannian structure arises  from the description of Teichm\"uller space via hyperbolic metrics presented in Section \ref{S:metric}.  This is the Weil-Petersson metric, and its  properties continue to be the subject of much research.  In this section, we present a short review of some of the aspects of Weil-Petersson geometry that will be relevant later on.  The basic definitions as well as  properties of the Weil-Petersson completion are discussed in Section \ref{S:wp}.  In Section \ref{S:mcg}, we introduce the mapping class group, and in  Section \ref{S:classification}, we indicate how  the classification of individual mapping classes follows from the structure of Weil-Petersson geodesics.

\subsubsection{The Weil-Petersson metric and its completion.}  \label{S:wp}

Recall from Section \ref{S:metric} that 
the cotangent space $T^\ast_{[\sigma]} \Teich$ is identified
with the space of holomorphic quadratic differentials on $(S,\sigma)$.
  The
complete hyperbolic metric on $(S,\sigma)$ can be expressed in local conformal coordinates as $ds^2=\sigma(z)|dz|^2$.  Similarly, a
quadratic differential has a local expression $\varphi=\varphi(z)dz^2$.
 Then for $\varphi\in
T^\ast_{[\sigma]} \Teich$, the
 Weil-Petersson cometric is  given by \index{Weil-Petersson metric}
$\Vert\varphi\Vert_{wp}=\Vert\varphi\Vert_2$ (see \eqref{E:L2}).
While there exist a wide variety of invariant metrics, the Weil-Petersson metric is in a real sense the most useful for applications.  We refer the reader to Wolpert's recent survey \cite{Wo5}.  The two most important facts for us here are that (1)  the Weil-Petersson metric has  \emph{negative sectional curvature} and (2) it  is \emph{incomplete}.

The curvature properties of Teichm\"uller space with the Weil-Petersson or Teichm\"uller metrics have an interesting history.   It was long thought that the Teichm\"uller metric had negative curvature in the sense of  triangle comparisons (see \cite{Kravetz}).  This was disproven by Masur in \cite{M1} (see also \cite{Linch}, and more recently \cite{Iv4,  MP2, MP3, MW1}).   For the Weil-Petersson metric, the  first step was taken by Ahlfors \cite{Ah2}, who showed that the first variation of the area element induced by the hyperbolic metric vanishes.  This implies the k\"ahlerity. 
He also established the negativity of the Ricci and holomorphic sectional curvatures. The following result was established later:
\begin{theorem}[Tromba \cite{Tromba1}, Wolpert \cite{Wo2}, see also \cite{J3,Siu}.]  \label{T:wpcurvature} \index{Weil-Petersson metric!curvature of}
The curvature of the Weil-Petersson metric has
\begin{enumerate}
\item  holomorphic sectional curvatures and Ricci curvatures  bounded above by $-1/2\pi(p-1)$, and
\item negative sectional curvature.
\end{enumerate}
\end{theorem}

Incompleteness is a consequence of the nature of degenerating Riemann surfaces. This was first recognized in the work of Bers, Chu, Wolpert and Masur  (cf.\ \cite{Chu, M2, Wo1}).  \index{Weil-Petersson metric!incompleteness of}
A model for degeneration is given by the ``plumbing construction."  Here is a simple version:  let $S_1$ and $S_2$ be compact surfaces of genera $p_1, p_2$.  Choose local coordinates $z_1, z_2$ centered at points $x_1\in S_1$, $x_2\in S_2$.  Fix $0<t<1$ and construct a new surface from the following three pieces:
$S_1\setminus\{|z_1|\leq 1\}$, $S_2\setminus\{|z_2|\leq 1\}$, and the annulus $\{ (z_1, z_2) : z_1z_2=t\}$.  The boundary of the annulus is identified in the obvious way with the boundaries of the surfaces with disks deleted.  In this way, one obtains a compact Riemann surface $S_t$ of genus $p_1+p_2$.

\setlength{\unitlength}{1cm}
\begin{picture}(14,3.75)
\put(.5,2){$S_t\ =$}
\put(3.3,2.2){$ S_1$}
\put(7.6,2.3){$ S_2$}
\ifx\pdftexversion\undefined
\put(1.75,1){
{\scalebox{.35}{\includegraphics{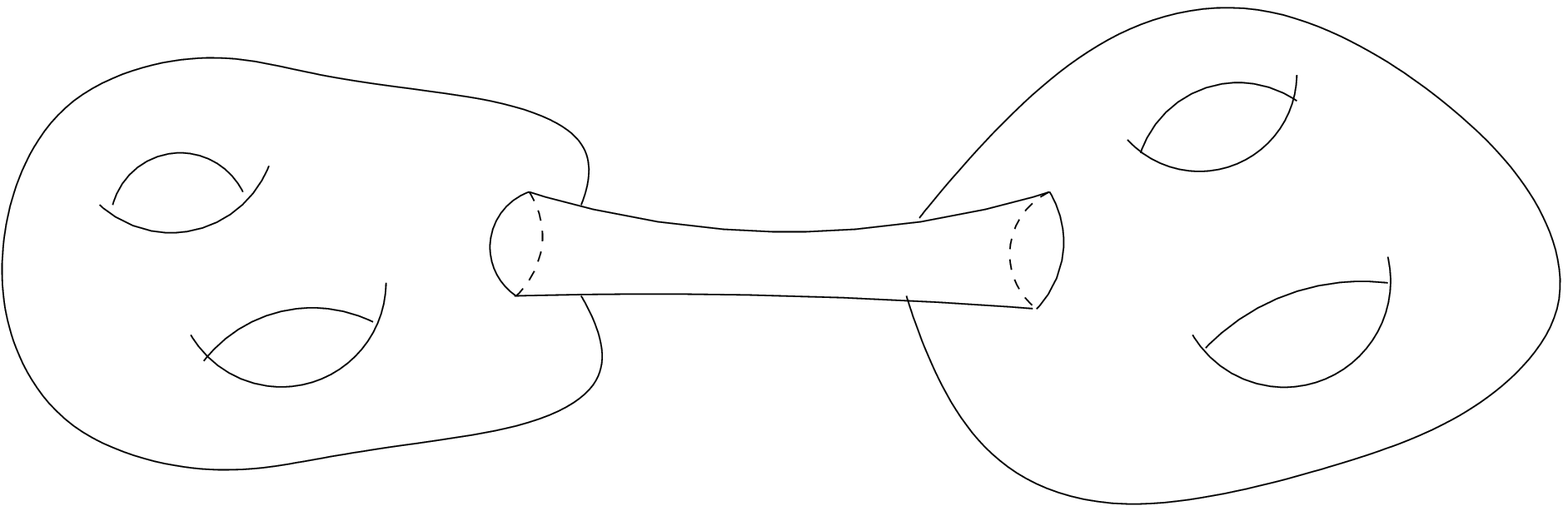}}}}
\else
\put(1.75,1){
{\scalebox{.35}{\includegraphics{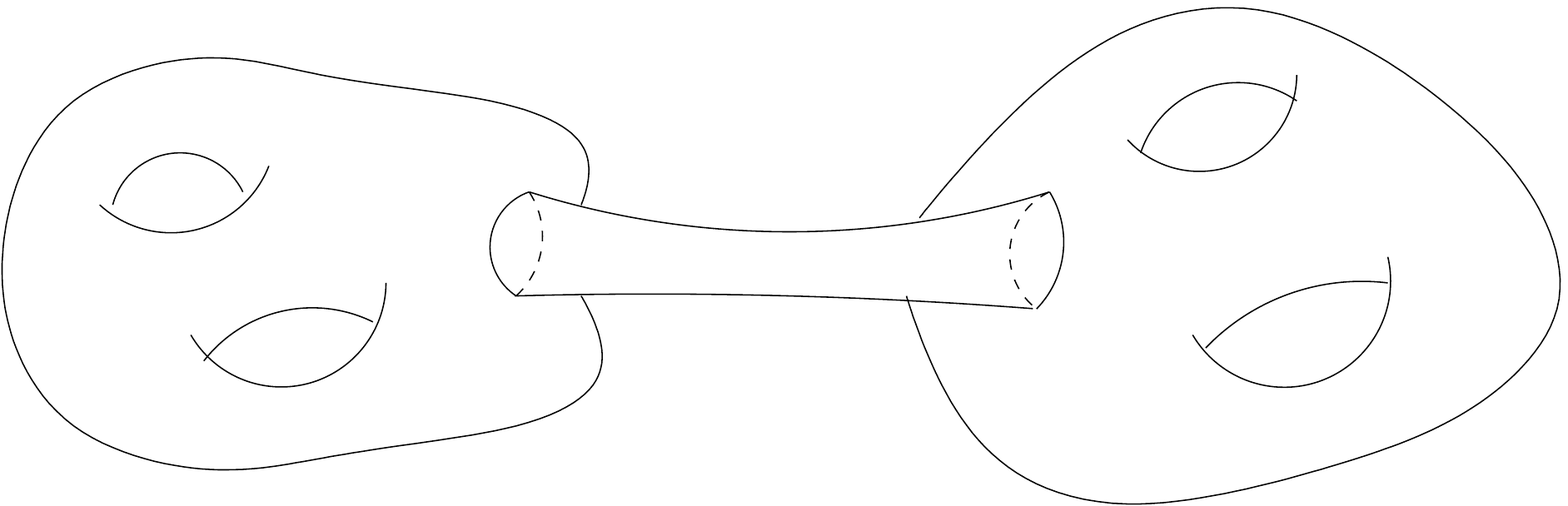}}}}
\fi
\put(4.5,0){Figure 6.}
\end{picture}
\smallskip

As $t\to 0$, the points in Teichm\"uller space corresponding to $S_t$ diverge, because the annulus is begin ``pinched."   This can also be seen from the hyperbolic geometry.  Using the maximal principle, one can approximate the behavior of the hyperbolic metric on $S_t$ (cf.\ \cite{Wo4}).  In the pinching region, it is roughly approximated by the hyperbolic metric on the annulus given by
$$
ds^2_t= \frac{|dz|^2}{|z|^2(\log|z|)^2} \frac{\Phi_{t}^2}{\left(\sin\Phi_t\right)^2}\ ,
$$
where
$
 \Phi_t=\pi\log|z|/\log|t|
$, and $z$ is either $z_1$ or $z_2$.
As a result the length $\ell$ of the ``waist" of the annulus is shrinking to zero as $t\to 0$.  In fact, the length is of order $\ell\sim 1/\log(1/t)$ (see \cite{Wo4}).  Notice that every curve passing through the annulus must then become rather long.  This is a general fact in the hyperbolic geometry of surfaces.
The following rough statement of the \emph{Collar Lemma}\index{collar lemma} indicates that around short geodesics on a hyperbolic surface one always can find long cylinders.  For a more precise statement, see \cite{Keen}.
\begin{lemma} \label{L:collar}
Let $(S,\sigma)$ be a hyperbolic surface and $c$ a simple closed geodesic of length $\ell\neq 0$.  Then  any simple closed essential curve having nonzero geometric intersection with $c$ has length on the order  $\sim\log(1/\ell)$.
\end{lemma}

The behavior of the  Weil-Petersson metric  at points in $\Teich$ described by these degenerations has the following model due to Masur, Yamada, and Wolpert (for a review, see \cite{Wo5}).
Define an incomplete metric space
\begin{equation} \label{E:modelmetric}
\model=\left\{ (\xi,\theta)\in \R^2 : \xi>0\right\} \ ,\qquad
ds^2_{\model}=4d\xi^2+\xi^6d\theta^2\ .
\end{equation}
The metric completion $\barmodel$ of $\model$ is obtained by adding a single point $\partial\model$  corresponding to
the entire real axis $\xi=0$.
The completion is then an NPC
space which is, however, not locally compact.   Indeed, an $\varepsilon$-neighborhood of $\partial \model$ contains all points of the form $(\xi,\theta)$, $\xi<\varepsilon$, and $\theta$ arbitrary.

The importance of $\barmodel$ is that it is a model for the normal space to the boundary strata. 
Let $\barTeich$ denote the metric completion of $\Teich$.
 We have  the following  local description (cf.\ \cite{M2}): \index{Weil-Petersson metric!completion of} $\partial\Teich=\barTeich\setminus\Teich$ is a disjoint union
of smooth connected strata formed  by collapsing a collection of disjoint simple closed essential
 curves  on
$S$  to  points.  
Associated to  the nodal surface is another
  Teichm\"uller
space which is by definition the set of equivalence classes of complex structures on the normalized (possibly
disconnected) surface, with the preimages of the nodes as additional marked points. It is therefore naturally isomorphic to a product of lower
dimensional Teichm\"uller spaces.
A neighborhood of a point in the boundary is then homeomorphic to an open set in the lower dimensional product crossed with as many multiples of $\barmodel$ as there are collapsed curves.  Metrically, 
the statement is that the Weil-Petersson metric in this neighborhood is equal to the product metric up to third order in the $\xi$ variables (see \cite{DW, M2, OW, Wo4, Wo5, WW, Y2}).   Moreover,  by Wolpert's theorem (see Theorem \ref{T:convexity} below) $\Teich$ with the Weil-Petersson metric is geodesically convex and the  boundary strata
 are totally geodesically embedded.    
 
The following observation is also due to Yamada: 

\begin{theorem} \label{T:wpcompletion}
  The completion $\barTeich$ of $\Teich$ is a complete NPC space.
  \end{theorem}
While this follows on general principles (cf.\ \cite{BH}), the identification of the boundary strata of the completion with lower dimensional Teichm\"uller spaces (and Weil-Petersson metrics) is especially useful.

 Let us point out  two  properties of the geometry of the Weil-Petersson
completion that are consequences of this expansion.  These were first stated by Yamada \cite{Y2}.   
The first  result, dubbed \emph{nonrefraction}\index{nonrefraction} by Wolpert,  is the statement that geodesics
 from points in Teichm\"uller space to the boundary touch the boundary only at their
endpoints (see Figure 7 (a)).
It is easy to see that this is true for the model space above.  Indeed, the equations for a unit speed geodesic $\alpha(t)=(\xi(t), \theta(t))$ in $\model$ are
\begin{align} 
 \xi\ddot \xi  &=
(3/4)\xi^6\dot{\theta}^2\ , \notag\\
\xi^6 \dot\theta &=\, \text{constant}\ , \label{E:modelgeodesic}
 \\
2|\dot\xi| \ & ,\ \xi^3|\dot\theta|
\ \leq 1\ .\notag
\end{align}
If $\xi(t)\to 0$ as $t\to 1$, say, then the second and third equations imply that the constant above must vanish.  In other words, $\theta(t)$ is constant and $\xi(t)$ is linear.
The proof of the statement for geodesics in $\barTeich$ involves a scaling argument to approximate geodesics in $\barTeich$ by corresponding geodesics in the model space.  The third order approximation of the Weil-Petersson metric by the model metric is sufficient to show that the approximation of geodesics is also to high order, and the qualitative behavior of geodesics in $\barTeich$ is the same as for the model space.

\setlength{\unitlength}{1cm}
\begin{picture}(14,4)
\ifx\pdftexversion\undefined
\put(.5,1.3){
{\scalebox{.5}{\includegraphics{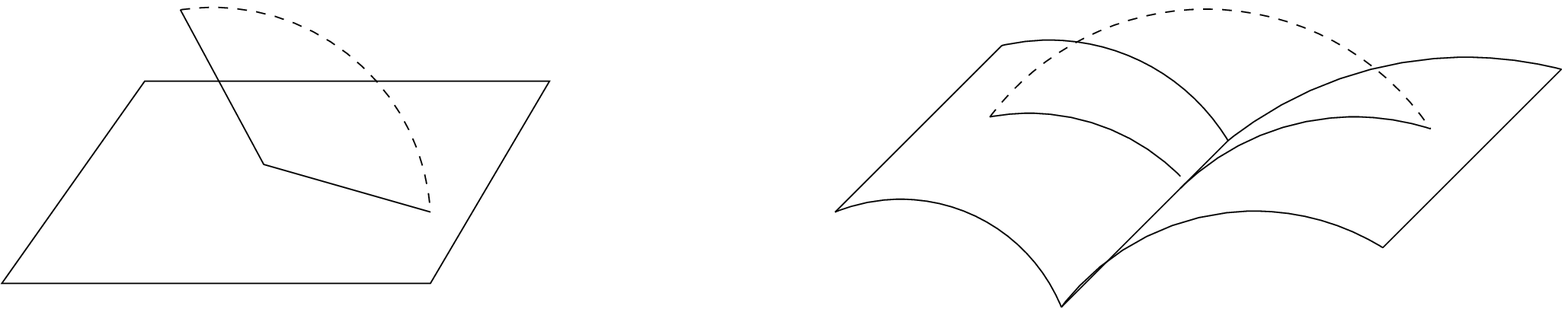}}}}
\else
\put(.5,1.3){
{\scalebox{.5}{\includegraphics{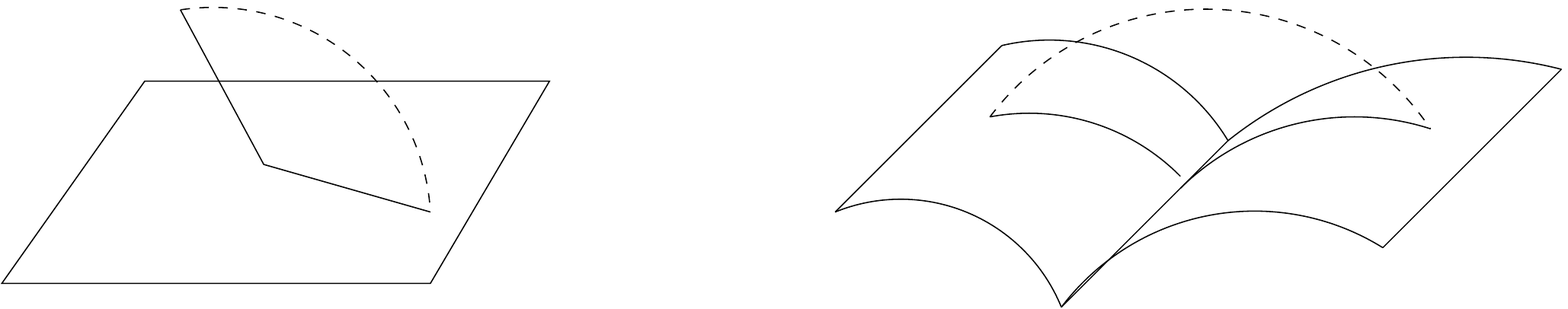}}}}
\fi
\put(0,2.3){$\partial\Teich$}
\put(5.7,2.5){${\mathcal D}[c_1]$}
\put(10.4,1.8){${\mathcal D}[c_2]$}
\put(2.3,.5){(a)}
\put(8.2,.5){(b)}
\put(4.5,0){Figure 7.}
\end{picture}

\smallskip

Another application of this approximation gives 
the second important result: the different strata of the boundary of $\barTeich$ intersect transversely.
For example,  consider disjoint nonisotopic simple closed essential curves $\{c_1, c_2\}$ on
a closed compact surface $S$ with isotopy classes $[c_1]$, $[c_2]$.  Let
$[\sigma_1]$ denote a point in the boundary component ${\mathcal D}[c_1]$ of $\barTeich$ corresponding to pinching $c_1$.  Similarly, let $[\sigma_2]$ denote a
point in the boundary component ${\mathcal D}[c_2]$ of $\barTeich$ corresponding to pinching $c_2$.   Since $c_1$ and $c_2$ are disjoint, the
intersection of the closures
$\overline {\mathcal D}[c_1]\cap\overline {\mathcal D}[c_2]$ is nonempty, and in fact contains
${\mathcal D}([c_1], [c_2])$, the stratum where both $c_1$ and $c_2$ are pinched.  In particular, there is a path in $\barTeich$ from $\sigma_1$ to $\sigma_2$, lying completely in the boundary, which corresponds to
first pinching $c_2$, and then ``opening up" $c_1$.   The theorem states that this path has a ``corner" at its intersection with ${\mathcal D}([c_1], [c_2])$, and  is therefore not length minimizing.  In fact, the geodesic from
$[\sigma_1]$ to $[\sigma_2]$ intersects the  boundary of $\barTeich$ only in its endpoints (see Figure 7 (b)).


\subsubsection{The mapping class group.}  \label{S:mcg}


Denote by $\Diff$ (resp.\ $\Diffp$) denote the group of smooth diffeomorphisms (resp.\  orientation preserving diffeomorphisms)  of $S$ with the
smooth topology.
Recall that $\Diffo$ denotes the identity component
of $\Diff$, that is, the group of all diffeomorphisms isotopic to the
identity. The {\em mapping class group\/} \index{mapping class!group}  of $S$ is the quotient
\begin{equation*}
\mcg = \Diffp/\Diffo\ .
\end{equation*}
See \cite{FMar}
for a recent survey on mapping class groups.
From any of the several definitions of Teichm\"uller space given previously, it is clear that $\mcg$ acts on $\Teich$.  The first important result about this action is the following
\begin{theorem} \label{T:fricke}
The mapping class group acts properly discontinuously on $\Teich$.
\end{theorem}
This result is commonly attributed to Fricke.  One method of proof follows  from the general fact that the 
action of $\Diff$ on the space $\met$ of smooth Riemannian metrics is properly discontinuous (cf.\ \cite{Ebin, EE}). In particular the restriction to the action of $\Diffo$  on $\hypmet$ is properly discontinuous, and Teichm\"uller space, which is  the quotient
$\metTeich=\hypmet/\Diffo$, inherits such an action of $\mcg$.

Diffeomorphisms of $S$ determine automorphisms of  $\Gamma=\pi_1(S)$ as follows.
Let $x_0\in S$ be a fixed basepoint.
A diffeomorphism $\phi:S\to S$ determines an automorphism of the fundamental
group $\pi_1(S,x_0)$ if $\phi(x_0)=x_0$.
Now any diffeomorphism
 is isotopic to  one which fixes $x_0$.
Different choices of isotopy define automorphisms of $\Gamma$ which
differ by an  inner automorphism.  Hence, there is a homomorphism
$
\Diff/\Diffo \to \Out
$
where $\Out$ is the outer automorphism group of $\pi_1(S)$. \index{mapping class!as outer automorphism}
\begin{theorem}[Dehn-Nielsen \cite{Nielsen1}]\label{T:dn}
The homomorphism described above gives an isomorphism $ \Diff/\Diffo \simeq\Out$.
\end{theorem}

The quotient $\Mod=\Teich/\mcg$ is the classical Riemann moduli space. \index{Riemann moduli space} Since by Teichm\"uller's Theorem  $\Teich$ is contractible (and in particular, simply connected), $\mcg$ may be regarded as the fundamental group of $\Mod$.  However, $\mcg$ does not quite act freely, so this interpretation holds only in the \emph{orbifold} sense.  Indeed,  $\Mod$ is actually simply connected \cite{Mac}.\index{Riemann moduli space!is simply connected}
 The compactification 
 $$\barMod=\barTeich/\mcg$$
  is homeomorphic to the \emph{Deligne-Mumford compactification}\index{Riemann moduli space!Deligne-Mumford compactification} of $\Mod$, and it is a projective algebraic variety (with orbifold singularities) \cite{DM}.  An important measure of the interior regions of $\Mod$ is given by the Mumford-Mahler compactness theorem:  

\begin{theorem}[Mumford \cite{Mumford}] \label{T:mm} \index{Mumford-Mahler compactness}
The set $\Mod_\varepsilon\subset\Mod$ consisting of equivalence classes of Riemann surfaces where the hyperbolic lengths of all closed geodesics are bounded below by $\varepsilon>0$ is compact.
\end{theorem}

Note the condition in the theorem is $\Diff$-invariant and so is valid on the moduli space $\Mod$.  The corresponding result is, of course,  not true for $\Teich$ because of the proper action of the infinite discrete group $\mcg$.  For example, the orbit of a point in $\Teich$ by $\mcg$ is unbounded, but projects to a single point in $\Mod$.  This, however, is the \emph{only} distinction between $\Teich$ and $\Mod$, and so the Mumford-Mahler compactness theorem can be used effectively to address convergence questions in $\Teich$ as well. 

We illustrate this by proving a fact that will be useful later on.
Given a simple closed curve $c\subset S$, let $\ell_c[\sigma]$ denote the length of the geodesic in the homotopy class of $c$ with respect to the hyperbolic metric $\sigma$.  Note that this is independent of the choice of $\sigma$ up to the action of $\Diffo$.  Hence, $\ell_c$ gives a well-defined function
\begin{equation} \label{E:length}
\ell_c: \Teich\lra \R^+
\end{equation}
Then we have the following

\begin{corollary} \label{C:mm}
If $[\sigma_j]$ is a sequence in $\Teich$ contained in no compact subset then there is a simple closed curved $c\subset S$ such that $\ell_c[\sigma_j]$ is unbounded.
\end{corollary}

\begin{proof}  For a point $[\sigma]\in \Teich$, let $[[\sigma]]\in\Mod$ denote the corresponding point in $\Mod$.  Without loss of generality, we may assume $[\sigma_j]$ has no convergent subsequence in $\Teich$.  The same may or may not be true for the sequence $[[\sigma_j]]\subset\Mod$.  Indeed,
by Theorem \ref{T:mm}, there are two cases:  (1) there are elements $[\phi_j]\in \mcg$ and a point $[\sigma_\infty]\in \Teich$  such that $[\phi_j][\sigma_j]\to [\sigma_\infty]$ (after passing to a subsequence); (2) there are simple closed curves $c_j$ such that $\ell_{c_j}[\sigma_j]\to 0$ (after passing to a subsequence).  In the first case, our assumptions imply that infinitely many $[\phi_j]$ are distinct.  It follows that there is a simple closed curve   $c$ such that $\ell_{f_j(c)}[\sigma_\infty]\to \infty$.  But then $\ell_c[\sigma_j]\to \infty$, as desired.  In the second case, 
 we may assume $c_j$ converges projectively to a nontrivial measured foliation $\mathcal F$ (see Section \ref{S:trees}).  If $c$ is any simple closed curve with $i([c],{\mathcal F})\neq 0$, then $i(c, c_j)\neq 0$ for $j$ large.  But since $\ell_{c_j}[\sigma_j]\to 0$, $\ell_c[\sigma_j]\to \infty$ by the Collar Lemma \ref{L:collar}.
\end{proof}

Thurston's classification\index{mapping class!Thurston classification} of surface diffeomorphisms may be described in terms of the natural action of $\mcg$ on $\MF$ and $\PMF$:
an element $[\phi]\in\mcg$ is called \emph{reducible}\index{mapping class!reducible} if $[\phi]$ fixes (up to isotopy) some collection of disjoint simple closed essential curves on $S$. 
 It is called
\emph{pseudo-Anosov}\index{mapping class!pseudo-Anosov} if there is  $r>1$ and transverse measured
foliations
${\mathcal F}_+$, ${\mathcal F}_-$ on $S$ such that $[\phi] {\mathcal F}_+$ is measure equivalent to $r {\mathcal F}_+$, and $[\phi] {\mathcal F}_-$ is measure equivalent to $r^{-1}{\mathcal F}_-$.  
 ${\mathcal F}_+$ and ${\mathcal F}_-$ are called the \emph{stable} and \emph{unstable} foliations of $[\phi]$, respectively. The classification states that
any
$[\phi]\in\mcg$ is either \emph{periodic} (i.e.\ finite order),\index{mapping class!pseudo-Anosov} infinite order and reducible, or pseudo-Anosov.  Moreover, these are mutually exclusive
possibilities.


\subsubsection{Classification of Weil-Petersson isometries.} \label{S:classification}


We now indicate how the Thurston classification of mapping classes is mirrored by the Weil-Petersson geometry.
The action of $\mcg$ on $\Teich$ is isometric with respect to the Weil-Petersson metric.  Conversely, every Weil-Petersson isometry is essentially given by a mapping class (see \cite{BM, MP2, MW2, Wo5}). \index{Weil-Petersson metric!isometry group of} \index{mapping class!as Weil-Petersson isometry} Since the Weil-Petersson metric has negative curvature it is a natural to classify individual mapping classes in a manner similar to isometries of Cartan-Hadamard manifolds.   

\begin{theorem}[Daskalopoulos-Wentworth \cite{DW}, Wolpert \cite{Wo5}] \label{T:dw} \index{Weil-Petersson metric!geodesics}
If $[\phi]\in\mcg$ is infinite order and irreducible, then there is a unique $[\phi]$-invariant complete Weil-Petersson geodesic in $\Teich$.
\end{theorem}

Here is a very rough idea of proof of this result.  Let $\tilde \alpha_j : [0,1]\to \Teich$, $\tilde\alpha_j(1)=[\phi] \tilde\alpha_j(0)$ be a sequence of curves minimizing the translation length of $[\phi]$, i.e.\
$$
\lim_{j\to\infty} \int_0^1 \Vert\dot{\tilde \alpha}_j\Vert_{wp} dt=\twp[\phi]\ .
$$
Let $\alpha_j: S^1\to \barMod$ be the projection of $\tilde\alpha_j$.  Since $\barMod$ is compact one can show using Ascoli's Theorem that, after passing to a subsequence, $\alpha_j$ converge uniformly to some curve $\alpha: S^1\to \barMod$.    The trick now is to show that this curve admits a lift $\tilde \alpha: [0,1]\to \barTeich$, $\tilde \alpha(1)=[\phi]\tilde \alpha(0)$.  Then $\tilde \alpha$ must be an invariant geodesic.  Since $[\phi]$ is irreducible, by the nonrefraction results we know that $\tilde\alpha$ must have image in $\Teich$.  The existence of a lift is not obvious, since $\barTeich\to \barMod$ is ``branched" to infinite order along the boundary $\partial\Teich$.  One needs to exploit the fact that $\alpha$ is the \emph{limit} of curves that are liftable.  We refer to \cite{DW} for more details.

The existence of invariant geodesics for infinite order irreducible mapping classes allows for the precise 
classification of  Weil-Petersson isometries in terms of translation length  that we have given in Table 1.  
For $[\phi]\in \mcg$, define the  \emph{Weil-Petersson  translation length} by
\begin{equation} \label{E:wptl} \index{translation length!Weil-Petersson}
\twp[\phi]=\inf_{[\sigma]\in \Teich} \dwp([\sigma],[\phi][\sigma])\ .
\end{equation}

\begin{table}[h] 
\centering
\begin{tabular} {|c||c|c|} \hline
 &  { semisimple} &  {  not semisimple} \\ \hline\hline
$ \twp=0$  &  {  periodic} &  {   strictly pseudoperiodic} \\ \hline 
$\twp\neq 0$  &   { infinite order irreducible} &   {  reducible but not pseudoperiodic} \\ 
 \hline
\end{tabular}
\bigskip
\caption{\emph{Classification of Weil-Petersson Isometries.}}
\end{table}

First, let us
clarify the terminology used there:  
$[\phi]\in\mcg$ is  \emph{pseudoperiodic} if it is either periodic, or it is reducible and periodic on the reduced components;\index{mapping class!pseudoperiodic}  it is
strictly pseudoperiodic if it is pseudoperiodic but
not periodic.  Furthermore, we say that $[\phi]$ is  \emph{semisimple}\index{semisimple!mapping class}\index{mapping class!semisimple} if there is  $[\sigma]\in\Teich$ such that
$\twp[\phi]=\dwp([\sigma],[\phi][\sigma])$.

Here is a sketch of the  proof: first, note that 
it is a consequence of Theorem \ref{T:fricke} that $\twp[\phi]=0$ if and only if $[\phi]$ is pseudoperiodic.  The first row of Table 1 then follows from this and the fact that $[\phi]$ has a fixed point in $\Teich$ if and only if
$[\phi]$ is periodic. 
 If $[\phi]$ is infinite order irreducible, then as a consequence of Theorem \ref{T:dw},
  $\twp[\phi]$ is
attained along an invariant geodesic, so these mapping classes are semisimple.
Conversely, suppose the translation length is attained at $[\sigma]\in\Teich$, but $[\phi][\sigma]\neq[\sigma]$.  
Then we argue as in Bers \cite{Bers2} (see also, \cite[p.\ 81]{BGS}) to show that the geodesic from $[\sigma]$ to $[\phi][\sigma]$, which exists by the geodesic convexity of the Weil-Petersson metric,  may be extended to a complete $[\phi]$-invariant
geodesic.  On the other hand, if there is a complete, nonconstant Weil-Petersson geodesic in $\Teich$ that is invariant with respect to a mapping class $[\phi]\in \mcg$, the negative curvature implies that 
$[\phi]$  must be infinite order and irreducible.

It is worth mentioning that no properties of  pseudo-Anosov's other than the fact that they have infinite order and are
irreducible were used in the proof above.  In particular, the description given in Table 1 is independent of Thurston's classification.

We point out 
a further property of the axes of pseudo-Anosov's.

\begin{theorem}[Daskalopoulos-Wentworth \cite{DW}, Wolpert \cite{Wo5}] \label{T:axes}
Let $A_{[\phi]}$ and  $A_{[\phi']}$ be the axes for independent
pseudo-Anosov mapping classes $[\phi]$ and $[\phi']$.  Then
$A_{[\phi]}$ and $A_{[\phi']}$ diverge.
\end{theorem}

This result is also not completely obvious because of the noncompleteness of $\Teich$.  More to the point, there exist  
\emph{flats}, i.e.\  a totally geodesically embedded copy of $\R^m\hookrightarrow\barTeich$.   which potentially hinder the divergence.
A much more detailed discussion of asymptotics of complete Weil-Petersson geodesics is forthcoming (see \cite{BMM}).

\subsection{Energy of Harmonic Maps}

In this section we return to harmonic maps and show how they can be used to probe the action of the mapping class group on Teichm\"uller space.  In Section \ref{S:nielsen}, we discuss Nielsen's realization problem for finite subgroups of the mapping class group.  In Section \ref{S:proper}, we introduce two classes of functions on Teichm\"uller space that are constructed using the energy of harmonic maps, and we indicate when these functions are proper.  In Section \ref{S:convex}, we discuss the convexity of one of the two classes and show how this resolves the Nielsen conjecture.  We also state Wolpert's result on convexity of length functions.  Finally, in Section \ref{S:applications}, we indicate some other applications of the energy functionals.

\subsubsection{Nielsen realization.}  \label{S:nielsen} \index{Nielsen realization problem}

  Here we discuss the classical question of Nielsen \cite{Nielsen2}.
The exact sequence 
\begin{equation} \label{E:exact}
1\lra \Diffo\lra \Diffp
{\buildrel \pi\over\longrightarrow}\
\mcg\lra 1
\end{equation}
which defines the mapping class group does not split in general (see \cite{Markovic, Morita1, Morita2}). The \emph{realization problem} asks for which subgroups  $G\subset \mcg$ does there exist a homomorphism $\jmath: G\to\Diffp$ such that $\pi\circ\jmath=\id$.

Let $S$ be a closed Riemann surface of negative Euler characteristic.  Then we have the following two important facts.  First, if $\phi$ is a holomorphic automorphism of $S$ homotopic to the identity, then $\phi$ is in fact \emph{equal} to the identity.  Indeed, if this were not the case then since complex curves in a complex surface intersect positively  the number of fixed points of $\phi$, counted with multiplicity, would necessarily be positive.  On the other hand, if $\phi\sim\id$, then by
 the Lefschetz fixed point theorem the total intersection number is just the Euler characteristic of $S$, which we have assumed is negative.

From this fact we arrive at Fenchel's observation that 
 if a subgroup $G\subset \mcg$ fixes a point $[j]\in\Teich$, then $G$ can be realized as the automorphism group of a Riemann surface $(S,j)$ with $j$ in the class $[j]$.   For if $ \phi_1, \ldots, \phi_m$ are holomorphic lifts to $\Diffp$ of generators $[\phi_1], \ldots, [\phi_m]$ of $G$, then  any relation on the $[\phi_j]$'s, applied to the $ \phi_j$'s, is a holomorphic map $\sim \id$, and so by the previous paragraph the relations in the group also lift.
  In particular, \eqref{E:exact} splits over $G$.

The second fact is that the automorphism group of a Riemann surface of genus $p\geq 2$ is finite.  This is because on the one hand it is the isometry group of the hyperbolic metric, which is compact, and on the other hand it is discrete, since there are no holomorphic vector fields.
Hence, any subgroup of the mapping class group which fixes a point in Teichm\"uller space is finite and \eqref{E:exact} splits over it.   These two facts motivate the  following result, which is known as the \emph{Nielsen Realization Theorem}.
\begin{theorem}[Kerckhoff \cite{Kf}] \label{T:nielsen}
The sequence \eqref{E:exact} splits over all  finite subgroups of $\mcg$\ .
\end{theorem}
From the discussion above, the idea of the  proof is to show the following
\begin{theorem}  \label{T:fixed}
Let $G\subset \mcg$ be a finite subgroup of the mapping class group.  Then $G$ has a fixed point in $\Teich$.
\end{theorem}
\noindent   
The complete  proof of Theorem \ref{T:nielsen} was first obtained by Kerckhoff in \cite{Kf} and later by Wolpert \cite{Wo3}. Both proofs proceed via Theorem \ref{T:fixed}.
 Partial results had been found earlier by Fenchel \cite{F1,F2} and Zieschang \cite{Z}.  See also Tromba \cite{Tromba3}.

\subsubsection{Properness of the energy.}  \label{S:proper}

Let $M$ be an arbitrary compact Riemannian manifold and $S$ a closed hyperbolic surface with negative Euler characteristic.  
Now if $\rho: \pi_1(M)\to \pi_1(S)$ is a given homomorphism it follows by Theorem \ref{T:es} that there is a harmonic map $f: M\to S$ such that the induced action $f_\ast: \pi_1(M)\to \pi_1(S)$ coincides with $\rho$.  The energy $E(f)$ then depends only on the equivalence class of hyperbolic metrics $[\sigma]\in \Teich$ (see Theorem \ref{T:homotopy}).  In other words, there is a well-defined function
$$
{\mathcal E}^+_\rho: \Teich\longrightarrow\R^+\ .
$$
The existence of a minimum is in turn a reflection of the homomorphism $\rho$.  One way to guarantee a minimum is to show that ${\mathcal E}^+_\rho$ diverges at infinity.
In this context, we have the following

\begin{proposition} \label{P:range}
 If $\rho$ is surjective  then the associated function ${\mathcal E}^+_\rho$ is proper.  
\end{proposition}

\begin{proof}
This is  easy to see, given the Lipschitz bound Proposition \ref{P:lipschitz} and the Mumford-Mahler Compactness Theorem \ref{T:mm} (or more precisely, Corollary \ref{C:mm}).  Indeed, if ${\mathcal E}^+_\rho$ is not proper, there is a sequence $[\sigma_j]$ and harmonic maps $f_j:M\to (S,\sigma_j)$ in the homotopy class defined by $\rho$,  such that $E[\sigma_j]\leq B$ for some constant $B$.  Furthermore, we may assume there is  a simple closed curve $c$ with $\ell_{c}[\sigma_j]\to \infty$. Let $s$ be a closed curve in $M$ with $f_j(s)$ homotopic to $c$.  Then since the $f_j$ are uniformly Lipschitz,
$$
\ell_c[\sigma_j]\leq \length(f_j(s))\leq \tilde B \length(s)\ .
$$
Since the right hand side is fixed independent of $j$ and the left hand side diverges with $j$, we derive a contradiction.
\end{proof}

The superscript $+$ on ${\mathcal E}^+_\rho$ is to remind us that this is a function of the hyperbolic metric on the \emph{target}.
 It is also interesting to consider the energy as a function of the domain metric (cf.\ \eqref{E:eminus}).
Let $M$ be a compact Riemannian manifold with nonpositive sectional curvature.  Let $S$ be a closed surface and let $\rho:\pi_1(S)\to \pi_1(M)$ is a homomorphism.  Then 
for each complex structure $\sigma$ on $S$ there is a harmonic map $f: (S,\sigma)\to M$ whose induced action on $\pi_1$ coincides with $\rho$.  
The energy of this map gives a well-defined function
$$ 
{\mathcal E}^-_\rho: \Teich\longrightarrow\R^+\ .
$$

Again, the existence of minima can be deduced from the properness of this functional.  The following can be proved using the same ideas as in the proof of Proposition \ref{P:range} and Corollary \ref{C:mm}.
\begin{proposition}[see \cite{SU2,SY2}] \label{P:domain}
If $\rho$ is injective then ${\mathcal E}^-_\rho$ is proper.
\end{proposition}

As we have seen in Theorem \ref{T:diff}, the function ${\mathcal E}^-_\rho $ is differentiable. From the discussion in Section \ref{S:metric} (see esp. \eqref{E:nondegenerate}), 
critical points correspond to \emph{conformal} harmonic maps, i.e.\ those for which the Hopf differential vanishes.  According to Sacks-Uhlenbeck \cite{SU2}, these are branched minimal surfaces in $M$.

There is a remarkable connection between these functionals and  the Weil-Petersson metric.  If we take $M=(S,\sigma_0)$ for some hyperbolic metric $\sigma_0$ and $\rho=\id$,  we have defined two functions ${\mathcal E}^\pm_{\id}$ on $\Teich$, both of which clearly have critical points  at $[\sigma_0]$.   We have 

\begin{theorem}[Tromba \cite{Tromba2}, Wolf \cite{Wf1}, Jost \cite{J2}] \label{T:trombawolf}
 The second variation of either ${\mathcal E}^\pm_{\id}$ at $[\sigma_0]$ is a positive definite hermitian form on $T_{[\sigma_0]}\Teich$ which \emph{coincides} with the Weil-Petersson metric.
\end{theorem}

A critical point of ${\mathcal E}^-_\id $ is a  holomorphic map $(S,\sigma)\to (S,\sigma_0)$ homotopic to the identity.  As argued  in Section \ref{S:nielsen}, this must be the identity and $\sigma=\sigma_0$.
We conclude that ${\mathcal E}^-_\id$ is a proper function on $\Teich$ with a unique critical point.   By Theorem \ref{T:trombawolf}, it is also nondegenerate.  It follows that $\Teich$ is diffeomorphic to $\R^n$; hence, we have a fourth  a proof of Theorem \ref{T:teich} (see \cite {FT}). \index{Teichm\"uller!theorem}


\subsubsection{Convexity of energy and length functionals.} \label{S:convex}


\begin{theorem}[Tromba \cite{Tromba3}, Yamada \cite{Y1}] \label{T:tromba}
 The energy ${\mathcal E}^+_\rho$ defined above is strictly convex along Weil-Petersson geodesics.
\end{theorem}
This result was first obtained by Tromba in the case where $M$ is homeomorphic to $S$.  It was later generalized to the statement above by Yamada.
It follows that the minimum of ${\mathcal E}^+_\rho$ is unique if it exists.

The conclusion is that there exists an abundance of convex exhaustion functions on Teichm\"uller space and an explicit method to construct them.  
  Any one of these gives a solution to the Nielsen problem!  \index{Nielsen realization problem}
For the average of such a function over a finite subgroup  $G\subset\mcg$ is again strictly convex \emph{and} $G$-invariant.  Hence, its unique minimum is also $G$-invariant, i.e.\ a fixed point of $G$, and Theorem \ref{T:fixed} is proven.  The easiest example is to take $M=(S,\sigma_0)$,  for any complex structure $\sigma_0$, and $\rho=\id$, as in the previous section.

It turns out that the analogous statement Theorem \ref{T:tromba} for  ${\mathcal E}^-_\rho$ is false.  For example, we could choose $M$ to be a fibered hyperbolic 3-manifold with $S$ a fiber and $\rho$ the homomorphism coming from the inclusion.  Then $\rho$ is invariant by conjugation of the monodromy of the fibration, which by a theorem of Thurston is represented by a pseudo-Anosov diffeomorphism (cf.\ \cite{Ot}).  In particular, it has infinite order.  This fact leads to infinitely many minima of ${\mathcal E}^-_\rho$, whereas if ${\mathcal E}^-_\rho$ were strictly convex, it would have a unique minimum.
 It is certainly an interesting question to find conditions where convexity holds.  
 
A very special case of the previous discussion is when $M$ is a circle.  Harmonic maps from a circle correspond to geodesics.  Historically, geodesic length functions were considered before the energy of harmonic maps from higher dimensional domains.  In particular, we have the following important result of Wolpert.

\begin{theorem}[Wolpert \cite{Wo3, Wo6}]  \label{T:convexity} 
For any simple closed curve $c$, the function $\ell_c:\Teich\to \R^+$ defined in \eqref{E:length} is strictly convex along Weil-Petersson geodesics.  The extension of the length function to geodesic currents is also strictly convex.
\end{theorem}

One consequence of this is the geodesic convexity of Teichm\"uller space, \index{Weil-Petersson metric!geodesic convexity of} i.e.\ between any two points in $\Teich$ there exists a unique Weil-Petersson geodesic.  One can also construct convex exhaustion functions, although in a manner slightly different from that of the previous section.  If we choose a collection $c_1,\ldots, c_m$ of simple closed curves which are \emph{filling} in  $S$\index{filling} in the sense that any other simple closed essential curve has nontrivial intersection with at least one $c_j$, then the function
 \begin{equation}\label{E:wolpert}
 \beta=\ell_{c_1}+\cdots+\ell_{c_m}\ ,
 \end{equation}
  is an exhaustion function.  This again follows by Mumford-Mahler compactness and the Collar Lemma.  Since $\beta$ is also strictly convex, this gives a solution to Nielsen's problem as above, and indeed this is Wolpert's method.

Finally, we point out that
Kerckhoff's proof of  Theorem \ref{T:nielsen} was the first to lay out this type of  argument.  The difference is that he proved convexity not with respect to the Weil-Petersson geometry but  along Thurston's earthquake deformations.


\subsubsection{Further applications.} \label{S:applications}


We now enumerate some other applications of the ideas developed in previous section.

\medskip\noindent $\bullet$ \emph{Convex cocompact representations.} \index{convex cocompact}
Note that Proposition \ref{P:domain} can also be adapted to the equivariant case and metric space targets.  Here, $\rho:\Gamma=\pi_1(S)\to \iso(X)$, where $X$ is  a simply connected NPC space.  Injectivity is replaced by the condition that the translation length of any  isometry in the image is bounded below by a uniform constant.

A discrete embedding $\rho:\Gamma\to \iso(X)$ is {\em convex cocompact\/}
if there exists a $\rho$-invariant closed geodesically convex subset
$N\subset X$ such that $N/\rho{\Gamma}$ is compact.  

\begin{theorem}[Goldman-Wentworth \cite{GW}] \label{T:gw}
$\mcg$ acts properly discontinuously  on the space of convex cocompact embeddings $\rho:\Gamma\to \iso(X)$.
\end{theorem}
When $\iso(X)=\PSL(2,\C)$, a convex cocompact representation is {\em quasi-Fuchsian,\/}\index{quasi-Fuchsian representation}
that is, a discrete embedding whose action on $S^2= \partial\HH^3$
is topologically conjugate to the action of a discrete subgroup of
$\PSL(2,\C)$. In this case,
Theorem \ref{T:gw} is just the known fact that $\mcg$ acts properly on
the space $\qf$ of quasi-Fuchsian embeddings. Indeed,
Bers' simultaneous uniformization theorem \cite{Bers} provides
a $\mcg$-equivariant homeomorphism
\begin{equation*}
\qf \longrightarrow \Teich\times \Teich.
\end{equation*}
Properness of the action of $\mcg$  on $\Teich$, Theorem \ref{T:fricke},  implies properness 
on $\qf$.

The idea of the proof of Theorem \ref{T:gw} is to show that if $\rho$ is convex cocompact
\begin{enumerate}
\item  then there exists a $\rho$-equivariant harmonic map $f:\widetilde S\to X$,
\item  and the corresponding energy functional ${\mathcal E}^-_\rho: \Teich\to \R^+$ is proper.
\end{enumerate}
Then one associates to each $\rho$ the compact subset of minima of ${\mathcal E}^-_\rho$ in $\Teich$, and properness of the action of $\mcg$ on $\Teich$ implies the result. See \cite{GW} for more details.

\medskip\noindent $\bullet$ \emph{Filling foliations.} Recall from Section \ref{S:skora} that by the Hubbard-Masur Theorem any measured foliation can be realized as the horizontal foliation of a holomorphic quadratic differential.  As a second application, consider the problem of realizing a \emph{pair} of measured foliations as the horizontal and vertical foliations of a single quadratic differential on some Riemann surface.  A pair ${\mathcal F}_+, {\mathcal F}_-$ of measured foliations on $S$ is called \emph{filling}\index{filling} if for any third measured foliation $\mathcal G$
$$
i({\mathcal F}_+,{\mathcal G})+i({\mathcal F}_-,{\mathcal G})\neq 0\ ,
$$
where $i(\cdot,\cdot)$ denotes the intersection number (see Section \ref{S:trees}).
 \begin{theorem}[Gardiner-Masur \cite{GM}]
 ${\mathcal F}_+, {\mathcal F}_-$ are filling if and only if there is a complex structure $j$ and a holomorphic quadratic differential $\varphi$ on $(S,j)$ such that ${\mathcal F}_+$ and ${\mathcal F}_-$ are measure equivalent to the vertical and horizontal foliations of $\varphi$, respectively.  Moreover, $[j]\in \Teich$, and $\varphi$ for each $j\in[j]$, are uniquely determined by ${\mathcal F}_\pm$.
\end{theorem}
It is relatively easy to see that the  horizontal and vertical trajectories of a holomorphic quadratic differential are filling (cf.\ \cite[Lemma 5.3]{GM}). 
The proof of the  converse follows by showing, using arguments similar to those in the proof of Proposition \ref{P:domain},  that $\ext_{{\mathcal F}_+}+\ext_{{\mathcal F}_-}$ is a proper function on $\Teich$.
The first variational formula Theorem \ref{T:diff} shows that a local minimum is a point at which the quadratic differentials for ${\mathcal F}_+$ and ${\mathcal F}_-$ are related by a minus sign.  
On the other hand, by the argument in Section \ref{S:skora},   ${\mathcal F}_\pm$ are therefore vertical and horizontal foliations of one and the same differential.    Uniqueness can also be proven by analytic methods (see \cite{W2} for more details).

\medskip\noindent $\bullet$ \emph{Holomorphic convexity of $\Teich$.}
The convex exhaustion functions constructed in the previous sections are, in particular, strictly plurisubharmonic (Tromba \cite{Tromba4} showed that this is true for ${\mathcal E}^-_\id$ as well).  This gives a new proof of the following

\begin{theorem}[Bers-Ehrenpreis \cite{BE}] \label{T:bers}\index{Teichm\"uller!space!is Stein}
Teichm\"uller space is a Stein manifold.
\end{theorem}

By a slight modification of length functions, we also have

\begin{theorem}[Yeung \cite{SKY}]  $\Teich$ admits a bounded strictly plurisubharmonic function.
\end{theorem}

Explicitly, one may take  $-\beta^{-\varepsilon}$, where $\beta$ is the function in \eqref{E:wolpert} and $0<\varepsilon<1$.  The existence of a bounded plurisubharmonic function has important implications for the equivalence of invariant metrics on Teichm\"uller space (see \cite{Chen, LSY, SKY}).


\section{Harmonic Maps to Teichm\"uller Space}

\begin{itemize}
\item 5.1 Existence of Equivariant Harmonic Maps
\item 5.2 Superrigidity
\end{itemize}


\vspace{-.25in}

\subsection{Existence of Equivariant Harmonic Maps}

In many ways this last chapter  combines ideas from all of the previous ones.  Because of the nonpositive curvature of  the Weil-Petersson metric, harmonic maps with Teichm\"uller space as a target have good regularity properties.   The isometry group is the mapping class group, so the equivariant problem gives a way to study representations of fundamental groups to $\mcg$.  Since the Weil-Petersson metric is not complete, we need to pass to the completion $\barTeich$ and use the theory of singular space targets of Gromov-Korevaar-Schoen. In Section \ref{S:maps}, we show how the results of Section \ref{S:classification} can be used to prove existence of equivariant harmonic maps to $\barTeich$, and in Section \ref{S:surface}, we state a result on the regularity of energy minimizing maps for surface domains.  
Finally, in Section \ref{S:holo}, we discuss the special case of holomorphic maps from surfaces to Teichm\"uller space.  An a priori bound on the energy of such maps gives rise to the Arakelov-Parsin finiteness result (see Theorem \ref{T:arakelov}).


\subsubsection{Maps to the completion.} \label{S:maps}


As an application of the previous results, we consider the problem of finding energy minimizing equivariant maps to Teichm\"uller space with the Weil-Petersson metric.  Recall the set-up:  let $M$ be a compact Riemannian manifold with universal cover $\widetilde M$, and let $\rho:\Gamma=\pi_1(M)\to \mcg$ be a homomorphism.  Since $\mcg$ acts on $\barTeich$ by isometries, we may ask under what conditions does there exist a $\rho$-equivariant energy minimizing map $f: \widetilde M\to \barTeich$.  

Note that these may be regarded as harmonic maps $M\to \barMod$, although there are two points of caution.  The first is that strictly speaking $\barMod$ is not a manifold, but has orbifold singularities at those points corresponding to Riemann surfaces with automorphisms.  Hence, the smoothness of the map, and the harmonic map equations, should be understood on a smooth finite (local) cover of $\barMod$.  The second (more important) point is that the homotopy class of a map $M\to \Mod$ should be taken in the orbifold sense (i.e.\ equivariantly with respect to a homomorphism $\Gamma\to \mcg$).  Indeed, by the simple connectivity of $\Mod$ remarked on in Section \ref{S:mcg}, homotopy classes of maps to $\Mod$ are very different from equivariant homotopy classes of maps to $\Teich$.

As in Section \ref{S:reductive}, the answer to the existence question depends on the asymptotic dynamics of the image subgroup $\rho(\Gamma)\subset \mcg$. In general, the asymptotic behavior of Weil-Petersson geodesics is quite complicated (see \cite{Brock, BMM}).  As an approximation, one can consider the action on the Thurston boundary $\PMF$ of projective measured foliations.  From this point of view
 there derives a complete classification, analogous to the Thurston classification, of subgroups of the mapping class group. 

\begin{theorem}[McCarthy-Papadopoulos \cite{MP1}] \label{T:mp} \index{mapping class!group!classification of subgroups}
A subgroup of $ \mcg$ is exactly one of the following types:
\begin{enumerate}
\item finite;
\item infinite irreducible and virtually cyclic;
\item infinite reducible;
\item sufficiently large.
\end{enumerate}
\end{theorem}
\noindent 
By \emph{sufficiently large} \index{sufficiently large}  we mean that the subgroup contains two pseudo-Anosov's with distinct fixed point sets in $\PMF$.  These groups contain free groups on two generators.

We apply this theorem to the image  $G=\rho(\Gamma)$ of the homomorphism $\rho$.
By the Nielsen Realization Theorem \ref{T:nielsen}, if  $G$ is finite then it fixes a point $[\sigma]$ in Teichm\"uller space.
Hence, the constant map $f(x)=[\sigma]$ is equivariant and clearly harmonic.  

Case (2) arises when $G$ has a finite index subgroup $\langle[\phi]\rangle\simeq\Z$  generated by a pseudo-Anosov $[\phi]$.  By Theorem \ref{T:dw} this stabilizes a complete Weil-Petersson geodesic  $A_{[\phi]}\subset\Teich$.
The corresponding finite index subgroup $\widehat \Gamma\subset\Gamma$ defines a finite cover $M_{\widehat\Gamma}\to M$, and the group of deck transformations then acts on $S^1$.  Hence, it suffices to find an equivariant harmonic map $M_{\widehat\Gamma}\to S^1\hookrightarrow\Teich/\langle[\phi]\rangle$.  This can be done using the heat equation approach, since  equivariance is preserved under the flow \eqref{E:heat}. 

In Case (3), $G$ fixes a stratum in the boundary $\partial\Teich$ isomorphic to a product of lower dimensional Teichm\"uller spaces.  Since the boundary strata are totally geodesically embedded, the problem of finding an energy minimizer to $\barTeich$ is reduced to Cases (1), (2), and (4) for lower dimensions.

Finally, we come to Case (4).

\begin{theorem}[Daskalopoulos-Wentworth \cite{DW}]   
If $\rho:\Gamma\to \mcg$ is sufficiently large then it is proper in the sense of Korevaar-Schoen (see Section \ref{S:reductive}).
\end{theorem}

This is a consequence of Theorem \ref{T:axes}.
Using Theorem \ref{T:ksproper}, it follows that there exist equivariant harmonic maps in this case as well.  Putting all of these considerations together, we have

\begin{corollary} \label{C:existence}
Let $\rho: \pi_1(M)\to \mcg$ be a homomorphism, where $M$ is a compact Riemannian manifold.  Then there exists a finite energy  $\rho$-equivariant harmonic map $f:\widetilde M\to \barTeich$.  Moreover, $f$ is uniformly Lipschitz.
\end{corollary}

We note that this statement, apparently stronger  than the one appearing in \cite{DW}, follows from considering the possibilities in Theorem \ref{T:mp}.

It is certainly expected that uniqueness holds in the corollary under certain assumptions.  Some generalization of Theorem \ref{T:meseunique} is needed.  Roughly speaking, one expects uniqueness to fail only if the image of $f$ lies in a flat.  Alternatively, one should be able to prove a priori that if $\rho$ is sufficiently large then some point of the image of $f$ lies in the interior $\Teich$.  Then the strictly negative curvature implies that the image is a geodesic, which again contradicts the assumption of sufficiently large.


\subsubsection{Surface domains.}    \label{S:surface}


A natural question arises from the statement of Corollary \ref{C:existence}.  Under what conditions does the image of a harmonic map to $\barTeich$ actually lie in $\Teich$?  This is an important issue, since if $f(x)\in\Teich$, then since $\Teich$ is a manifold $f$ is smooth near the point $x$.  More generally, one would at least like to have control over the size of the singular set (cf.\ Theorem \ref{T:sun}).

The first result in this direction is the following

\begin{theorem}[Wentworth \cite{W1}] \label{T:reg}
Let $\Omega\subset \R^2$ be a bounded domain, and suppose $f:\Omega\to \barTeich$ is energy minimizing with respect to its boundary conditions.  If $f(x)\in \Teich$ for some $x\in\Omega$, then $f(\Omega)\subset \Teich$.
\end{theorem}

To give a rough idea of why this should be the case, we again consider the model for the Weil-Petersson geometry near the boundary $\partial\Teich$ discussed in Section \ref{S:wp}.  
Let $f:B_1(0)\to\barmodel$ be a finite energy harmonic map. By Theorem \ref{T:ks1}, $f$ is uniformly Lipschitz.
The generalization of \eqref{E:modelgeodesic} are the equations
\begin{align*}
\xi\Delta \xi &= \tfrac{3}{4} \xi^6|\nabla\theta|^2 , \\
\Div(\xi^6\nabla \theta)&=0 ,\\
|\nabla \xi|\ ,&  \ \xi^3|\nabla \theta|\ \text{are locally bounded.}
\end{align*}
Because of the singularities, $\xi(x,y)$ and $\theta(x,y)$ are only weak solutions of these equations.
 We may assume that $f$ is nonconstant with  $f(0)\in \partial\model$.  Furthermore, suppose the origin is not a zero of the Hopf differential $\varphi$.  
It is not hard to show that the singular set,\index{singular set} i.e.\  $f^{-1}(\partial\model)$ is a leaf of the horizontal foliation of $\varphi$.  In $\varphi$-coordinates $(x,y)$, one shows with some more analysis of the situation that $\xi(x,y)\sim y$.  Then the second equation above becomes essentially $\Div(y^6\nabla \theta)=0$.
This kind of degenerate equation appears in the study of the porous medium equation \cite{Koch}, and one can show that $\theta(x,y)$ itself is locally bounded.  Then a scaling argument using the monotonicity formula  \eqref{E:monotonicity} can be used to derive a contradiction.

To go from a regularity result for harmonic maps to $\barmodel$ to a result for maps to $\barTeich$ requires an approximation of harmonic maps to targets with asymptotically product metrics.  This is similar to the discussion of geodesics above.  For more details we refer to \cite{W1}.

The local regularity implies 
\begin{corollary} \label{C:reg}
 Let $\rho: \pi_1(B)\to \mcg$ be irreducible, where $B$ is a compact Riemann surface. Then there exists a smooth  $\rho$-equivariant harmonic map $f:\widetilde B\to \Teich$.  Moreover, if $\rho$ is sufficiently large, then $f$  is unique.
\end{corollary}

An interesting potential application of this result pertains to the following

\medskip\noindent {\bf Question.}  Let $B$ be a closed surface.  Does there exist an injective homomorphism $\rho:\pi_1(B)\to\mcg$ such that the image of $\rho$ consists entirely of pseudo-Anosov's ?

\medskip
Examples of all pseudo-Anosov subgroups of $\mcg$ have been constructed in \cite{Wh}, but these are not surface groups.  Such groups, should they exist, would admit minimal surface representations in $\mcg$:

\begin{corollary} \label{C:minimal}
Let $B$ be a closed surface and $\rho:\pi_1(B)\to \mcg$.  In
addition, we assume that for every simple closed essential curve in $B$, the image by $\rho$ of the associated
conjugacy class in $\pi_1(B)$ is pseudo-Anosov.  Then there is 
 a conformal harmonic
$\rho$-equivariant map $f:(\widetilde B, j)\to \Teich$ for some complex structure $j$ on $B$.
\end{corollary}

The argument proceeds as in  the proof of Proposition \ref{P:domain}.  Note that there is a lower bound, depending only on the genus, of the Weil-Petersson translation length of any pseudo-Anosov (see \cite{DW}).


\subsubsection{Holomorphic maps from Riemann surfaces.}  \label{S:holo}


By Proposition \ref{P:holo} (see esp.\ \eqref{E:holo}), since the Weil-Petersson metric is K\"ahler, equivariant holomorphic maps from surfaces to $\Teich$ are examples of energy minimizers; in particular, harmonic maps.   \index{energy!and holomorphic maps} These are given by holomorphic curves in $\Mod$ that are locally liftable to $\Teich$.  Alternatively, consider a family $X\to B$, where $B$ is a compact Riemann surface, and $X$ is a locally liftable holomorphic fibration of genus $p$ Riemann surfaces.  Associated to this is a \emph{monodromy} homomorphism $\rho: \pi_1(B)\to \mcg$.  By Corollary \ref{C:reg}, if $\rho$ is irreducible there is a $\rho$-equivariant harmonic map $f:\widetilde B\to \Teich$.   In general, this will \emph{not} be holomorphic for any choice of complex structure on $B$.  By the essential uniqueness of the harmonic map, we see that the issue of holomorphicity is a property of the (conjugacy class) of the monodromy representation $\rho$. Let us call a homomorphism $\rho:\pi_1(B)\to \mcg$   \emph{holomorphic} if there exists a $\rho$-equivariant holomorphic map $\widetilde B \to \Teich$.

A simple example occurs when the monodromy has finite image.  Then by the Nielsen realization theorem, $\rho$ fixes a point in $\Teich$.  In particular, there is a (constant) holomorphic map.  In terms of the family $X\to B$, this is precisely the case where the lift of the fibration $p^\ast X\to \widehat B$ to some finite cover  $p:\widehat B\to B$ is trivial.  Such a fibration is called \emph{isotrivial}. \index{isotrivial}

The harmonic map point of view provides a tool to study holomorphic families.  Here is one property:

\begin{theorem} \label{T:isotrivial}
 If $\rho:\pi_1(B)\to \mcg$ is holomorphic and nonisotrivial,  then $\rho$ is sufficiently large.
\end{theorem}

\begin{proof}
Suppose not.
By the classification of subgroups of the mapping class group Theorem \ref{T:mp}, $\rho$ is either reducible or virtually cyclic.  In the former case, there is a proper totally geodesic stratum $S\subset \partial \Teich$ that is invariant under $\rho$.  Since projection to $S$ from the interior $\Teich$ is \emph{strictly} distance decreasing, the geodesic homotopy of $f$ to $S$ is both $\rho$-equivariant and strictly energy decreasing.  This contradicts the fact that $u$ is the energy minimizer.    If $\rho$ is virtually cyclic, then the energy minimizer maps onto a geodesic.  Since the image is one dimensional, this contradicts holomorphicity.
\end{proof}

The following is also a consequence of the uniqueness of harmonic maps to $\barTeich$ discussed in the proof above.  This is sometimes called the \emph{rigidity} theorem.

\begin{theorem} \label{T:rigidity}
Holomorphic families with the same monodromy (up to conjugation) are equivalent.
\end{theorem}

The main finiteness result  is the following

\begin{theorem}[Arakelov \cite{Ar}, Parsin \cite{Parsin}] \label{T:arakelov} \index{Arakelov-Parsin theorem}
Fix a closed Riemann surface $B$, and let $\mcg$ denote the mapping class group of a compact surface  of genus $p\geq 2$.  Then there are at most finitely many conjugacy classes of non-isotrivial holomorphic homomorphisms $\rho: \pi_1(B)\to \mcg$.
\end{theorem}

  We note that this can be extended to the case where $B$ is a Riemann surface with punctures.  The punctures correspond to singularities in the surface fibration, and in the holomorphic case the local monodromy around the punctures is pseudoperiodic.  Finite energy maps always exist in this case (see \cite{DKW}).  

The key to Theorem \ref{T:arakelov} is a uniform bound on the energy.  Since $\Teich$ has holomorphic sectional curvature bounded above by a negative constant (see Theorem \ref{T:wpcurvature}), Royden's version of the Yau-Schwartz lemma implies that if $f:\widetilde B\to \Teich$ is holomorphic, then 
$$
f^\ast ds^2_{wp}\leq C ds^2_{\widetilde B}\ ,
$$
for a uniform constant $C$ (see \cite{Ro2}).  In particular, by \eqref{E:holo}, the energy of a holomorphic map is uniformly bounded.  Since by Proposition \ref{P:lipschitz} the Lipschitz constant of harmonic maps is bounded by the total energy, a sequence of holomorphic maps to $\Mod$ is necessarily equicontinuous (see also \cite{GR}).  This allows one to construct convergent subsequences for the maps $\Mod$. 
As in the argument in Section \ref{S:classification} there is the issue of lifting the limiting map.
 In this way,  one derives a contradiction to the existence of infinitely many distinct conjugacy classes of holomorphic $\rho$.  For a fuller account of this approach to the Arakelov-Parsin Theorem, we refer to \cite{JY2} and \cite{IS}.

\subsection{Superrigidity}

In this final section we briefly describe how  equivariant harmonic map theory can be used to study homomorphisms of fundamental groups of compact manifolds to the mapping class group.
The link between superrigidity and harmonic maps uses a technique which can be traced back to Bochner and Calabi-Weil and was first fully utilized in connection with the Margulis superrigidity theorem.  In fact, as mentioned earlier, many of the ideas in this paper were inspired by the attempt to give a harmonic maps proof of superrigidity.  In Section \ref{S:fm}, we state the Ivanov-Farb-Kaimanovich-Masur theorem for homomorphisms of superrigid lattices into mapping class groups.  In Section \ref{S:comb}, we describe two approaches in generalizing harmonic maps by allowing the domain to be singular as well.  the first is the analytic approach along the lines for smooth domains described in this article.  The second is the combinatorial approach.  As an application one can prove a statement on the non-Archimedean superigidity of lattices in mapping class groups.

\subsubsection{The Ivanov-Farb-Kaimanovich-Masur Theorem.}  \label{S:fm} Harvey originally asked whether the mapping class group could be isomorphic to a lattice in a symmetric space \cite{Harvey}.  This was shown not to be the case by Ivanov \cite{Iv1,Iv3}.
For some of the similarities and differences between $\mcg$ and arithmetic lattices, see \cite{FLM,  Iv1, MC} and Ivanov's survey article \cite{Iv2}.   Indeed, a stronger statement is true:

\begin{theorem} \label{T:fm}  \index{Farb-Kaimanovich-Masur theorem} \index{mapping class!group!superrigidity of}
Let $\Gamma$ be a cocompact lattice in any symmetric space with nonpositive curvature other than the real or complex hyperbolic spaces.  Then any homomorphism $\Gamma\to \mcg$ has finite image.
\end{theorem}

For symmetric spaces of rank $\geq 2$ this result is due to Farb-Masur \cite{FM}, following  earlier work of Kaimanovich-Masur \cite{KM}.  Ivanov has announced an independent proof.  Bestvina-Fujiwara \cite{BF} gave a proof using bounded cohomology, and for hermitian symmetric spaces an independent proof can be found in Hain \cite{Hain}.
Using the method of \cite{C2,JY3,MSY} the remaining  rank $1$ cases were proven by S.-K. Yeung \cite{SKY2}.

Geometric superrigidity uses harmonic maps to prove results of this type. 
The basic philosophy is to show that equivariant harmonic maps $f: G/K\to N$, where $G/K$ is a symmetric space of higher rank and $N$ has nonpositive curvature,  would necessarily be totally geodesic.  Recall from Section \ref{S:def}  that the harmonic map equations are of the form $\Tr\nabla df=0$, whereas the equations for a totally geodesic map are $\nabla df=0$.  Curvature conditions must be used to show that the stronger (overdetermined) set of equations are automatically satisfied.  One then attempts to use geometric considerations to rule out the existence of nonconstant totally geodesic maps.

To give a simple example of how this might come about, consider the following

\begin{theorem}[Eells-Sampson \cite{ES}] \label{T:sr}
If $f:\widetilde M\to N$ is an equivariant harmonic map, $N$ is a Riemannian manifold of nonpositive curvature, and $M$ is closed compact with non-negative Ricci curvature, then $f$ is totally geodesic.  If the Ricci curvature of $M$ is positive at one point, $f$ is constant.  If the sectional curvature of $N$ is negative then $f$ is either constant or maps to a geodesic.
\end{theorem}

Indeed,   the statement easily follows by 
 integrating both sides of the Bochner formula \eqref{E:bochner2} and using the divergence theorem.
 When the domain does not satisfy this curvature restriction, the proof fails.  Nevertheless, more sophisticated forms of the Bochner formulas have been derived in the case of domains with Einstein metrics, or more generally, certain parallel tensors.  For more details, we refer to \cite{C2, JY3, MSY}.

 In light of Corollary \ref{C:existence}, one is tempted to prove Theorem \ref{T:fm} using harmonic maps to $\barTeich$.  The difficulty is in the singular nature of the NPC space $\barTeich$.  However,
 the idea that  these techniques could be generalized to singular space targets is one of the major contributions of \cite{GS}.  The argument based on the Bochner formula given above continues to be valid, so long as the singular set\index{singular set} of $f$ is relatively small, e.g.\ has codimension at least $2$, so that the integration by parts needed to apply the divergence theorem holds.  All of this is motivation to extend the regularity result of Theorem \ref{T:reg} to higher dimensional domains.

\subsubsection{Harmonic maps from singular domains.} \label{S:comb}

Thus far we have discussed the theory of harmonic maps from smooth domains into (possibly singular) metric space targets.  These included singular surfaces, $\R$-trees, and the Weil-Petersson completion of Teichm\"uller space.  In this section we sketch two generalizations of this study to the case where the domain is also allowed to be singular.

We start with an analytic approach closely related to the techniques discussed above.  Let $\Sigma$ be a finite $2$-dimensional  simplicial complex.  The restriction to two dimensions is not essential and most of the following results hold in general.  It is important, however, to assume that $\Sigma$ is \emph{admissible}\index{admissible complex} (cf. \cite{JChen, EF}), meaning that it satisfies the following conditions:
\begin{enumerate}
\item  Every simplex is contained in a face (i.e. a $2$-simplex);
\item Every pair of faces can be joined by a sequence of pairwise adjacent faces;
\item $\Sigma$ has no boundary, i.e.\ every edge is contained in at least two faces;
\item $\Sigma$ is flat in that every open face is isometric to an equilateral triangle in $\R^2$.
\end{enumerate}
We also allow ourselves a choice $w$ of weights $w(F)>0$ for each face.  This is an important technical point. Given an NPC space $(X,d)$ and a map $f:\Sigma\to X$,  define the $w$-energy 
$$
E_w(f)=\frac{1}{2}\sum_{F} w(F)\int_F |\nabla f|^2(x)dx\ ,
$$
where the sum is over all faces $F$ of $\Sigma$.  A map $f$ is called \emph{$w$-harmonic} if it is locally energy minimizing among all maps of finite $w$-energy.\index{harmonic map!singular domain}
As before, we also consider the equivariant theory, where $f$ is a map from the universal cover $\widetilde \Sigma$  of $\Sigma$ that is equivariant with respect to a homomorphism $\rho:\Gamma=\pi_1(\Sigma)\to \iso(X)$.  The existence Theorem \ref{T:ks2} then holds for domains $\Sigma$ as well (cf.\ \cite{EF, DM1}).  

In the following, we will assume a fixed choice of weights and omit $w$ from the notation.  Perhaps the most interesting feature of harmonic maps from simplicial domains is the H\"older continuity.  This was first proven by J. Chen \cite{JChen} for flat metrics  and in a more general context by Eells-Fuglende \cite{EF}.  The following stronger version describes the singular behavior near the vertices.

\begin{theorem}[Daskalopoulos-Mese \cite{DM1}]  \label{T:dm1}
Let $f:\Sigma\to X$ be harmonic.  Then for  domains $U\subset\subset \Omega\subset \Sigma$,
\begin{enumerate}
\item  $f$ is Lipschitz continuous on $U$ away from the vertices of $\Sigma$, with the Lipschitz constant depends only on $U$,  the total energy on $\Omega$, and the distance to the vertex set;
\item  Let $v$ be a vertex with $\alpha=\ord_v(f)$, where the order is defined as in \eqref{E:order}.  Then there exists $r_0>0$ and $C$ depending only on the energy of $f$ such that
$$
\sup_{x\in B_r(v)}|\nabla f|^2(x)\leq C r^{2\alpha-2}
$$
for all $0<r\leq r_0$.
\end{enumerate}
\end{theorem}

The important point here is that, unlike the case of smooth domains, $\alpha$ need not be $\geq 1$.
One  application of Theorem \ref{T:dm1} is the compactification of character varieties for arbitrarily finitely presented groups along the lines of Theorem \ref{T:ddw}. Indeed, one can always realize such a group as the fundamental group of an admissible $2$-complex.
 Other potential consequences use the notion of a Hopf differential.  Clearly, for  energy minimizers, $\varphi=\Hopf(f)$ is a holomorphic quadratic differential on the interior of each face (cf.\ Section \ref{S:skora}).  For points $x$ on  an edge $e$, we have the following \emph{balancing condition}\index{balancing condition}:
$$
\imag\sum_{F} \varphi_F(x)=0\ ,
$$
where the sum is over all faces $F$ adjacent to $e$ at $x$.  An important open question is whether zeros of $\varphi$ can accumulate along the edges.  If not, then the Hopf differentials of $w$-harmonic maps  define \emph{geometric} or \emph{track foliations} on $\Sigma$ (cf.\ \cite{Bow, Dun, LP}).  Another important issue is the asymptotic behavior of the induced foliation on $\widetilde \Sigma$.  More generally, one might ask under what conditions one can generalize to this setting the results  for surface groups discussed previously in this paper.  

We now return to the relationship between regularity and rigidity.  We have the following
\begin{theorem}[Daskalopoulos-Mese \cite{DM1}]  \label{T:dm2}
Let $f:\Omega\subset\Sigma\to X$ be energy minimizing, where $X$ is a smooth manifold of nonpositive curvature.  For any $x\in \Omega$ which is not a vertex, then there is a neighborhood $U$ of $x$ such that for any face $F$ the restriction of $f$ to $\overline F\cap U$ is smooth.
\end{theorem}

Using this, one has a nontrivial generalization of Theorem \ref{T:sr} to the case of singular domains:
\begin{theorem}[Daskalopoulos-Mese \cite{DM2}] \label{T:dm3}
Suppose $\Sigma$ is an admissible $2$-simplex and $X$ is a complete Riemannian manifold of nonpositive curvature.  If $f:\Sigma\to X$ is harmonic and $|\nabla f|^2$ bounded, then $f$ is totally geodesic on each simplex of $X$.  If the sectional curvature of $X$ is strictly negative, then either $f$ is constant or it maps to a geodesic.
\end{theorem}

This is a kind of rigidity result for the group $\Gamma=\pi_1(\Sigma)$, and a combinatorial version was first proven by M.-T. Wang (see below).  See also \cite{BS}.  The result follows by the Bochner formula \eqref{E:bochner2}, the vanishing of the Ricci curvature on the domain,  and the fact that $|\nabla f|^2$ allows us to integrate by parts.  Global boundedness of the energy density is guaranteed by a combinatorial condition on $\Sigma$.  Namely, the first eigenvalue of the discrete Laplacian on the link of every vertex with the induced weights should be  $\geq 1/2$ (see \cite{DM2}).    This condition is a generalization of the notion of p-adic curvature that first appeared in the work of Garland (cf.\ \cite{Garland}).

The second approach discretizes the notion of an energy minimizer.  Let $\Sigma$ be an admissible $2$-complex and $X$ an NPC space as above and $\rho:\Gamma\to \iso(X)$ a homomorphism.  Given a system of weights on the faces of $\Sigma$ there is a standard way to induce weights on the lower dimensional simplices.  For example, the weight of an edge is the sum of the weights of adjacent faces.  Let $\Sigma^i$, $\widetilde \Sigma^i$ denote the $i$-skeletons.  Given a $\rho$-equivariant map $f:\widetilde\Sigma^0\to X$ define its energy by 
$$
E_{comb.}(f)=\frac{1}{2}\sum_{e_{xy}\in \Sigma^1}w(e_{xy}) d^2(f(\tilde x), f(\tilde y))\ ,
$$
where $e_{xy}$ denotes an edge with adjacent vertices $x$ and $y$, and $\tilde x$, $\tilde y$ are adjacent vertices of a lift of $e_{xy}$ to $\widetilde \Sigma$.  We say that $f$ is a \emph{$\rho$-equivariant combinatorial harmonic map}  if it minimizes $E_{comb.}(f)$.\index{harmonic map!combinatorial} \index{harmonic map!combinatorial}  Under the assumption that $X$ is locally compact and that $\rho(\Gamma)$ does not fix a point in $\partial X$ one can prove the existence of combinatorial harmonic maps (see Wang \cite{Wang1,Wang2}).  Furthermore, assuming the first eigenvalue of the combinatorial Laplacian of the link of every vertex with the induced weights is $>1/2$, one can deduce rigidity results as in the first approach.  
This  can be used to deduce   non-Archimedean generalizations of Theorem \ref{T:fm} \cite{SWang}.
%

\printindex

\end{document}